\documentclass{article}
\usepackage{fullpage}
\usepackage{fancyhdr}
\usepackage{amsmath, amssymb, latexsym, amscd, amsthm,amsfonts,amstext}
\usepackage{subfigure}
\usepackage{enumerate}
\usepackage{xcolor}
\usepackage{textcomp}
\usepackage{mathtools}

\usepackage{graphicx}
\newtheorem{theorem}{Theorem}[section]
\newtheorem{lemma}[theorem]{Lemma}
\newtheorem{proposition}[theorem]{Proposition}

\numberwithin{equation}{section}
\newtheorem{definition}[theorem]{Definition}
\newtheorem{remark}[theorem]{Remark}

\newenvironment{indeed}[0]{Indeed, }{}
\numberwithin{equation}{section}

\usepackage[mathscr]{eucal}
\usepackage{graphicx}
\usepackage{subfigure}
\usepackage{cite}
\usepackage{stfloats}
\usepackage{fancyhdr}
\usepackage{enumerate}
\usepackage{xcolor}

\allowdisplaybreaks

\allowdisplaybreaks
\title{On the justification of the Foldy-Lax approximation for the acoustic scattering by small rigid bodies of arbitrary shapes} 

\author{
Durga Prasad Challa
\thanks{RICAM, Austrian Academy of Sciences,
Altenbergerstrasse 69, A-4040, Linz, Austria.
(Email: durga.challa@oeaw.ac.at)
\newline
Supported by the Austrian Science Fund (FWF): P22341-N18.} 
\and  Mourad Sini
\thanks{RICAM, Austrian Academy of Sciences,
Altenbergerstrasse 69, A-4040, Linz, Austria.
(Email:mourad.sini@oeaw.ac.at)
\newline
Partially supported by the Austrian Science Fund (FWF): P22341-N18.}
}

\begin{document}
\maketitle
\begin{abstract}
  We are concerned with the acoustic scattering problem by many small rigid obstacles of arbitrary shapes. 
 We give a sufficient condition on the number $M$ and the diameter $a$ of the obstacles as well as 
 the minimum distance $d$ between them under which the Foldy-Lax approximation is valid. 
 Precisely, if we use single layer potentials for the representation of the scattered fields, as it is done sometimes in the literature, then this condition is
 $(M-1)\frac{a}{d^2} <c$, with an appropriate constant $c$, while if we use double layer potentials then a weaker condition of the form
 $\sqrt{M-1}\frac{a}{d} <c$ is enough. 
 In addition, we derive the error in this approximation explicitly in terms of the parameters $M, a$ and $d$.
The analysis is based, in particular, on the precise scalings of the boundary integral operators between the corresponding
 Sobolev spaces. As an application, we study the inverse scattering
by the small obstacles in the presence of multiple scattering. 
\end{abstract}




\section{Introduction and statement of the results}\label{Introduction-smallac-sdlp}

 Let $B_1, B_2,\dots, B_M$ be $M$ open, bounded and simply connected sets in $\mathbb{R}^3$ with Lipschitz boundaries
\footnote{Let us recall that the surface $\partial B_j$ is of Lipschitz class with
the Lipschitz constants $r_j, L_j > 0$ if for any $P \in \partial B_j$, there exists a rigid transformation of coordinates
under which we have $P = 0$ and $B_j \cap  B^3_{r_j}(0) = \{x \in B_{r_j}(0) : x_3 > \varphi_j(x_1, x_2
)\}$, $B^3_{r_j}(0)$ being the ball of center $0$ and radius $r_j$ ,
where $\varphi_j$ is a Lipschitz continuous function on the disc of center $0$ and radius $r_j$, i.e. $B^2_{r_j}(0)$, satisfying $\varphi_j(0) = 0$ and
$\Vert \varphi_j \Vert_{ 
C^0_1(B^2_{r_j}(0))}
\leq L_j.$} containing the origin.
We assume that the Lipschitz constants of $B_j$, $j=1,..., M$ are uniformly bounded.  
We set $D_m:=\epsilon B_m+z_m$ to be the small bodies characterized by the parameter 
$\epsilon>0$ and the locations $z_m\in \mathbb{R}^3$, $m=1,\dots,M$. Let $U^{i}$ be a solution of the Helmholtz equation $(\Delta + \kappa^{2})U^{i}=0 \mbox{ in } \mathbb{R}^{3}$.  
We denote by  $U^{s}$ the acoustic field scattered by the $M$ small bodies $D_m\subset \mathbb{R}^{3}$ due to 
the incident field $U^{i}$. We restrict ourselves to (1.) the plane incident waves, $U^{i}(x,\theta):=e^{ikx\cdot\theta}$, 
with the incident direction $\theta \in \mathbb{S}^2$, with $\mathbb{S}^2$ being the unit sphere, and (2.) the scattering by rigid bodies. Hence the total field $U^{t}:=U^{i}+U^{s}$ satisfies the following exterior Dirichlet problem of the acoustic waves
\begin{equation}
(\Delta + \kappa^{2})U^{t}=0 \mbox{ in }\mathbb{R}^{3}\backslash \left(\mathop{\cup}_{m=1}^M \bar{D}_m\right),\label{acimpoenetrable}
\end{equation}
\begin{equation}
U^{t}|_{\partial D_m}=0,\, 1\leq m \leq M, \label{acgoverningsupport}  
\end{equation}
\begin{equation}
\frac{\partial U^{s}}{\partial |x|}-i\kappa U^{s}=o\left(\frac{1}{|x|}\right), |x|\rightarrow\infty, ~(\text{S.R.C}) \label{radiationc}
\end{equation}
where  $\kappa>0$ is the wave number, $\kappa=2\pi\slash \lambda$, $\lambda$ is the wave length and S.R.C stands for the Sommerfield radiation condition.  The scattering problem 
(\ref{acimpoenetrable}-\ref{radiationc}) is well posed in the H\"{o}lder or Sobolev spaces, see \cite{C-K:1983, C-K:1998, Mclean:2000} for instance, and the scattered field $U^s(x, \theta)$ has the following asymptotic expansion:
\begin{equation}\label{far-field}
 U^s(x, \theta)=\frac{e^{i \kappa |x|}}{|x|}U^{\infty}(\hat{x}, \theta) + O(|x|^{-2}), \quad |x|
\rightarrow \infty,
\end{equation}
with $\hat{x}:=\frac{x}{\vert x\vert}$, where the function
$U^{\infty}(\hat{x}, \theta)$ for $(\hat{x}, \theta)\in \mathbb{S}^{2} \times\mathbb{S}^{2}$  is called the far-field pattern.
We recall that the fundamental solution, $\Phi_\kappa(x,y)$, of the Helmholtz equation in $\mathbb{R}^3$
with the fixed wave number $\kappa$ is given by
\begin{eqnarray}\label{definition-ac-small-fundamentalkappa}
 \Phi_\kappa(x,y)&:=&\frac{e^{i\kappa|x-y|}}{4\pi|x-y|},\quad \text{for all } x,y\in\mathbb{R}^3.
\end{eqnarray}
\begin{definition} 
\label{Def1}
We define 
\begin{enumerate}
 \item $a$ as the maximum among the diameters, $diam$, of the small bodies $D_m$, i.e.
\begin{equation}\label{def-a} 
a:=\max\limits_{1\leq m\leq M } diam (D_m) ~~\big[=\epsilon \max\limits_{1\leq m\leq M } diam (B_m)\big],
\end{equation}

 \item  $d$ as the minimum distance  between the small bodies $\{D_1,D_2,\dots,D_m\}$, i.e.
\begin{equation}\label{def-d}
d:=\min\limits_{\substack{m\neq j\\1\leq m,j\leq M }} d_{mj},
\end{equation}
$\text{where}\,d_{mj}:=dist(D_m, D_j)$. We assume that
\begin{equation}\label{def-dmax}
0\,<\,d\,\leq\,d_{\max},
\end{equation}
and $d_{\max}$ is given.
\item $\kappa_{\max}$ as the upper bound of the used wave numbers, i.e. $\kappa\in[0,\,\kappa_{\max}]$.
\end{enumerate} 
\end{definition}
\bigskip

The scattering by a collection of point-like scatterers $z_m$, $m=1,\dots,M,$ has been studied since the works by L. L. Foldy \cite{LLF:PR1945} and M. Lax \cite{Lax-M:RMP1951}, 
see also \cite{A-G-H-H:AMS2005 , Martin:2006} for a detailed account on this problem and \cite{DPC-HG-SM:M3AS2014,DPC-SM:InvPrbm2012,H-S} for the corresponding results for the Maxwell and the Lam\'e models respectively.
In their modelling, the total field is of the form:
\begin{equation}\label{Foldy-Phy-acintro-pointlike}
 U^{t}(x) =  U^i(x)+\sum_{m=1}^M\,\Phi_{\kappa}(x,z_m)\,A_m,
\end{equation}
where $A_m$ are solutions of the (so-called Foldy-Lax) linear algebraic system
\begin{equation}\label{Foldy-Phy-ext1-acintro-pointlike}
A_m-\sum_{\substack{j=1 \\ j\neq m}}^M\,g_m\Phi_{\kappa}(z_m,z_j)\,A_j = g_mU^i(z_m),
\end{equation}
where $g_m$, $m=1,\dots,M$, are the given scattering coefficients. From \eqref{Foldy-Phy-acintro-pointlike}, we derive the following representation of the far-field
\begin{equation}\label{Foldy-Phy-acintro-pointlike-farfield}
U^\infty(\hat{x},\theta) =  \sum_{m=1}^M\,e^{-i\kappa\hat{x}\cdot z_m}\,A_m.
\end{equation}

\par 
The goal of our work is to dervie an asymptotic expansion of the scattered field by the collection of the small scatterers $D_m,m=1,\dots,M$, taking into account the parameters $a, d$ and $M$. 
This is the object of the following theorem.
\bigskip

\begin{theorem}\label{Maintheorem-ac-small-sing}
 There exist two positive constants $a_0$ and $c_0$ depending only on the 
Lipschitz character of $B_m,m=1,\dots,M$, $d_{\max}$ and $\kappa_{\max}$ such that
if 
\begin{equation} \label{conditions}
a \leq a_0 ~~ \mbox{and} ~~ \sqrt{M-1}\frac{a}{d}\leq c_0
\end{equation} 
then the far-field pattern $U^\infty(\hat{x},\theta)$ has the following asymptotic expansion
\begin{equation}\label{x oustdie1 D_m farmain}
U^\infty(\hat{x},\theta)=\sum_{m=1}^{M}e^{-i\kappa\hat{x}\cdot z_m}Q_m+
O\left(M\kappa a^2+M(M-1)\left(\frac{\kappa a^3}{d}+\frac{a^3}{d^2}\right)+
M(M-1)^2\frac{a}{d}\left(\frac{\kappa a^3}{d}+\frac{a^3}{d^2}\right)\right),
 \end{equation}
uniformly in $\hat{x}$ and $\theta$ in $\mathbb{S}^2$. The constant appearing in the estimate $O(.)$ depends only on 
the Lipschitz character of the obstacles, $a_0$, $c_0$, $d_{\max}$ and $\kappa_{max}$. The coefficients $Q_m$, $m=1,..., M,$ are the solutions of the following linear algebraic system
\begin{eqnarray}\label{fracqcfracmain}
 Q_m +\sum_{\substack{j=1 \\ j\neq m}}^{M}C_m \Phi_\kappa(z_m,z_j)Q_j&=&-C_mU^{i}(z_m, \theta),~~
\end{eqnarray}
for $ m=1,..., M,$ with $C_m:=\int_{\partial D_m}\sigma_m(s)ds$ and $\sigma_{m}$ is 
the solution of the integral equation of the first kind
\begin{eqnarray}\label{barqcimsurfacefrm1main}
\int_{\partial D_m}\frac{\sigma_{m} (s)}{4\pi|t-s|}ds&=&1,~ t\in \partial D_m.
\end{eqnarray}
The algebraic system \eqref{fracqcfracmain} is invertible under the conditions:
\footnote{If $\Omega$ is a domain containing the small bodies, and $diam(\Omega)$ denotes its diameter, 
then one example for the validity of the second condition in \eqref{invertibilityconditionsmainthm} is $diam(\Omega)< \frac{\pi}{2\kappa}$.}
\begin{eqnarray}\label{invertibilityconditionsmainthm}
 \frac{a}{d}\leq c_1 \text{ and } 
 \min_{j\neq m}\cos(\kappa \vert z_j- z_m\vert)\geq 0,
\end{eqnarray}
where $c_1$ depends only on Lipschitz character of the obstacles $B_j$, $j=1, ..., M$. 

\end{theorem}
\bigskip

The asymptotic expression in \eqref{x oustdie1 D_m farmain} provides us with an approximation of the far-fields by the dominant term $\sum_{m=1}^{M}e^{-i\kappa\hat{x}\cdot z_m}Q_m$. This term 
is reminiscent to the right hand side of \eqref{Foldy-Phy-acintro-pointlike-farfield} replacing $A_m$ by $Q_m$ and taking the scattering coefficient $g_m$ as the capacitance $C_m$ of the small scatterer $D_m$ with negative sign. 
This suggests naturally to call \eqref{x oustdie1 D_m farmain} as the Foldy-Lax approximation for the scattering by the collection of the small obstacles $D_m$.


Before discussing this result compared to the existing literature, let us first mention the following remark.
\begin{remark}\label{remarkof-Maintheorem-ac-small-sing}
  The second condition of \eqref{conditions} can be replaced by the stronger one:
\begin{equation} \label{conditions1}
 (M-1)\frac{a}{d}\leq c_2,
\end{equation} 
with $c_2$ depending only on the 
Lipschitz character of $B_m,m=1,\dots,M$ and $\kappa_{\max}$, under which \eqref{x oustdie1 D_m farmain} is reduced to: 
\begin{equation}\label{x oustdie1 D_m farmain1}
U^\infty(\hat{x},\theta)=\sum_{m=1}^{M}e^{-i\kappa\hat{x}\cdot z_m}Q_m+O\left(M\kappa a^2+M(M-1)\left(\frac{\kappa a^3}{d}+\frac{a^3}{d^2}\right)\right),
 \end{equation}
and the algebraic system \eqref{fracqcfracmain} is invertible as well. Remark that in this case the condition \eqref{def-dmax} is not required.

Due to the condition \eqref{conditions1}, the approximation \eqref{x oustdie1 D_m farmain1} can also be reduced to 
\begin{equation}\label{Ramm-model}
U^\infty(\hat{x},\theta)=\sum_{m=1}^{M}e^{-i\kappa\hat{x}\cdot z_m}Q_m+O\left(M\left(\kappa+\frac{1}{d}\right)a^2\right).
\end{equation}
However this models only the Born approximation while the approximation in \eqref{x oustdie1 D_m farmain1} takes into account some multiple scattering as the first order interaction, 
see Section \ref{The inverse problem} for more details.
\end{remark}
\bigskip

\par This type of result, precisely the estimate (\ref{Ramm-model}), is known by A. Ramm since mid 1980's, 
see \cite{MRAMAG-book1, RAMM:2007, RAMM:2011} and the references therein for his recent related results. However, he used the (rough) condition 
$\frac{a}{d}\ll 1$ and no mention has been made  on the number of obstacles $M$. In his arguments, he used the single layer potential representation (SLPR) of the scattered field. Using the same 
representation, we realize that a condition of the type $(M-1)\frac{a}{d^2}\leq c$ is needed, see Proposition \ref{normofsigmastmt}.  However, using (naturally)
 the double layer potential representation (DLPR) we need only the weaker condition $\sqrt{M-1}\frac{a}{d}\leq c_0$, see Section \ref{DLPR}.
This condition appears naturally since we use the Neumann series expansion to estimate the inverse of the boundary operator, 
see the proof of Proposition \ref{normofsigmastmtdbl}.

The particular but important case where the obstacles have circular shapes has been considered recently by M. Cassier and C. Hazard in \cite{C-H:2013} where error estimates are obtained 
replacing the condition $\sqrt{M-1}\frac{a}{d} \leq c_0$ by the weaker one of the form $\frac{a}{d} \leq c_0$ (appearing implicitly 
in their analysis). 
This is possible due to the use of the Fourier series expansion of the scattered fields with which they could avoid the use of the Neumann series expansion 
to estimate the inverse of the corresponding boundary operator. However, as it is mentioned in \cite{C-H:2013}, they did not provide 
quantitative estimate of the errors in terms of the density of the obstacles (i.e. $M$ and $d$).

Let us also mention the approach by V. Maz'ya and A. Movchan \cite{M-M:MathNach2010} and  by V. Maz'ya, A. Movchan and M. Nieves \cite{M-M-N:MMS:2011} where 
asymptotic methods are used to study
boundary value problems for the Laplacian with source terms in bounded domains. 
They obtain estimates in forms similar to the previous theorem with weaker conditions of the form $\frac{a}{d} \leq c$, or 
$\frac{a}{d^2} \leq c$, (where, here and in \cite{C-H:2013}, $d$ is the smallest distance between the centers of the scatterers).
In their analysis, they rely on the maximum principle to treat the boundary estimates in addition to the fact that the source terms 
are assumed to be supported away from the small obstacles, a condition that can not be satisfied for scattering by incident plane waves.
To avoid the use of the maximum principle, which is not valid due to the presence of the wave number $\kappa$, 
we use boundary integral equation methods. The price to pay is the need of the stronger assumption $\sqrt{M-1}\frac{a}{d}\leq c_0$.

The integral equation methods are widely used in such a context, see for instance the series of works by H. Ammari and H. Kang 
 and their collaborators, as \cite{A-K:2007} and the references therein. They combine layer potential techniques with the series expansion of the Green's 
functions of the background medium to derive the full asymptotic expansion in terms of the polarization tensors. The difference between their
 asymptotic expansion
and the one described in the previous theorem is that their polarization tensors are built up from densities which are solutions of a system of 
integral equations while in the previous theorem the approximating terms are built up from the linear algebraic system (\ref{fracqcfracmain}). 
Related to the results stated in Theorem \ref{Maintheorem-ac-small-sing}, we cite the following works \cite{A-H-K-L, AH-KH-LH:AMS2009, A-K:JMPA2003, AH-KH:JMAA2004,AH-GP-GBL-KH:JDE2010} 
where different models are also studied. The new part in Theorem \ref{Maintheorem-ac-small-sing} is the precise
estimate of the reminder in the asymptotic expansion in terms of $M$, $a$ and $d$.

Let us finally mention that the linear algebraic system (\ref{fracqcfracmain}) appeared also in the works \cite{MRAMAG-book1, M-M:MathNach2010, M-M-N:MMS:2011} we cited above
and in the work \cite{GL-RV-JCP1987,GL-RV-JCP1997} by L. Greengard and V. Rokhlin where the fast multipole method is designed to solve it numerically, see also \cite{BoardJr-JCP1997}.  

Before concluding the introduction, we find it worth mentioning the following remark. 
\begin{remark}\label{corMaintheorem-ac-small-sing}
 If we assume, in addition to the conditions of Theorem \ref{Maintheorem-ac-small-sing}, that $D_m$ are balls with same diameter $a$ for $m=1,\dots,M$, then we have the following asymptotic expansion of the far-field pattern:
  \begin{eqnarray}\label{x oustdie1 D_m farmain-recent-near}
 \hspace{-1cm}U^\infty(\hat{x},\theta)
&=&\sum_{m=1}^{M}e^{-i\kappa\hat{x}\cdot z_m}Q_m\nonumber\\
&&+O\left(M\left[a^2+\frac{a^3}{d^{5-3\alpha}}+\frac{ a^4}{d^{9-6\alpha}}\right]+M(M-1)\left[\frac{a^3}{d^{2\alpha}}+\frac{a^4}{d^{4-\alpha}}+\frac{a^4}{d^{5-2\alpha}}\right]
+M(M-1)^2\frac{a^4}{d^{3\alpha}}\right)
 \end{eqnarray}
 where $0<\alpha\leq1$.
\par
Consider now the special case $d=a^t,\,M=a^{-s}$ with  $t,s>0$. Then the asymptotic expansion \eqref{x oustdie1 D_m farmain-recent-near} can be rewritten as 
\begin{eqnarray*}\label{x oustdie1 D_m farmain-recent**}
U^\infty(\hat{x},\theta)
\hspace{-.05cm}&=&\hspace{-.1cm}\sum_{m=1}^{M}e^{-i\kappa\hat{x}\cdot z_m}Q_m\hspace{-.03cm}+\hspace{-.03cm}O\left(a^{2-s}\hspace{-.03cm}+\hspace{-.03cm}a^{3-s-5t+3t\alpha}\hspace{-.03cm}+\hspace{-.03cm}a^{4-s-9t+6t\alpha}\hspace{-.03cm}+\hspace{-.03cm}a^{3-2s-2t\alpha}\hspace{-.03cm}+\hspace{-.03cm}a^{4-3s-3t\alpha}\hspace{-.03cm}+\hspace{-.03cm}a^{4-2s-5t+2t\alpha}\right).
 \end{eqnarray*}
As the diameter $a$ tends to zero the error term tends to zero for $t$ and $s$ such that 
$0<t<1$ and $0<s<\min\{2(1-t),\,\frac{7-5t}{4},\,\frac{12-9t}{7},\frac{20-15t}{12},\frac{4}{3}-t\alpha\}$. 
In particular for $t=\frac{1}{3}$, $s=1$, we have 
\begin{eqnarray}\label{x oustdie1 D_m farmain-recent***}
U^\infty(\hat{x},\theta)
&=&\sum_{m=1}^{M}e^{-i\kappa\hat{x}\cdot z_m}Q_m+O\left(a+a^{2\alpha}+a^{1-\alpha}+a^{\frac{1+2\alpha}{3}}\right)\nonumber\\
&=&\sum_{m=1}^{M}e^{-i\kappa\hat{x}\cdot z_m}Q_m+O\left(a^{\frac{1}{2}}\right)\qquad[\mbox{obtained for } \alpha=\frac{1}{4}].
 \end{eqnarray}
 This particular case is used to derive the effective medium by perforation using many small bodies, see \cite{RAMM:2007, RAMM:2011} for instance.  The result of Remark 
 \ref{corMaintheorem-ac-small-sing}, in particular \eqref{x oustdie1 D_m farmain-recent***}, ensures the rate of the error in deriving such an effective medium.
\par
 Actually, the results \eqref{x oustdie1 D_m farmain-recent-near} and \eqref{x oustdie1 D_m farmain-recent***} are valid 
  for the non-flat Lipschitz obstacles $D_m=\epsilon{B}_m+z_m, m=1,\dots,M$ with the same diameter $a$, i.e. $D_m$'s are Lipschitz obstacles and there exist constants $t_m \in (0, 1]$ such that
 \begin{equation}\label{non-flat-condition}
 B^{3}_{t_m\frac{a}{2}}(z_m)\subset\,D_m\subset\,B^{3}_{\frac{a}{2}}(z_m),
 \end{equation}
 where $t_m $ are assumed to be uniformly bounded from below by a positive constant. 
 
 Finally, we observe that the distribution of the scatterers such that $M=a^{-1}$ and $d=a^\frac{1}{3}$, as $a \rightarrow 0$, used in the effective medium theory, 
 fits into the condition, $\sqrt{M-1}\frac{a}{d} \leq c_0$, we needed in the double layer representation,  
 but not into the condition, $(M-1)\frac{a}{d^2} \leq c$, we needed in the single layer representation.
 \end{remark}
 \bigskip
 
The rest of the paper is organized as follows. In Section \ref{Proof of Theorem Small}, we prove Theorem \ref{Maintheorem-ac-small-sing} in two steps. In Section \ref{SLPR}, we use single
layer potential representation of the scattered field and show how the stronger condition $(M\hspace{-.07cm}-\hspace{-.07cm}1)\frac{a}{d^2}\leq c$ appears naturally. 
In Section \ref{DLPR}, we use double layer potentials and reduce this condition to the weaker one $\sqrt{M\hspace{-.07cm}-\hspace{-.07cm}1}\frac{a}{d}\leq c_0$. In Section  \ref{The inverse problem}, we deal with 
the inverse scattering problem which consists of recovering the location and the size of the small scatterers from far-field patterns corresponding
to finitely many incident plane waves. We present numerical tests using the MUSIC algorithm with an emphasis on the multiple scattering effect 
due to closely spaced scatterers.

\section{Proof of Theorem \ref{Maintheorem-ac-small-sing}}\label{Proof of Theorem Small}
We wish to kindly warn the reader that in our analysis we use sometimes the parameter $\epsilon$ 
and some other times the parameter $a$ as they appear naturally in the estimates.
But we bear in mind the relation (\ref{def-a}) between $a$ and $\epsilon$. 
\subsection{Using SLP representation}\label{SLPR}
\subsubsection{The representation via SLP}\label{SLPR-1}
We start with the following proposition on the solution of the problem (\ref{acimpoenetrable}-\ref{radiationc}) via the single layer potential representation.
\begin{proposition}\label{existence-of-sigmas}
For $m=1,2,\dots,M$, there exists $\sigma_m\in L^2(\partial D_m)$ such that the problem (\ref{acimpoenetrable}-\ref{radiationc}) has a unique solution of the form
 \begin{equation}\label{qcimprequiredfrm1}
  U^{t}(x)=U^{i}(x)+\sum_{m=1}^{M}\int_{\partial D_m}\Phi_\kappa(x,s)\sigma_{m} (s)ds,~x\in\mathbb{R}^{3}\backslash\left(\mathop{\cup}_{m=1}^M \bar{D}_m\right), 
\end{equation}
\end{proposition}
\begin{proofn}{\it{of Proposition \ref{existence-of-sigmas}}.}
  We look for the solution of the problem (\ref{acimpoenetrable}-\ref{radiationc}) of the form \eqref{qcimprequiredfrm1}, 
  then from the Dirichlet boundary condition \eqref{acgoverningsupport}, we obtain
\begin{equation}\label{qcimprequiredfrm1bd}
\sum_{m=1}^{M}\int_{\partial D_m}\Phi_\kappa(s_j,s)\sigma_{m} (s)ds=-U^i{(s_j)},\,\forall s_j\in \partial D_j,\, j=1,\dots,M.
\end{equation}
One can write it in a compact form as 
$(L+K)\sigma=-U^I$ with $L:=(L_{mj})_{m,j=1}^{M}$ and $K:=(K_{mj})_{m,j=1}^{M}$, where
\begin{eqnarray}\label{definition-L_K}
L_{mj}=\left\{\begin{array}{ccc}
            \mathcal{S}_{mj} & m=j\\
            0 & else
           \end{array}\right.,
&  \ &K_{mj}=\left\{\begin{array}{ccc}
            \mathcal{S}_{mj} & m\neq j\\
            0 & else
           \end{array}\right., 
\end{eqnarray}

$U^I=U^I(s_1,\dots,s_M):=\left(U^i(s_1),\dots,U^i(s_M)\right)^T$ and $\sigma=\sigma(s_1,\dots,s_M):=\left(\sigma_1(s_1),\dots,\sigma_M(s_M)\right)^T$. 
Here, for the indices  $m$ and $j$ fixed, $\mathcal{S}_{mj}$ is the integral operator acting as

\begin{eqnarray}\label{defofSmjed}
 \mathcal{S}_{mj}(\sigma_j)(t):=\int_{\partial D_j}\Phi_\kappa(t,s)\sigma_j(s)ds,\quad t\in\partial D_m. 
\end{eqnarray}
Then the operator $\mathcal{S}_{mm}:H^{-s}(\partial D_m)\rightarrow H^{1-s}(\partial D_m)$ is isomorphic and hence
 Fredholm with zero index and for $m\neq j$, $\mathcal{S}_{mj}:H^{-s}(\partial D_j)\rightarrow H^{1-s}(\partial D_m)$ 
is compact for $0\leq s\leq 1$, 
see \cite[Theorem 4.1]{MD:InteqnsaOpethe1997}.\footnote{This property is proved for the case $\kappa=0$. 
By a perturbation argument, we can obtain the same results for every $\kappa$ such that $\kappa^2$ is not a Dirichlet-eigenvalue of $-\Delta$ in $D_m$. 
This last condition is satisfied for every $\kappa$ fixed, $\kappa\leq\kappa_{\max}$, if we take $a<\frac{1}{\kappa_{\max}}\sqrt[3]{\frac{4\pi}{3}}{\rm j}_{1/2,1}$.
Here ${\rm j}_{1/2,1}$ is the 1st positive zero of the Bessel function ${\rm J}_{1/2}$.
}  In our case  we consider $s=0$. Remark here that, for the scattering by single obstacle $K$ is zero  operator.
So, $(L+K):\prod\limits_{m=1}^{M}L^{2}(\partial D_m)\rightarrow \prod\limits_{m=1}^{M}H^{1}(\partial D_m)$ is Fredholm with zero index.  
We induce the product of spaces by the maximum of the norms of the spaces.
To show that $(L+K)$ is invertible it is enough to show that it is injective, i.e. $(L+K)\sigma=0$ implies $\sigma=0$.
\newline
We write 
$$\tilde{U}(x)=\sum_{m=1}^{M}\int_{\partial D_m}\Phi_\kappa(x,s)\sigma_m(s)ds, \mbox{ in } \mathbb{R}^3\backslash 
\left(\mathop{\cup}_{m=1}^{M}\bar{D}_m\right)$$
and 
$$\tilde{\tilde{U}}(x)=\sum_{m=1}^{M}\int_{\partial D_m}\Phi_\kappa(x,s)\sigma_m(s)ds, \mbox{ in } \mathop{\cup}_{m=1}^{M}D_m.$$

Then $\tilde{U}$ satisfies $\Delta\tilde{U}+\kappa^2\tilde{U}=0$ for $x\in\mathbb{R}^{3}\backslash\left(\mathop{\cup}\limits_{m=1}^{M}\bar{D}_m\right)$,  
with S.R.C and $\tilde{U}(x)=0$ on $\mathop{\cup}\limits_{m=1}^{M}\partial D_m$.  
Similarly, $\tilde{\tilde{U}}$ satisfies $\Delta\tilde{\tilde{U}}+\kappa^2\tilde{\tilde{U}}=0$ for $x\in\mathop{\cup}\limits_{m=1}^{M}D_m$ with  $\tilde{\tilde{U}}(x)=0$ 
on $\mathop{\cup}\limits_{m=1}^{M}\partial D_m$.
From the uniqueness of the interior problem, which is true under the condition $a<\frac{1}{\kappa_{\max}}\sqrt[3]{\frac{4\pi}{3}}{\rm j}_{1/2,1}$, 
and the exterior problem we deduce that $\tilde{U}=0$ and $\tilde{\tilde{U}}=0$ in their respective domains. Hence
$\frac{\partial\tilde{U}}{\partial \nu}(x)=0$ and $\frac{\partial\tilde{\tilde{U}}}{\partial \nu}(x)=0$ for $x\in\mathop{\cup}\limits_{m=1}^{M}\partial D_m$. \\
By the jump relations, we have
\begin{eqnarray}\label{nbc1}
\frac{\partial\tilde{U}}{\partial \nu}(x)=0&\Longrightarrow&(\mathbf{K}^{*}\sigma_m)(x)-\frac{\sigma_m(x)}{2}+\sum\limits^{M}_{\substack{j=1\\j\neq m}}\frac{\partial}{\partial \nu_x}\mathcal{S}_{mj}(\sigma_j)(x)=0
\end{eqnarray} 
and 
\begin{eqnarray}\label{nbc2}
\frac{\partial\tilde{\tilde{U}}}{\partial \nu}(x)=0&\Longrightarrow&(\mathbf{K}^{*}\sigma_m)(x)+\frac{\sigma_m(x)}{2}+\sum\limits^{M}_{\substack{j=1\\j\neq m}}\frac{\partial}{\partial \nu_x}\mathcal{S}_{mj}(\sigma_j)(x)=0
\end{eqnarray}
 for $x\in \partial D_m$ and for $m=1,\dots,M$. Here, $\mathbf{K}^{*}$ is the adjoint of the double layer operator $\mathbf{K}$,
\[ (\mathbf{K}\sigma_m)(x):=\int_{\partial D_m}\frac{\partial}{\partial \nu_s}\Phi_\kappa(x,s) \sigma_m(s) ds,\mbox{ for } m=1,\dots,M.\]
Difference between \eqref{nbc1} and \eqref{nbc2} provides us, $\sigma_m=0$ for all $m$. \\
We conclude then that $L+K:=\mathcal{S}:\prod\limits_{m=1}^{M}L^{2}(\partial D_m)\rightarrow \prod\limits_{m=1}^{M}H^{1}(\partial D_m)$ is invertible.
\end{proofn}
\subsubsection{An appropriate estimate of the densities $\sigma_m,\,m=1,\dots,M$}\label{SLPR-2}
From the above theorem, we have the following representation of $\sigma$:
\begin{eqnarray}\label{invLplusK}
 \sigma&=&(L+K)^{-1}U^I \nonumber\\
       &=&L^{-1}(I+L^{-1}K)^{-1}U^I \nonumber\\
       &=&L^{-1}\sum_{l=0}^{\infty}\left(-L^{-1}K\right)^{l}U^I,\,\text{if}\,\left\|L^{-1}K\right\|<1.
\end{eqnarray}
\noindent
The operator $L$ is invertible since it is Fredholm of index zero and injective ( by the assumption $a<\frac{1}{\kappa_{\max}}\sqrt[3]{\frac{4\pi}{3}}{\rm j}_{1/2,1}$). 
This implies that
\begin{eqnarray}\label{nrminvLplusK}
 \left\|\sigma\right\| 
                                  &\leq&\frac{\left\|L^{-1}\right\|}{1-\left\|L^{-1}\right\|\left\|K\right\|}\left\|U^I\right\|.
\end{eqnarray}
Here we use the following notations: 
\begin{eqnarray}
\left\|K\right\|:= \left\|K\right\|_{\mathcal{L}\left(\prod\limits_{m=1}^{M}L^{2}(\partial D_m),\prod\limits_{m=1}^{M}H^{1}(\partial D_m)\right)}
    &\equiv&\max\limits_{1\leq m \leq M}\sum_{j=1}^{M}\left\|K_{mj}\right\|_{\mathcal{L}\left(L^{2}(\partial D_j),H^{1}(\partial D_m)\right)}\nonumber\\
    &=&\max\limits_{1\leq m \leq M}\sum_{\substack{j=1\\j\neq\,m}}^{M}\left\|\mathcal{S}_{mj}\right\|_{\mathcal{L}\left(L^{2}(\partial D_j),H^{1}(\partial D_m)\right)},\label{Knrm}\\
\left\|L^{-1}\right\|:= \left\|L^{-1}\right\|_{\mathcal{L}\left(\prod\limits_{m=1}^{M}H^{1}(\partial D_m),\prod\limits_{m=1}^{M}L^{2}(\partial D_m)\right)}
     &\equiv&\max\limits_{1\leq m \leq M}\sum_{j=1}^{M}\left\|{(L^{-1})}_{mj}\right\|_{\mathcal{L}\left(H^{1}(\partial D_m),L^{2}(\partial D_j)\right)}\nonumber\\
     &=&\max\limits_{1\leq m \leq M}\left\|\mathcal{S}_{mm}^{-1}\right\|_{\mathcal{L}\left(H^{1}(\partial D_m),L^{2}(\partial D_m)\right)},\label{invLnrm}\\
\left\|\sigma\right\|:= \left\|\sigma\right\|_{\prod\limits_{m=1}^{M}L^{2}(\partial D_m)}
    &\equiv&\max\limits_{1\leq m \leq M}\left\|\sigma_{m}\right\|_{L^{2}(\partial D_m)}\label{sigmaU^Inrm1}\\ 
\mbox{and}\qquad\left\|U^I\right\|:= \left\|U^I\right\|_{\prod\limits_{m=1}^{M}H^{1}(\partial D_m)}
    &\equiv&\max\limits_{1\leq m \leq M}\left\|U^{i}\right\|_{H^{1}(\partial D_m)}.\label{sigmaU^Inrm}
\end{eqnarray}

In the following proposition, we provide conditions under which $\left\|L^{-1}\right\|\left\|K\right\|<1$ and then estimate $\left\|\sigma\right\|$ via \eqref{nrminvLplusK}.
\begin{proposition}\label{normofsigmastmt}
There exists a constant $c$ depending only on the 
Lipschitz character of $B_m,m=1,\dots,M$, $d_{\max}$ and $\kappa_{\max}$ such that if $(M-1)\epsilon<cd^2$, 
 then the $L^2$-norm of densities $\sigma_m$ appearing in the solution \eqref{qcimprequiredfrm1} 
of the problem (\ref{acimpoenetrable}-\ref{radiationc}) are uniformly bounded by a positive constant.
 \end{proposition}
\noindent

\textit{Proof of Proposition \ref{normofsigmastmt}}.
 
For any functions $f,g$ defined on $\partial D_\epsilon$ and $\partial B$ respectively, we use the following notations; 
\begin{eqnarray}\label{acsmalla-wedgevee-defntn}
 (f)^\wedge(\xi)\,:=\,\hat{f}(\xi)\,:=\,f(\epsilon\xi+z)&\mbox{ and }& (g)^\vee(x)\,:=\, \check{g}(x)\,:=\,g\left(\frac{x-z}{\epsilon}\right).
\end{eqnarray}
Let $T_1$ and $T_2$ be an orthonormal basis for the tangent plane to $\partial D_\epsilon$ at $x$ and let $\partial\slash\partial\,T=\sum\limits_{l=1}^{2}\partial\slash{\partial\,T_l}\,~T_l$,
denote the tangential derivative on $\partial D_\epsilon$. We recall that the space $H^{1}(\partial D_\epsilon)$ is defined as
\begin{eqnarray}\label{defofH1spaceonbdry}
 H^{1}(\partial D_\epsilon)&:=&\{\phi\in\,L^{2}(\partial D_\epsilon);\partial\phi\slash\partial\,T\in\,L^{2}(\partial D_\epsilon)\}.
\end{eqnarray}
We start with the following lemma which is similar to Lemma 4.1 of \cite{AM-GR-HM:Mathcomp2007}.
\begin{lemma} \label{L2H1estimates}
Suppose $0<\epsilon\leq1$ and $D_\epsilon:=\epsilon B+z\subset\mathbb{R}^n$. Then for each $\phi\in H^{1}(\partial D_\epsilon)$ 
and $\psi \in L^{2}(\partial D_\epsilon)$, we have
 \begin{equation}\label{habib2*}
 \|\psi \|_{L^{2}(\partial D_\epsilon)}=\epsilon^\frac{n-1}{2} \|\hat{\psi} \|_{L^{2}(\partial B)}
\end{equation}
\noindent
and
\begin{equation}\label{habib1*}
\epsilon^{\frac{n-1}{2}}\vert\vert \hat{\phi} \vert\vert_{H^1(\partial B)}~\leq ~\vert\vert \phi \vert\vert_{H^{1}(\partial D_\epsilon)}
~\leq~\epsilon^{\frac{n-3}{2}}\vert\vert \hat{\phi} \vert\vert_{H^{1}(\partial B)}.
\end{equation}
\end{lemma}
\begin{proofn}{\it{of Lemma \ref{L2H1estimates}}.}
\ ~ \ \\

 \par Consider, $\psi \in L^{2}(\partial D_\epsilon)$. Then \eqref{habib2*} is derived as follows
\begin{eqnarray*}
  \|\psi \|^2_{L^2(\partial D_\epsilon)}&=&\int_{\partial D_\epsilon}|\psi(x)|^2 dx \\
                                   &=&\epsilon^{n-1}\int_{\partial B}|\psi(\epsilon\xi+z)|^2 d\xi\\
                                   &=&\epsilon^{n-1}\int_{\partial B}|\hat{\psi}(\xi)|^2 d\xi, (\hat{\psi}(\xi):=\psi(\epsilon\xi+z))\\
                                   &=&\epsilon^{n-1} \|\hat{\psi} \|^2_{L^2(\partial B)}.
\end{eqnarray*}
\par Now consider, $\phi\in H^{1}(\partial D_\epsilon)$. Then
\begin{eqnarray*} 
  \|\phi \|^2_{H^1(\partial D_\epsilon)}&=&\int_{\partial D_\epsilon}(|\phi(x)|^2+|\partial_{T_x}\phi(x)|^2) dx \\
                                   &=&\epsilon^{n-1}\int_{\partial B}\left(|\phi(\epsilon\xi+z)|^2 +\frac{1}{\epsilon^2}|\partial_{T_\xi}\phi(\epsilon\xi+z)|^2\right)d\xi\\
                                   &=&\epsilon^{n-1}\int_{\partial B}|\hat{\phi}(\xi)|^2 d\xi+\epsilon^{n-3}\int_{\partial B}|\partial_{T_\xi}\hat{\phi}(\xi)|^2d\xi.                             
\end{eqnarray*}
Using the fact that $\epsilon^{n-1}\leq\epsilon^{n-3}$, we obtain \eqref{habib1*}.
\end{proofn}
\noindent
We divide the rest of the proof of Proposition \ref{normofsigmastmt} into two steps. In the first step, we assume we have a single obstacle and then in the second step we deal with the multiple obstacle case.
\paragraph{The case of a single obstacle}\label{case of a single obstacle-smallac-sdlp} 
Let us consider a single obstacle $D_\epsilon:=\epsilon B+z$. 
Then define the operator 
$\mathcal{S}_{ D_\epsilon}:L^2(\partial D_{\epsilon})\rightarrow H^1(\partial D_{\epsilon})$ by
\begin{eqnarray}\label{defofSpartialDe}
\left(\mathcal{S}_{ D_\epsilon}\psi\right) (x):=\int_{\partial D_\epsilon}\Phi_\kappa(x,y)\psi(y)dy.
\end{eqnarray}
Following the arguments in the proof of Proposition \ref{existence-of-sigmas}, the integral operator 
$\mathcal{S}_{ D_\epsilon}:L^2(\partial D_{\epsilon})\rightarrow H^1(\partial D_{\epsilon})$ is invertible. 
If we consider the problem (\ref{acimpoenetrable}-\ref{radiationc}) 
in $\mathbb{R}^{3}\backslash \bar{D}_\epsilon$, we obtain
  $$\sigma=\mathcal{S}_{ D_\epsilon}^{-1}U^{i},\text{ where } \mathcal{S}_{ D_\epsilon}=L,$$
and then
\begin{eqnarray}\label{estsigm1}
 \|\sigma \|_{L^2(\partial D_{\epsilon})}\leq  \|\mathcal{S}_{ D_\epsilon}^{-1} \|_{\mathcal{L}\left(H^1(\partial D_{\epsilon}),L^2(\partial D_{\epsilon})\right)} \|U^i \|_{L^2(\partial D_{\epsilon})}.
\end{eqnarray}
\begin{lemma}\label{rep1singulayer}
 Let $\phi\in H^{1}(\partial D_\epsilon)$ and $\psi \in L^{2}(\partial D_\epsilon)$. Then, 
\begin{equation}\label{rep1singulayer1}
 \mathcal{S}_{ D_\epsilon}\psi=\epsilon~ (\mathcal{S}^\epsilon_B \hat{\psi})^\vee,
\end{equation}
\begin{equation}\label{rep1singulayer2}
 \mathcal{S}_{ D_\epsilon}^{-1}\phi=\epsilon^{-1} ({\left(\mathcal{S}^\epsilon_B\right)}^{-1} \hat{\phi})^\vee
\end{equation}
and
\begin{equation}\label{nrm1singulayer2}
 \left\|\mathcal{S}_{ D_\epsilon}^{-1}\right\|_{\mathcal{L}\left(H^1(\partial D_\epsilon), L^2(\partial D_\epsilon) \right)}\leq 
\epsilon^{-1}\left\|{\left(\mathcal{S}^\epsilon_B\right)}^{-1}\right\|_{\mathcal{L}\left(H^1(\partial B), L^2(\partial B) \right)}
\end{equation}
with $\mathcal{S}^\epsilon_B \hat{\psi}(\xi):=\int_{\partial B} \frac{e^{i\kappa\epsilon|\xi-\eta|}}{4\pi|\xi-\eta|}\hat{\psi}(\eta) d\eta$.
\end{lemma}
\begin{proofn}{\it{of Lemma \ref{rep1singulayer}}.}
The identities \eqref{rep1singulayer1} and \eqref{rep1singulayer2} are derived respectively as follows
 \begin{eqnarray*}
  \mathcal{S}_{ D_\epsilon}\psi(x)&=&\int_{\partial D_\epsilon}\frac{e^{i\kappa|x-y|}}{4\pi|x-y|}\psi(y) dy
          \,=\,\int_{\partial B}\frac{e^{i\kappa\epsilon|\xi-\eta|}}{4\pi\epsilon|\xi-\eta|}\psi(\epsilon\eta+z) \epsilon^{2}d\eta
          \,=\,\epsilon ~\mathcal{S}^\epsilon_B \hat{\psi}(\xi)
 \end{eqnarray*}
and
\begin{eqnarray*}
 \mathcal{S}_{ D_\epsilon}({\left(\mathcal{S}^\epsilon_B\right)}^{-1} \hat{\phi})^\vee \,\substack{=\\\eqref{rep1singulayer1}}\, \epsilon~(\mathcal{S}^\epsilon_B{\left(\mathcal{S}^\epsilon_B\right)}^{-1} \hat{\phi})^\vee 
                                        \,=\, \epsilon~\hat{\phi}^{\vee}
                                        \,=\, \epsilon~\phi.                                       
\end{eqnarray*}
To derive the estimate \eqref{nrm1singulayer2}, we proceed as follows
\begin{eqnarray*}
 \left\|\mathcal{S}_{ D_\epsilon}^{-1}\right\|_{\mathcal{L}\left(H^1(\partial D_\epsilon), L^2(\partial D_\epsilon) \right)}&:=&
\substack{Sup\\ \phi(\neq0)\in H^{1}(\partial D_\epsilon)} \frac{ \|\mathcal{S}_{ D_\epsilon}^{-1}\phi \|_{L^2(\partial D_\epsilon)}}{ \|\phi \|_{H^1(\partial D_\epsilon)}}\\
&\substack{\leq\\ \eqref{habib2*},\eqref{habib1*}}\,&\substack{Sup\\ \phi(\neq0)\in H^{1}(\partial D_\epsilon)} 
\frac{\epsilon~\left\|(\mathcal{S}_{ D_\epsilon}^{-1}\phi)^{\wedge}\right\|_{L^2(\partial B)}}{\epsilon~ \|\hat{\phi} \|_{H^1(\partial B)}}\\
&\substack{=\\\eqref{rep1singulayer2}}\,&\substack{Sup\\ \hat{\phi}(\neq0)\in H^{1}(\partial D_\epsilon)} 
\frac{\left\|\left(\epsilon^{-1} ({\left(\mathcal{S}^\epsilon_B\right)}^{-1} \hat{\phi})^\vee\right)^{\wedge}\right\|_{L^2(\partial B)}}{ \|\hat{\phi} \|_{H^1(\partial B)}}\\
&=&\substack{Sup\\ \hat{\phi}(\neq0)\in H^{1}(\partial D_\epsilon)} 
\frac{\epsilon^{-1}\left\| {\left(\mathcal{S}^\epsilon_B\right)}^{-1} \hat{\phi}\right\|_{L^2(\partial B)}}{ \|\hat{\phi} \|_{H^1(\partial B)}} \\
&=&\epsilon^{-1}\left\|{\left(\mathcal{S}^\epsilon_B\right)}^{-1}\right\|_{\mathcal{L}\left(H^1(\partial B), L^2(\partial B) \right)}. 
\end{eqnarray*}
\end{proofn}
From the explicit form of $\mathcal{S}^\epsilon_B$, we can estimate the left hand side of \eqref{nrm1singulayer2} by $\epsilon^{-1}C$ using Banach's theorem of uniform boundedness. 
However this theorem provides only an existence of $C$ with no information on its dependence on $B$ and $\kappa$. The next lemma provides such an estimate.
\begin{lemma}\label{lemmanrm1singulayer31}
  The operator norm of the inverse of $\mathcal{S}_{ D_\epsilon}:L^{2}(\partial D_\epsilon)\rightarrow H^{1}(\partial D_\epsilon)$, defined 
  in \eqref{defofSpartialDe}, is estimated by $\epsilon^{-1}$,
i.e.\begin{eqnarray}\label{nrm1singulayer31}
 \left\|\mathcal{S}_{ D_\epsilon}^{-1}\right\|_{\mathcal{L}\left(H^1(\partial D_\epsilon), L^2(\partial D_\epsilon) \right)}
&\leq&C_6 ~\epsilon^{-1},
\end{eqnarray}
with 
$C_6:=\frac{2\pi\left\|{\mathcal{S}^{i_{\kappa}}_B}^{-1}\right\|_{\mathcal{L}\left(H^1(\partial B), L^2(\partial B) \right)}}
{2\pi-(1+\kappa)\kappa\epsilon|\partial B|\left\|{\mathcal{S}^{i_{\kappa}}_B}^{-1}\right\|_{\mathcal{L}\left(H^1(\partial B), L^2(\partial B) \right)}}$. 
Here, $\mathcal{S}^{i_{\kappa}}_B:L^2(\partial B)\rightarrow H^1(\partial B)$ is the single layer potential with the wave number zero.
\end{lemma}
\noindent

Here we should mention that if 
$\epsilon\leq\frac{\pi}{(1+\kappa)\kappa|\partial B|\left\|{\mathcal{S}^{i_{\kappa}}_B}^{-1}\right\|_{\mathcal{L}\left(H^1(\partial B), L^2(\partial B) \right)}}$,
then $C_6$ is bounded by $2\left\|{\mathcal{S}^{i_{\kappa}}_B}^{-1}\right\|_{\mathcal{L}\left(H^1(\partial B), L^2(\partial B) \right)}$, 
which is a constant depending only on $\partial B$ through its Lipschitz character, see Remark \ref{Lipschitz-character}.
\begin{proofn}{\it{of Lemma \ref{lemmanrm1singulayer31}}.} To estimate the operator norm of $\mathcal{S}_{ D_\epsilon}^{-1}$ we decompose 
$\mathcal{S}_{ D_\epsilon}=:\mathcal{S}_{ D_\epsilon}^\kappa=\mathcal{S}_{ D_\epsilon}^{i_{\kappa}}+\mathcal{S}_{ D_\epsilon}^{d_{\kappa}}$ into 
two parts $\mathcal{S}_{ D_\epsilon}^{i_{\kappa}}$ ( independent of $\kappa$ ) and  $\mathcal{S}_{ D_\epsilon}^{d_{\kappa}}$ ( dependent of $\kappa$ ) given by
\begin{eqnarray}
\mathcal{S}_{ D_\epsilon}^{i_{\kappa}} \psi(x):=\int_{\partial D_\epsilon} \frac{1}{4\pi|x-y|}\psi(y) dy,\label{def-of_Sik}\\
\mathcal{S}_{ D_\epsilon}^{d_{\kappa}} \psi(x):=\int_{\partial D_\epsilon} \frac{e^{i\kappa|x-y|}-1}{4\pi|x-y|}\psi(y) dy.\label{def-of_Sdk}
\end{eqnarray}
 With this definition, $\mathcal{S}_{ D_\epsilon}^{i_{\kappa}}:L^2(\partial D_\epsilon)\rightarrow H^1(\partial D_\epsilon)$ is invertible, see \cite{Mclean:2000,MD:InteqnsaOpethe1997}. Hence, 
$\mathcal{S}_{ D_\epsilon}=\mathcal{S}_{ D_\epsilon}^{i_{\kappa}}\left(I+{\mathcal{S}_{ D_\epsilon}^{i_{\kappa}}}^{-1}\mathcal{S}_{ D_\epsilon}^{d_{\kappa}}\right)$ and so
\begin{eqnarray}
 \left\|\mathcal{S}_{ D_\epsilon}^{-1}\right\|_{\mathcal{L}\left(H^1(\partial D_\epsilon), L^2(\partial D_\epsilon) \right)}
&\leq& \left\|\left(I+{\mathcal{S}_{ D_\epsilon}^{i_{\kappa}}}^{-1}\mathcal{S}_{ D_\epsilon}^{d_{\kappa}}\right)^{-1}\right\|_{\mathcal{L}\left(L^2(\partial D_\epsilon), L^2(\partial D_\epsilon) \right)} 
                            \left\|{\mathcal{S}_{ D_\epsilon}^{i_{\kappa}}}^{-1}\right\|_{\mathcal{L}\left(H^1(\partial D_\epsilon), L^2(\partial D_\epsilon) \right)}.\label{nrm1singulayer3}
\end{eqnarray}
So, to estimate the operator norm of $\mathcal{S}_{ D_\epsilon}^{-1}$ one needs to estimate the 
operator norm of $\left(I+{\mathcal{S}_{ D_\epsilon}^{i_{\kappa}}}^{-1}\mathcal{S}_{ D_\epsilon}^{d_{\kappa}}\right)^{-1}$, 
in particular one needs to have the knowledge about the operator norms of
 ${\mathcal{S}_{ D_\epsilon}^{i_{\kappa}}}^{-1}$ and $\mathcal{S}_{ D_\epsilon}^{d_{\kappa}}$ to apply the Neumann series. 
For that purpose,  we can estimate the operator norm of  ${\mathcal{S}_{ D_\epsilon}^{i_{\kappa}}}^{-1}$ from \eqref{nrm1singulayer2} by
\begin{eqnarray}\label{optnormsik}
 \left\|{\mathcal{S}_{ D_\epsilon}^{i_{\kappa}}}^{-1}\right\|_{\mathcal{L}\left(H^1(\partial D_\epsilon), L^2(\partial D_\epsilon) \right)}
   &\leq&\epsilon^{-1} \left\|{\mathcal{S}^{i_{\kappa}}_B}^{-1}\right\|_{\mathcal{L}\left(H^1(\partial B), L^2(\partial B) \right)}.
\end{eqnarray}
Here $\mathcal{S}^{i_{\kappa}}_B \hat{\psi}(\xi):=\int_{\partial B} \frac{1}{4\pi|\xi-\eta|}\hat{\psi}(\eta) d\eta$. 
From the definition of the operator $S_{ D_\epsilon}^{d_{k}}$ in \eqref{def-of_Sdk}, we deduce that
\begin{eqnarray}\label{S^{d_{k}}}
\mathcal{S}_{ D_\epsilon}^{d_{\kappa}} \psi(x)
             &=&\epsilon~\int_{\partial B} \frac{e^{i\kappa\epsilon|\xi-\eta|}-1}{4\pi|\xi-\eta|}\hat{\psi}(\eta) d\eta
\end{eqnarray}
 and so the tangential derivative with respect to basis vector $T~(T_1\text{ or } T_2)$ is given by
\begin{eqnarray}\label{TanS^{d_{k}}}
\frac{\partial}{\partial T}\mathcal{S}_{ D_\epsilon}^{d_{\kappa}} \psi(x)&=&\int_{\partial D_\epsilon} 
     \left(\frac{e^{i\kappa|x-y|}}{4\pi|x-y|}\left[i\kappa-\frac{1}{|x-y|}\right]\frac{x-y}{|x-y|}+\frac{1}{4\pi}\frac{x-y}{|x-y|^3}\right)\cdot T_x \psi(y) dy\nonumber \\
&=&\int_{\partial B} 
     \left(\frac{e^{i\kappa\epsilon|\xi-\eta|}}{4\pi\epsilon|\xi-\eta|}\left[i\kappa-\frac{1}{\epsilon|\xi-\eta|}\right]\frac{\xi-\eta}{|\xi-\eta|}
      +\frac{1}{4\pi}\frac{\xi-\eta}{\epsilon^2|\xi-\eta|^3}\right)\cdot T_\xi \hat{\psi}(\eta) \epsilon^{2}d\eta\nonumber\\
&=&\int_{\partial B} 
     \left(\frac{e^{i\kappa\epsilon|\xi-\eta|}}{4\pi|\xi-\eta|}i\kappa\epsilon-\frac{e^{i\kappa\epsilon|\xi-\eta|}-1}{4\pi|\xi-\eta|^2}\right)\frac{\xi-\eta}{|\xi-\eta|}
      \cdot T_\xi \hat{\psi}(\eta) d\eta\nonumber\\
&=&\int_{\partial B} 
    \frac{(i\kappa\epsilon)^2}{4\pi} \left[\sum_{l=1}^{\infty}\left(\frac{1}{l!}-\frac{1}{l+1!}\right)(i\kappa\epsilon)^{l-1}|\xi-\eta|^{l-1}\right]
\frac{\xi-\eta}{|\xi-\eta|} \cdot T_\xi \hat{\psi}(\eta) d\eta\nonumber\\
&=&\frac{-\kappa^2\epsilon^{2}}{4\pi}\int_{\partial B} 
     \left[\sum_{l=1}^{\infty}\left(\frac{1}{l!}-\frac{1}{l+1!}\right)(i\kappa\epsilon)^{l-1}|\xi-\eta|^{l-1}\right]
\frac{\xi-\eta}{|\xi-\eta|}\cdot T_\xi \hat{\psi}(\eta) d\eta,
\end{eqnarray}
where we used the following expansions
\begin{eqnarray*}
\frac{e^{i\kappa\epsilon|\xi-\eta|}}{4\pi|\xi-\eta|}i\kappa\epsilon&=&\frac{i\kappa\epsilon}{4\pi|\xi-\eta|}+\frac{1}{1!}\frac{(i\kappa\epsilon)^2}{4\pi}+\frac{1}{2!}\frac{(i\kappa\epsilon)^3}{4\pi}|\xi-\eta|+\dots \\
 &=&\frac{i\kappa\epsilon}{4\pi|\xi-\eta|}+\frac{(i\kappa\epsilon)^2}{4\pi}\sum_{l=1}^{\infty}\frac{(i\kappa\epsilon)^{l-1}}{l!}|\xi-\eta|^{l-1}
\end{eqnarray*}
and
\begin{eqnarray*}
\frac{e^{i\kappa\epsilon|\xi-\eta|}-1}{4\pi|\xi-\eta|^2}&=&\frac{1}{1!}\frac{i\kappa\epsilon}{4\pi|\xi-\eta|}+\frac{1}{2!}\frac{(i\kappa\epsilon)^2}{4\pi}+\frac{1}{3!}\frac{(i\kappa\epsilon)^3}{4\pi}|\xi-\eta|+\dots \nonumber \\
 &=&\frac{i\kappa\epsilon}{4\pi|\xi-\eta|}+\frac{(i\kappa\epsilon)^2}{4\pi}\sum_{l=1}^{\infty}\frac{(i\kappa\epsilon)^{l-1}}{l+1!}|\xi-\eta|^{l-1}.
\end{eqnarray*}
From \eqref{S^{d_{k}}}, we observe that
 \begin{eqnarray}\label{modS^{d_{k}}}
  |\mathcal{S}_{ D_\epsilon}^{d_{\kappa}} \psi(x)|&\leq&\epsilon~\left\|\frac{e^{i\kappa\epsilon|\xi-\cdot|}-1}{4\pi|\xi-\cdot|}\right\|_{L^2(\partial B)}\|\hat{\psi}\|_{L^2(\partial B)}\nonumber\\
&=&\frac{\epsilon}{4\pi}\left\|i\kappa\epsilon\sum_{l=1}^{\infty}\frac{(i\kappa\epsilon)^{l-1}}{l!}|\xi-\cdot|^{l-1}\right\|_{L^2(\partial B)}
\|\hat{\psi}\|_{L^2(\partial B)}\nonumber\\
&\leq&\frac{\kappa\epsilon^{2}}{4\pi}
\sum_{l=1}^{\infty}\frac{(\kappa\epsilon)^{l-1}}{l!}\left\|\xi-\cdot\right\|^{l-1}_{L^2(\partial B)}|\partial B|^{\frac{2-l}{2}}\|\hat{\psi}\|_{L^2(\partial B)}\nonumber\\
&&(\mbox{Since, }\left\|~|x|^p\right\|_{L^2(D)}\leq\left\|x\right\|^p_{L^2(D)}|D|^{\frac{1-p}{2}})\nonumber\\
&\leq&\frac{\kappa\epsilon^{2}}{4\pi}|\partial B|^{\frac{1}{2}}
\sum_{l=1}^{\infty}\left[\frac{1}{2}\kappa\epsilon\left\|\xi-\cdot\right\|_{L^2(\partial B)}|\partial B|^{\frac{-1}{2}}\right]^{l-1}\|\hat{\psi}\|_{L^2(\partial B)}\nonumber\\
&=&\frac{\kappa\epsilon^{2}}{4\pi}|\partial B|^{\frac{1}{2}}
\frac{\|\hat{\psi}\|_{L^2(\partial B)}}{1-\frac{1}{2}\kappa\epsilon\left\|\xi-\cdot\right\|_{L^2(\partial B)}|\partial B|^{\frac{-1}{2}}},
\,\text{if}\,\left(\epsilon<\frac{2}{\kappa_{\max}\max\limits_{1\leq m \leq M}diam(B_m)}\right)\equiv (\kappa_{\max}\,a<2)\nonumber\\
&\leq&\frac{\kappa\epsilon^{2}}{2\pi}|\partial B|^{\frac{1}{2}}
\|\hat{\psi}\|_{L^2(\partial B)},\,\text{for}\,\left(\epsilon\leq\frac{1}{\kappa_{\max}\max\limits_{1\leq m \leq M}diam(B_m)}\right)\equiv (\kappa_{\max}\,a\leq1).
 \end{eqnarray}
We set $C_1:=\frac{|\partial B|^{\frac{1}{2}}}{2\pi}$.
From this we get,
\begin{eqnarray}
  \|\mathcal{S}_{ D_\epsilon}^{d_{\kappa}} \psi \|^2_{L^2(\partial D_\epsilon)}&=&\int_{\partial D_\epsilon}|\mathcal{S}_{ D_\epsilon}^{d_{\kappa}} \psi(x)|^2 dx\nonumber \\
                                          &\substack{\leq\\ \eqref{modS^{d_{k}}}}&\int_{\partial D_\epsilon}\left[C_1\kappa\epsilon^{2}
\|\hat{\psi}\|_{L^2(\partial B)}\right]^2 dx\nonumber \\
&=&C_1^2\kappa^2\epsilon^{4}
\|\hat{\psi}\|^2_{L^2(\partial B)}\int_{\partial D_\epsilon}dx \nonumber \\
&=&C_1^2\kappa^2|\partial B|\epsilon^{6}
\|\hat{\psi}\|^2_{L^2(\partial B)}, \nonumber 
\end{eqnarray}
which gives us 
\begin{eqnarray}\label{normS^{d_{k}}}
  \|\mathcal{S}_{ D_\epsilon}^{d_{\kappa}} \psi \|_{L^2(\partial D_\epsilon)}&\leq&C_1\kappa\epsilon^{3}|\partial B|^{\frac{1}{2}}
\|\hat{\psi}\|_{L^2(\partial B)}.
\end{eqnarray}
From  \eqref{TanS^{d_{k}}}, we have
$$\frac{\partial}{\partial T}\mathcal{S}_{ D_\epsilon}^{d_{\kappa}} \psi(x)=\frac{-\kappa^2\epsilon^{2}}{4\pi}\int_{\partial B} 
     \left[\sum_{l=1}^{\infty}\left(\frac{1}{l!}-\frac{1}{l+1!}\right)(i\kappa\epsilon)^{l-1}|\xi-\eta|^{l-1}\right]
\frac{\xi-\eta}{|\xi-\eta|}
      \cdot T_\xi~ \hat{\psi}(\eta) d\eta,$$
which gives us
 \begin{eqnarray}\label{modTanS^{d_{k}}}
\left|\frac{\partial}{\partial T}\mathcal{S}_{ D_\epsilon}^{d_{\kappa}} \psi(x)\right|&\leq&\frac{\kappa^2\epsilon^{2}}{4\pi}
 \left\|\sum_{l=1}^{\infty}\left(\frac{1}{l!}-\frac{1}{l+1!}\right)(i\kappa\epsilon)^{l-1}|\xi-\cdot|^{l-1}\right\|_{L^2(\partial B)} \|\hat{\psi}\|_{L^2(\partial B)} \nonumber\\
&=&\frac{\kappa^2\epsilon^{2}}{4\pi}
 \left\|\sum_{l=1}^{\infty}\frac{l}{l+1!}(i\kappa\epsilon)^{l-1}|\xi-\cdot|^{l-1}\right\|_{L^2(\partial B)} \|\hat{\psi}\|_{L^2(\partial B)} \nonumber\\
&\leq&\frac{\kappa^2\epsilon^2}{4\pi}
\sum_{l=1}^{\infty}\frac{l}{l+1!}(\kappa\epsilon)^{l-1}\vert|\xi-\cdot\vert|^{l-1}_{L^2(\partial B)}|\partial B|^{\frac{2-l}{2}} \|\hat{\psi}\|_{L^2(\partial B)} \nonumber\\
&\leq&\frac{\kappa^2\epsilon^2}{4\pi}|\partial B|^{\frac{1}{2}}\frac{1}{2}
\sum_{l=1}^{\infty}\left[\frac{1}{2}\kappa\epsilon\vert|\xi-\cdot\vert|_{L^2(\partial B)}|\partial B|^{\frac{-1}{2}}\right]^{l-1} \|\hat{\psi}\|_{L^2(\partial B)} \nonumber\\
&=&\frac{\kappa^2\epsilon^2}{8\pi}|\partial B|^{\frac{1}{2}}
\frac{\|\hat{\psi}\|_{L^2(\partial B)}}{1-\frac{1}{2}\kappa\epsilon\vert|\xi-\cdot\vert|_{L^2(\partial B)}|\partial B|^{\frac{-1}{2}}},
\,\text{if}\,\left(\epsilon<\frac{2}{\kappa_{\max}\max\limits_{1\leq m \leq M}diam(B_m)}\right)\equiv (\kappa_{\max}~a<2)  \nonumber\\
&\leq&\frac{\kappa^2\epsilon^{2}}{4\pi}|\partial B|^{\frac{1}{2}}
\|\hat{\psi}\|_{L^2(\partial B)},\,\text{for}\,\left(\epsilon\leq\frac{1}{\kappa_{\max}\max\limits_{1\leq m \leq M}diam(B_m)}\right)\equiv (\kappa_{\max}~a\leq1),\nonumber\\
&=&\frac{1}{2}C_1\kappa^2\epsilon^{2}
\|\hat{\psi}\|_{L^2(\partial B)}.
 \end{eqnarray}
From this we obtain,
\begin{eqnarray*}
 \left\|\frac{\partial}{\partial T}\mathcal{S}_{ D_\epsilon}^{d_{\kappa}} \psi\right\|^2_{L^2(\partial D_\epsilon)}
            &=&\int_{\partial D_\epsilon}\left|\frac{\partial}{\partial T}\mathcal{S}_{ D_\epsilon}^{d_{\kappa}} \psi(x)\right|^2 dx\nonumber \\
            &\substack{\leq\\ \eqref{modTanS^{d_{k}}}}&\int_{\partial D_\epsilon}\left[\frac{1}{2}C_1 \kappa^2\epsilon^{2}
 \|\hat{\psi}\|_{L^2(\partial B)}\right]^2 dx\nonumber\\  
 &=&\frac{1}{4}C_1^2\kappa^4\epsilon^{4}
 \|\hat{\psi}\|^2_{L^2(\partial B)}\int_{\partial D_\epsilon} dx\nonumber\\    
 &=&\frac{1}{4}C_1^2\kappa^4\epsilon^{6}|\partial B|
 \|\hat{\psi}\|^2_{L^2(\partial B)}.\nonumber\\ 
\end{eqnarray*}
Hence
\begin{eqnarray}\label{normTanS^{d_{k}}}
 \left\|\frac{\partial}{\partial T}\mathcal{S}_{ D_\epsilon}^{d_{\kappa}} \psi\right\|_{L^2(\partial D_\epsilon)}&\leq&
\frac{1}{2}C_1 \kappa^2\epsilon^{3}|\partial B|^\frac{1}{2}
  \|\hat{\psi}\|_{L^2(\partial B)}.
\end{eqnarray}
Now, we have
\begin{eqnarray*}
\left\|\mathcal{S}_{ D_\epsilon}^{d_{\kappa}} \psi\right\|_{H^1(\partial D_\epsilon)}&\leq&\left\|\mathcal{S}_{ D_\epsilon}^{d_{\kappa}} \psi\right\|_{L^2(\partial D_\epsilon)}
+\sum_{l=1}^{2}\left\|\frac{\partial}{\partial T_l}\mathcal{S}_{ D_\epsilon}^{d_{\kappa}} \psi\right\|_{L^2(\partial D_\epsilon)}
\end{eqnarray*}
and so from \eqref{normS^{d_{k}}} and  \eqref{normTanS^{d_{k}}}, we can write
\begin{eqnarray}\label{normTanS^{d_{k}}H1}
 \left\|\mathcal{S}_{ D_\epsilon}^{d_{\kappa}} \psi\right\|_{H^1(\partial D_\epsilon)}&\leq&
C_1\kappa\epsilon^{3}|\partial B|^{\frac{1}{2}}\|\hat{\psi}\|_{L^2(\partial B)}
+C_1 \kappa^2\epsilon^{3}|\partial B|^\frac{1}{2} \|\hat{\psi}\|_{L^2(\partial B)}\nonumber\\
&=&C_1(1+\kappa)\kappa\epsilon^{3}|\partial B|^\frac{1}{2}\|\hat{\psi}\|_{L^2(\partial B)}.
\end{eqnarray}
We estimate the norm of the operator $\mathcal{S}_{ D_\epsilon}^{d_{\kappa}}$ as
\begin{eqnarray}\label{fnormS^{d_{k}}H1}
 \left\|\mathcal{S}_{ D_\epsilon}^{d_{\kappa}}\right\|_{\mathcal{L}\left(L^2(\partial D_\epsilon), H^1(\partial D_\epsilon)\right)}&=&
\substack{Sup\\ \psi(\neq0)\in L^{2}(\partial D_\epsilon)} \frac{ \|\mathcal{S}_{ D_\epsilon}^{d_{\kappa}}\psi \|_{H^1(\partial D_\epsilon)}}{ \|\psi \|_{L^2(\partial D_\epsilon)}}\nonumber\\
&\leq&
 \substack{Sup\\ \hat{\psi}(\neq0)\in L^{2}(\partial B)} \frac{C_1(1+\kappa)\kappa\epsilon^{3}|\partial B|^\frac{1}{2}\|\hat{\psi}\|_{L^2(\partial B)}}
{\epsilon~\|\hat{\psi}\|_{L^2(\partial B)}}~~(\mbox{ By } \eqref{normTanS^{d_{k}}H1},\eqref{habib2*})\nonumber\\
&=&C_1(1+\kappa)\kappa\epsilon^{2}|\partial B|^\frac{1}{2}.
\end{eqnarray}
Hence, we get
 \begin{eqnarray}\label{norms0sk}
  \left\|{\mathcal{S}_{ D_\epsilon}^{i_{\kappa}}}^{-1}\mathcal{S}_{ D_\epsilon}^{d_{\kappa}}\right\|_{\mathcal{L}\left(L^2(\partial D_\epsilon), L^2(\partial D_\epsilon) \right)}
&\leq&
\left\|{\mathcal{S}_{ D_\epsilon}^{i_{\kappa}}}^{-1}\right\|_{\mathcal{L}\left(H^1(\partial D_\epsilon), L^2(\partial D_\epsilon) \right)}
\left\|\mathcal{S}_{ D_\epsilon}^{d_{\kappa}}\right\|_{\mathcal{L}\left(L^2(\partial D_\epsilon), H^1(\partial D_\epsilon) \right)}\nonumber\\
&\substack{\leq\\
\eqref{optnormsik},\eqref{fnormS^{d_{k}}H1}}&\epsilon^{-1}\left\|{\mathcal{S}^{i_{\kappa}}_B}^{-1}\right\|_{\mathcal{L}\left(H^1(\partial B), L^2(\partial B) \right)}C_1(1+\kappa)\kappa\epsilon^{2}|\partial B|^\frac{1}{2}\nonumber\\
&=&C_2(1+ \kappa) \kappa\epsilon ,
 \end{eqnarray}
where $C_2:=C_1|\partial B|^\frac{1}{2}\left\|{\mathcal{S}^{i_{\kappa}}_B}^{-1}\right\|_{\mathcal{L}\left(H^1(\partial B), L^2(\partial B) \right)}
=\frac{|\partial B|}{2\pi}\left\|{\mathcal{S}^{i_{\kappa}}_B}^{-1}\right\|_{\mathcal{L}\left(H^1(\partial B), L^2(\partial B) \right)}.$
Assuming $\epsilon$ to satisfy the condition $\epsilon<\frac{1}{C_2(1+ \kappa_{\max})\kappa_{\max}}$, 
then $\left\|{\mathcal{S}_{ D_\epsilon}^{i_{\kappa}}}^{-1}\mathcal{S}_{ D_\epsilon}^{d_{\kappa}}\right\|_{\mathcal{L}\left(L^2(\partial D_\epsilon), L^2(\partial D_\epsilon) \right)}<1$ 
and hence by using the Neumann series we obtain the following
\begin{eqnarray*}
 \left\|\left(I+{\mathcal{S}_{ D_\epsilon}^{i_{\kappa}}}^{-1}\mathcal{S}_{ D_\epsilon}^{d_{\kappa}}\right)^{-1}\right\|_{\mathcal{L}\left(L^2(\partial D_\epsilon), L^2(\partial D_\epsilon) \right)}
 &\leq&
 \frac{1}{1-\left\|{\mathcal{S}_{ D_\epsilon}^{i_{\kappa}}}^{-1}\mathcal{S}_{ D_\epsilon}^{d_{\kappa}}\right\|_{\mathcal{L}\left(L^2(\partial D_\epsilon), L^2(\partial D_\epsilon) \right)}}\\
 &\substack{\leq\\  \eqref{norms0sk}}&C_3:=\frac{1}{1-C_2( 1+\kappa) \kappa\epsilon}.\\
\end{eqnarray*}
 By substituting the above and \eqref{optnormsik} in \eqref{nrm1singulayer3}, we obtain the required result \eqref{nrm1singulayer31}.
\end{proofn}

\paragraph{The multiple obstacle case}\label{gbmulobscase}
\ ~ \ \\
\begin{lemma} \label{propsmjsmmest}
 For $m,j=1,2,\dots,M$, the operator $\mathcal{S}_{mj}:L^2(\partial D_j)\rightarrow H^{1}(\partial D_m)$ defined 
in Proposition \ref{existence-of-sigmas}, see \eqref{defofSmjed}, 
satisfies the following estimates,
\begin{itemize}
\item
For $j=m$,
\begin{eqnarray}\label{estinvsmm}
\left\|\mathcal{S}_{mm}^{-1}\right\|_{\mathcal{L}\left(H^1(\partial D_m), L^2(\partial D_m) \right)}&\leq&C_{6m}~ \epsilon^{-1}, \label{invfnormS_{ii}1}
\end{eqnarray}
where $C_{6m}:=\frac{2\pi\left\|{\mathcal{S}^{i_{\kappa}}_{B_m}}^{-1}\right\|_{\mathcal{L}\left(H^1(\partial B_m), L^2(\partial B_m) \right)}}{2\pi-(1+\kappa)\kappa\epsilon|\partial B_m|\left\|{\mathcal{S}^{i_{\kappa}}_{B_m}}^{-1}\right\|_{\mathcal{L}\left(H^1(\partial B_m), L^2(\partial B_m) \right)}}$.
\item For $j\neq m$,
\begin{eqnarray}\label{estinvsmj}
\left\|\mathcal{S}_{mj}\right\|_{\mathcal{L}\left(L^2(\partial D_j), H^1(\partial D_m)\right)}
                          &\leq&\frac{1}{4\pi}\left(\frac{2\kappa+1}{d}+\frac{2}{d^2}\right)\left|\partial \c{B} \right|\epsilon^{2},\label{fnormS_{ij}21}
\end{eqnarray}
where $\left|\partial \c{B} \right|:=\max\limits_m \left|\partial B_m\right|$.

\end{itemize}
\end{lemma}
\begin{proofn}{\it{of Lemma \ref{propsmjsmmest}}.}
The estimate \eqref{estinvsmm} is nothing else but \eqref{nrm1singulayer31} of Lemma \ref{lemmanrm1singulayer31}, replacing $B$ by $B_m$, $z$ by $z_m$ and $D_\epsilon$ by $D_m$ respectively. It remains to prove the estimate \eqref{estinvsmj}.

We have
\begin{eqnarray}\label{fnormS_{ij}}
 \left\|\mathcal{S}_{mj}\right\|_{\mathcal{L}\left(L^2(\partial D_j), H^1(\partial D_m)\right)}&=&
\substack{Sup\\ \psi(\neq0)\in L^{2}(\partial D_j)} \frac{ \|\mathcal{S}_{mj}\psi \|_{H^1(\partial D_m)}}{ \|\psi \|_{L^2(\partial D_j)}}\nonumber\\
&\leq&\substack{Sup\\ \psi(\neq0)\in L^{2}(\partial D_j)} \frac{ \|\mathcal{S}_{mj}\psi \|_{L^2(\partial D_m)}+\sum\limits_{l=1}^{2} \|\partial_{T_l}\mathcal{S}_{mj}\psi \|_{L^2(\partial D_m)}}{ \|\psi \|_{L^2(\partial D_j)}}.
\end{eqnarray}
\noindent
Let $\psi\in L^2(\partial D_j)$ then for $x\in \partial D_m$,  we have
\begin{eqnarray}\label{modS_{ij}}
\left|\mathcal{S}_{mj}\psi(x)\right|&=&\left|\int_{\partial D_j}\Phi_\kappa(x,s)\psi(s)ds\right|\nonumber\\
                                    &\leq&\int_{\partial D_j}\frac{1}{4\pi\left|x-s\right|}\left|\psi(s)\right|ds\nonumber\\
                                    &\leq&\frac{1}{4\pi d_{mj}}\int_{\partial D_j}\left|\psi(s)\right|ds\hspace{1cm} \left(d_{mj}:=dist(D_m, D_j)\right)\nonumber\\          
                                    &\leq&\frac{1}{4\pi d_{mj}}\epsilon~\left|\partial B_j\right|^\frac{1}{2} \|\psi \|_{L^2(\partial D_j)},                           
\end{eqnarray}
\noindent
from which, we obtain
\begin{eqnarray}\label{fL2normS_{ij}}
\left\|\mathcal{S}_{mj}\psi\right\|_{L^2(\partial D_m)}&=&\left(\int_{\partial D_m}\left|\mathcal{S}_{mj}\psi(x)\right|^2 dx\right)^\frac{1}{2}\nonumber\\
                                                                  &\substack{\leq\\ \eqref{modS_{ij}}}&\frac{1}{4\pi d_{mj}}\epsilon~\left|\partial B_j\right|^\frac{1}{2} \|\psi \|_{L^2(\partial D_j)}
                                                                         \left(\int_{\partial D_m} dx\right)^\frac{1}{2}\nonumber\\
                                                                  &=&\frac{1}{4\pi d_{mj}}\epsilon^{2}\left|\partial B_j\right|^\frac{1}{2}
                                                                           \left|\partial B_m\right|^\frac{1}{2} \|\psi \|_{L^2(\partial D_j)}.                                                                
\end{eqnarray}
Now for the basis tangential vector $T$ ($T_1$ or $T_2$), we have
\begin{eqnarray}\label{modTS_{ij}}
\left|\frac{\partial}{\partial T}\mathcal{S}_{mj}\psi(x)\right|&=&\left|\int_{\partial D_j}\frac{\partial}{\partial T}\Phi_\kappa(x,s)\psi(s)ds\right|\nonumber\\
                                                               &=&\left|\int_{\partial D_j}\nabla_x\Phi_\kappa(x,s)\cdot T_x\psi(s) ds\right|\nonumber\\
                                                               &\leq&\int_{\partial D_j}\left|\nabla_x\Phi_\kappa(x,s)\right|~\left|\psi(s)\right| ds\nonumber\\
                                                               &=&\int_{\partial D_j}\left|\frac{e^{i\kappa|x-s|}}{4\pi|x-s|}\left(i\kappa-\frac{1}{|x-s|}\right)\frac{x-s}{|x-s|}\right|
                                                                                         ~\left|\psi(s)\right| ds\nonumber\\
                                                               &\leq&\int_{\partial D_j}\frac{1}{4\pi|x-s|}\left(\kappa+\frac{1}{|x-s|}\right)
                                                                                         ~\left|\psi(s)\right| ds\nonumber\\    
                                                               &\leq&\frac{1}{4\pi}\left(\frac{\kappa}{d_{mj}}+\frac{1}{d_{mj}^2}\right)\int_{\partial D_j}\left|\psi(s)\right| ds\nonumber\\       
                                                               &\leq&\frac{1}{4\pi}\left(\frac{\kappa}{d_{mj}}+\frac{1}{d_{mj}^2}\right)\epsilon
                                                                             \left|\partial B_j\right|^\frac{1}{2} \|\psi \|_{L^2(\partial D_j)},                  
\end{eqnarray}
\noindent
from which, we obtain
\begin{eqnarray}\label{fL2normTS_{ij}}
\left\|\frac{\partial}{\partial T}\mathcal{S}_{mj}\psi\right\|_{L^2(\partial D_m)}
                                                                  &=&\left(\int_{\partial D_m}\left|\frac{\partial}{\partial T}\mathcal{S}_{mj}\psi(x)\right|^2 dx\right)^\frac{1}{2}\nonumber\\
                                                                  &\substack{\leq\\ \eqref{modTS_{ij}}}&\frac{1}{4\pi}\left(\frac{\kappa}{d_{mj}}+\frac{1}{d_{mj}^2}\right)\epsilon
                                                                             \left|\partial B_j\right|^\frac{1}{2} \|\psi \|_{L^2(\partial D_j)}\left(\int_{\partial D_m} dx\right)^\frac{1}{2}\nonumber\\
                                                                  &=&\frac{1}{4\pi}\left(\frac{\kappa}{d_{mj}}+\frac{1}{d_{mj}^2}\right)\epsilon^{2}
                                                                                             \left|\partial B_j\right|^\frac{1}{2}\left|\partial B_m\right|^\frac{1}{2} \|\psi \|_{L^2(\partial D_j)}.                                                                
\end{eqnarray}
\noindent
From \eqref{defofH1spaceonbdry}, \eqref{fL2normS_{ij}} and \eqref{fL2normTS_{ij}}, we derive
\begin{eqnarray}\label{fH1normS_{ij}}
\left\|\mathcal{S}_{mj}\psi\right\|_{H^1(\partial D_m)}&\leq&\frac{1}{4\pi}\left(\frac{2\kappa+1}{d_{mj}}+\frac{2}{d_{mj}^2}\right)\epsilon^{2}
                                                                                             \left|\partial B_j\right|^\frac{1}{2}\left|\partial B_m\right|^\frac{1}{2} \|\psi \|_{L^2(\partial D_j)}.
\end{eqnarray}
\noindent
Substitution of \eqref{fH1normS_{ij}} in \eqref{fnormS_{ij}} gives us
\begin{eqnarray*}\label{fnormS_{ij}1}
 \left\|\mathcal{S}_{mj}\right\|_{\mathcal{L}\left(L^2(\partial D_j), H^1(\partial D_m)\right)}
                          &\leq&\frac{1}{4\pi}\left(\frac{2\kappa+1}{d_{mj}}+\frac{2}{d_{mj}^2}\right)\epsilon^2
                                                                   \left|\partial B_j\right|^\frac{1}{2}\left|\partial B_m\right|^\frac{1}{2}
                          \,\leq\,\frac{1}{4\pi}\left(\frac{2\kappa+1}{d}+\frac{2}{d^2}\right)\left|\partial \c{B}\right|\epsilon^{2}.
\end{eqnarray*}
\end{proofn}

\begin{proofe}
\textbf{\textit{End of the proof of Proposition \ref{normofsigmastmt}}.}

By substituting \eqref{invfnormS_{ii}1} in \eqref{invLnrm} and \eqref{fnormS_{ij}21} in \eqref{Knrm}, we obtain
\begin{eqnarray}
\left\|K\right\|
    &\equiv&
\max\limits_{1\leq m \leq M}\sum_{\substack{j=1\\j\neq\,m}}^{M}\left\|\mathcal{S}_{mj}\right\|_{\mathcal{L}\left(L^{2}(\partial D_j),H^{1}(\partial D_m)\right)}\nonumber\\
    &\leq&\frac{M-1}{4\pi}\left(\frac{2\kappa+1}{d}+\frac{2}{d^2}\right)\left|\partial \c{B}\right|\epsilon^{2}\label{Knrm1}
\end{eqnarray}
and 
\begin{eqnarray}
\left\|L^{-1}\right\|
     &\equiv&\max\limits_{1\leq m \leq M}\left\|\mathcal{S}_{mm}^{-1}\right\|_{\mathcal{L}\left(H^{1}(\partial D_m),L^{2}(\partial D_m)\right)}\nonumber\\ 
     &\equiv&\left(\max\limits_{1\leq m \leq M}C_{6m}\right) \epsilon^{-1}.  \label{invLnrm1}
\end{eqnarray}
Hence, \eqref{Knrm1} and \eqref{invLnrm1} jointly provide
\begin{eqnarray}
\left\|L^{-1}\right\|\left\|K\right\|
    &\leq&\underbrace{\frac{M-1}{4\pi}\left(\max\limits_{1\leq m \leq M}C_{6m}\right)\left|\partial \c{B}\right|\left(\frac{2\kappa+1}{d}+\frac{2}{d^2}\right)\epsilon}_{=:C_s} ,\label{invKnrm1}
\end{eqnarray}
By imposing the condition $\left\|L^{-1}\right\|\left\|K\right\|<1$, we have from \eqref{nrminvLplusK} and (\ref{sigmaU^Inrm1}-\ref{sigmaU^Inrm});
\begin{eqnarray}\label{nrminvLplusK2}
 \left\|\sigma_m\right\|_{L^{2}(\partial D_m)}\leq\left\|\sigma\right\| 
                                  &\leq&\frac{\left\|L^{-1}\right\|}{1-\left\|L^{-1}\right\|\left\|K\right\|}\left\|U^I\right\|\nonumber\\
                   &\leq&C_p\left\|L^{-1}\right\| \max\limits_{1\leq m \leq M}\left\|U^{i}\right\|_{H^{1}(\partial D_m)}\hspace{.25cm} \left( C_p\geq\frac{1}{1-C_s},\mbox{ a  positive constant}\right)\nonumber\\
                   &\substack{\leq\\ \eqref{invLnrm1} }&\mathrm{C} \, \epsilon^{-1}\max\limits_{1\leq m \leq M}\left\|U^{i}\right\|_{H^{1}(\partial D_m)}\hspace{.25cm} \left(\mathrm{C}:=C_p\max\limits_{1\leq m \leq M}C_{6m}\right),
\end{eqnarray}
for all $m\in\{1,2,\dots,M\}$.
But,
\begin{eqnarray}\label{normU^i}
\left\|U^{i}\right\|_{H^{1}(\partial D_m)}&\leq&\left\|U^{i}\right\|_{L^{2}(\partial D_m)}+\sum_{l=1}^{2}\left\|\partial_{T_l}U^{i}\right\|_{L^{2}(\partial D_m)}\nonumber\\
                                            &=&\epsilon~\left|\partial B_m\right|^{\frac{1}{2}}+2k\epsilon\left|\partial B_m\right|^{\frac{1}{2}}\quad\left(\mbox{Since },U^{i}(x,\theta)=e^{i\kappa{x}\cdot\theta}\right) \nonumber\\
                                            &\leq&(1+2\kappa)\epsilon\left|\partial \c{B}\right|^{\frac{1}{2}}, \forall m=1,2,\dots,M.
\end{eqnarray}
Now by substituting \eqref{normU^i} in \eqref{nrminvLplusK2}, for each $m=1,\dots,M$, we obtain
\begin{eqnarray}\label{nrmsigmaf}
\left\|\sigma_m\right\|_{L^{2}(\partial D_m)}&\,\leq\,&  \mathcal{C}(\kappa),
\end{eqnarray}
where $\hspace{.25cm}\mathcal{C}(\kappa):=\mathrm{C} \left|\partial \c{B}\right|^{\frac{1}{2}}(1+2\kappa)$. 
\par The condition $\left\|L^{-1}\right\|\left\|K\right\|<1$ is satisfied if 
\begin{eqnarray}\label{Accond-invLK-singl-small}
C_s=\frac{M-1}{4\pi}\left|\partial \c{B}\right|\left(\frac{2\kappa+1}{d}+\frac{2}{d^2}\right)\left(\max\limits_{1\leq m \leq M}C_{6m}\right)\epsilon<1.
\end{eqnarray}
Since $(2\kappa+1)d<\bar{c}$ for $\bar{c}:=(2\kappa_{\max}+1)d_{\max}$, then \eqref{Accond-invLK-singl-small} 
reads as $(M-1)\epsilon<cd^2$, where we set $$c:=\left[\frac{(\bar{c}+2)}{4\pi}|\partial\c{B}|\max\limits_{1\leq m \leq M}C_{6m}\right]^{-1}.$$

\end{proofe}
\subsubsection{Further estimates on the total charge $\int_{ \partial D_m} \sigma_m(s) ds,\,m=1,\dots\,M$}\label{SLPR-3}
\ ~ \
\begin{definition}
\label{Qmdef}
Following \cite{MRAMAG-book1}, we call $\sigma_m\in L^2(\partial D_m)$ used in \eqref{qcimprequiredfrm1}, the solution of the problem (\ref{acimpoenetrable}-\ref{radiationc}), the surface charge distributions. 
For these surface charge distributions, we define the total charge on each surface $\partial D_m$ denoted by $Q_m$ as
\begin{eqnarray}\label{defofQm}
 Q_m:=\int_{ \partial D_m} \sigma_m(s) ds.
\end{eqnarray}
\end{definition}
\begin{lemma}\label{Qmestbigo}
For $m=1,2,\dots,M$, we have
\begin{eqnarray}\label{estofQm}
 |Q_m|&\leq&\tilde{c}\epsilon,
\end{eqnarray}
where $\tilde{c}:=|\partial \c{B}|\mathrm{C} (1+2\kappa)$ with $\partial\c{B}$ and $\mathrm{C}$ are defined in \eqref{fnormS_{ij}21} and \eqref{nrminvLplusK2} respectively.
\end{lemma}
\begin{proofn}{\it{of Lemma \ref{Qmestbigo}}.}
From Proposition \ref{normofsigmastmt}, we know that the surface charge distributions $\sigma_m\in L^2(\partial D_m)$ are bounded by a constant 
$\mathcal{C}(\kappa):=\mathrm{C} \left|\partial \c{B}\right|^{\frac{1}{2}}(1+2\kappa)$. Hence
\begin{eqnarray*}\label{Q2}
 |Q_m|&=&\left|\int_{ \partial D_m} \sigma_m(s) ds\right|\nonumber\\
    &\leq& \| 1  \|_{L^2(\partial D_m)} \| \sigma_m  \|_{L^2(\partial D_m)}\nonumber\\
    &\leq& \| 1  \|_{L^2(\partial D_m)}\mathrm{C} \left|\partial \c{B}\right|^{\frac{1}{2}}(1+2\kappa) \nonumber\\
    &\leq&|\partial \c{B}|\mathrm{C} (1+2\kappa)\epsilon.
\end{eqnarray*}
\end{proofn}
\begin{proposition}\label{farfldthm} 
The far-field pattern $U^\infty$ of the scattered solution of the problem (\ref{acimpoenetrable}-\ref{radiationc}) has the following asymptotic expansion
\begin{equation}\label{x oustdie1 D_m}
U^\infty(\hat{x})=\sum_{m=1}^{M}[e^{-i\kappa\hat{x}\cdot z_{m}}Q_m+O(\kappa\,a^2)],
\end{equation}
with $Q_m$ given by (\ref{defofQm}), if $\kappa\,a<1$ where $O(\kappa\,a^2)\,\leq\,C\kappa\,a^2$ and $C:=\frac{|\partial\c{B}|\mathrm{C} (1+2\kappa)}{\max\limits_{1\leq m \leq M} diam(B_m)}$.

\end{proposition}
\begin{proofn}{\it{of Proposition \ref{farfldthm}}.}
From \eqref{qcimprequiredfrm1}, we have 
\begin{eqnarray*}
 U^{s}(x)&=&\sum_{m=1}^{M}\int_{\partial D_m}\Phi_\kappa(x,s)\sigma_{m} (s)ds,\text{ for }x\in\mathbb{R}^{3}\backslash\left(\mathop{\cup}\limits_{m=1}^M \bar{D}_m\right). \nonumber
\end{eqnarray*}
Hence
\begin{eqnarray}\label{xfarawayimpnt}
U^{\infty}(\hat{x})&=&\sum_{m=1}^{M}\int_{\partial D_m}e^{-i\kappa\hat{x}\cdot s}\sigma_{m} (s)ds\nonumber\\
 &=&\sum_{m=1}^{M}\left(\int_{\partial D_m}e^{-i\kappa\hat{x}\cdot\,z_{m}}\sigma_{m}(s)ds+\int_{\partial D_m}[e^{-i\kappa\hat{x}\cdot\,s}-e^{-i\kappa\hat{x}\cdot\,z_{m}}]\sigma_{m}(s)ds\right)\nonumber\\
 &=&\sum_{m=1}^{M}\left(e^{-i\kappa\hat{x}\cdot\,z_{m}}Q_m+\int_{\partial D_m}[e^{-i\kappa\hat{x}\cdot\,s}-e^{-i\kappa\hat{x}\cdot\,z_{m}}]\sigma_{m}(s)ds\right).
\end{eqnarray}
As in Lemma \ref{Qmestbigo}, we have from Proposition \ref{normofsigmastmt};
\begin{eqnarray}\label{estimationofintsigma}
\int_{\partial D_m}\vert\sigma_{m}(s)\vert ds
&\leq&Ca,
\end{eqnarray}
\begin{equation}
 \text{with}\quad C:=\frac{\tilde{c}}{\max\limits_{1\leq m \leq M} diam(B_m)}=\frac{|\partial \c{B}|\mathrm{C} (1+2\kappa)}{\max\limits_{1\leq m \leq M} diam(B_m)}.
\end{equation}
It gives us the following estimate;
\begin{eqnarray}\label{acestimateforexponentsdif-pre}
\left|\int_{\partial D_m}[e^{-i\kappa\hat{x}\cdot\,s}-e^{-i\kappa\hat{x}\cdot\,z_{m}}]\sigma_{m}(s)ds\right|
                         &\leq&\int_{\partial D_m}\left|e^{-i\kappa\hat{x}\cdot\,s}-e^{-i\kappa\hat{x}\cdot\,z_{m}}\right \|\sigma_{m}(s)|ds\nonumber\\
                         &\leq&\int_{\partial D_m}\sum_{l=1}^{\infty}\kappa^l|s-z_{m}|^l|\sigma_{m}(s)|ds\nonumber\\
                         &\leq&\int_{\partial D_m}\sum_{l=1}^{\infty}\kappa^l\left(\frac{a}{2}\right)^l|\sigma_{m}(s)|ds\nonumber\\
                         &\substack{\leq\\ \eqref{estimationofintsigma}}&Ca\sum_{l=1}^{\infty}\kappa^l\left(\frac{a}{2}\right)^l\nonumber\\
                         &=&\frac{1}{2}C\kappa\,a^2\frac{1}{1-\kappa\frac{a}{2}},\text{ if }
a<\frac{2}{\kappa_{\max}}\left(\leq\frac{2}{\kappa}\right)
\end{eqnarray}
which means
\begin{eqnarray}\label{acestimateforexponentsdif}
\int_{\partial D_m}[e^{-i\kappa\hat{x}\cdot\,s}-e^{-i\kappa\hat{x}\cdot\,z_{m}}]\sigma_{m}(s)ds&\leq&C\kappa\,a^2,\,\text{for}\,a\leq\frac{1}{\kappa_{\max}}.
\end{eqnarray}
Now substitution of \eqref{acestimateforexponentsdif} in \eqref{xfarawayimpnt} gives the required result \eqref{x oustdie1 D_m}.
\end{proofn}

\noindent
Let us derive a formula for $Q_m$. For $s_m\in \partial D_m$, using the Dirichlet boundary condition \eqref{acgoverningsupport}, we have
\begin{eqnarray}\label{Q_mint}
 0&=&U^{t}(s_m)=U^{i}(s_m)+\sum_{j=1}^{M}\int_{\partial D_j}\Phi_\kappa(s_m,s)\sigma_{j} (s)ds\nonumber \\
&=&U^{i}(s_m)+\sum_{\substack{j=1 \\ j\neq m}}^{M}\left(\Phi_\kappa(s_m,z_j)Q_j+\int_{\partial D_j}[\Phi_\kappa(s_m,s)-\Phi_\kappa(s_m,z_j)]\sigma_{j} (s)ds\right)+\int_{\partial D_m}\Phi_\kappa(s_m,s)\sigma_{m} (s)ds.\nonumber\\
\end{eqnarray}
To estimate $\int_{\partial D_j}[\Phi_\kappa(s_m,s)-\Phi_\kappa(s_m,z_j)]\sigma_{j} (s)ds$ for $j\neq\,m$, we write from Taylor series that
\begin{eqnarray}\label{taylorphifar}
\Phi_\kappa(s_m,s)-\Phi_\kappa(s_m,z_j)=(s-z_j)\cdot R(s_m,s),&\ &R(s_m,s)=\int_{0}^{1}\nabla_2\Phi_\kappa(s_m,s-\alpha(s-z_j))d\alpha.
\end{eqnarray}
\begin{itemize}
\item From the definition of $\Phi_\kappa(x,y)
$ in \eqref{definition-ac-small-fundamentalkappa},
we have $\nabla_y\Phi_\kappa(x,y)=\Phi_\kappa(x,y)\left[\frac{1}{|x-y|}-i\kappa\right]\frac{x-y}{|x-y|}$ 
and hence, for $s\in \bar{D}_j$, we obtain
\begin{eqnarray}\label{estimatofRxs}
 |R(s_m,s)|\leq\max\limits_{y\in \bar{D}_j}\left|\nabla_y\Phi_\kappa(s_m,y)\right|<\,\frac{1}{d}\left(\kappa+\frac{1}{d}\right).
\end{eqnarray}
\end{itemize}
For $m,j=1,\dots,M$, and $j\neq\,m$, by making use of \eqref{estimatofRxs} and \eqref{estimationofintsigma} we obtain the estimate below;
\begin{eqnarray}\label{estphismzjdif}
 \left|\int_{\partial D_j}[\Phi_\kappa(s_m,s)-\Phi_\kappa(s_m,z_j)]\sigma_{j} (s)ds\right|&=&\left|\int_{\partial D_j}(s-z_j)\cdot R(s_m,s)\sigma_{j} (s)ds\right|\nonumber\\
                                                               &\leq&\int_{\partial D_j}\left|s-z_j\right| \left|R(s_m,s)\right| \left|\sigma_{j} (s)\right|ds\nonumber\\
                                                               &<\,&\frac{a}{d}\left(\kappa+\frac{1}{d}\right)\int_{\partial D_j}\left|\sigma_{j} (s)\right|ds\nonumber\\
                                                               &<\,&C\frac{a}{d}\left(\kappa+\frac{1}{d}\right)a.
\end{eqnarray}
\noindent
Then \eqref{Q_mint} can be written as
\begin{equation}\label{qcimsurfacefrm}
\begin{split}
 0=U^{i}(s_m)+\sum_{\substack{j=1 \\ j\neq m}}^{M}\Phi_\kappa(s_m,z_j)Q_j&+O\left((M-1)\left(\frac{\kappa a^2}{d}+\frac{a^2}{d^2}\right)\right) \\
 &+\int_{\partial D_m}\Phi_0(s_m,s)\left[1+(e^{i\kappa|s_m-s|}-1)\right]\sigma_{m}(s)ds.
\end{split}
\end{equation}
By using the Taylor series expansions of the exponential term $e^{i\kappa|s_m-s|}$, the above can also be written as,
\begin{eqnarray}\label{qcimsurfacefrm1}
\int_{\partial D_m}\Phi_0(s_m,s)\sigma_{m} (s)ds+O(\kappa a)&=&-U^{i}(s_m)-\sum_{\substack{j=1 \\ j\neq m}}^{M}\Phi_\kappa(s_m,z_j)Q_j+O\left((M-1)\left(\frac{\kappa a^2}{d}+\frac{a^2}{d^2}\right)\right).
\end{eqnarray}
\begin{indeed}
\begin{itemize}
\item For $m=1,\dots,M$, we have
\begin{eqnarray}\label{phinotexp-1est}
\left|\int_{\partial D_m}\Phi_0(s_m,s)\left(e^{i\kappa|s_m-s|}-1\right)\sigma_{m}(s)ds\right|
                             &=&\left|\int_{\partial D_m}\frac{1}{4\pi|s_m-s|}\left(\sum\limits_{l=1}^\infty\frac{(i\kappa)^l}{l!}|s_m-s|^{l}\right)\sigma_{m}(s)ds\right|\nonumber\\
                             &=&\left|\frac{1}{4\pi}\int_{\partial D_m}\left(\sum\limits_{l=1}^\infty\frac{(i\kappa)^l}{l!}|s_m-s|^{l-1}\right)\sigma_{m}(s)ds\right|\nonumber\\
                             &<&\frac{1}{2}\sum\limits_{l=1}^\infty\frac{\kappa^l}{l!}a^{l-1}\int_{\partial D_m}\left|\sigma_{m}(s)\right|ds\nonumber\\
                             &\substack{\leq\\\eqref{estimationofintsigma}}&\frac{1}{2}\sum\limits_{l=1}^\infty\frac{\kappa^l}{l!}a^{l-1}\cdot Ca\nonumber\\
                             &\leq&C\sum\limits_{l=1}^\infty\frac{\kappa^l}{2^l}a^{l}\nonumber\\
                             &\leq&C \kappa\,a, \hspace{.25cm} \mbox{for } a\leq\frac{1}{\kappa_{\max}}.
\end{eqnarray}
\end{itemize}
\end{indeed}
\noindent
Define $U_m:=\int_{\partial D_m}\Phi_0(s_m,s)\sigma_{m} (s)ds, s_m\in \partial D_m$. Then \eqref{qcimsurfacefrm1} can be written as
\begin{eqnarray}\label{qcimsurfacefrm2}
U_m&=&-U^{i}(s_m)-\sum_{\substack{j=1 \\ j\neq m}}^{M} \Phi_\kappa(s_m,z_j)Q_j+ O(\kappa a)+O\left((M-1)\left(\frac{\kappa a^2}{d}+\frac{a^2}{d^2}\right)\right).
\end{eqnarray}
For $m=1,\dots,M$, we set 
\begin{eqnarray}\label{def-acoust-singl-revise-defntn-barum}
\bar{U}_m:=-U^{i}(z_m)-\sum\limits_{\substack{j=1 \\ j\neq m}}^{M}\Phi_\kappa(z_m,z_j)Q_j.
\end{eqnarray}
Let $\bar{\sigma}_m\in L^2(\partial D_m)$ be the corresponding surface charge distributions, i.e.,
\begin{equation}\label{barqcimsurfacefrm1}
\int_{\partial D_m}\Phi_0(s_m,s)\bar{\sigma}_{m} (s)ds=\bar{U}_m, s_m\in\,\partial D_m.
\end{equation}
The total charge on the surface $\partial D_m$ is given by 
$$\bar{Q}_m:=\int_{\partial D_m}\bar{\sigma}_m(s)ds.$$
Now, we set the electrical capacitance $\bar{C}_m$ for $1\leq m\leq M$ as 
$$\bar{C}_m:=\frac{\bar{Q}_m}{\bar{U}_m}.$$ 
\begin{lemma}\label{lemmadifssbQQbCCb}
 We have the following estimates
 \begin{eqnarray}
  \|\sigma_{m}-\bar{\sigma}_{m}\|_{H^{-1}(\partial D_m)}&=&O\left(\kappa a+(M-1)\left(\frac{\kappa a^2}{d}+\frac{a^2}{d^2}\right)\right),\label{difssb}\\
  Q_m-\bar{Q}_m&=&O\left(\kappa a^2+(M-1)\left(\frac{\kappa a^3}{d}+\frac{a^3}{d^2}\right)\right),\label{difQmd}
 \end{eqnarray}
where the constants appearing in $O(.)$ depend only on the Lipschitz character of $B_m$.
\end{lemma}
\begin{proofn}{\it{of Lemma \ref{lemmadifssbQQbCCb}}.}
By taking the difference between \eqref{qcimsurfacefrm2} and \eqref{barqcimsurfacefrm1}, we obtain
\begin{eqnarray}\label{dbarqcimsurfacefrm1}
\hspace{-.1cm}U_m-\bar{U}_m=\int_{\partial D_m}\Phi_0(s_m,s)\left(\sigma_{m}-\bar{\sigma}_{m}\right) (s)ds=O(\kappa a)+O\left((M-1)\left(\frac{\kappa a^2}{d}+\frac{a^2}{d^2}\right)\right),~ s_m\in \partial D_m.
\end{eqnarray} 
\begin{indeed} by using Taylor series, we have
\begin{itemize}
 \item $U^{i}(s_m)-U^{i}(z_m)=O(\kappa a)$.
  \item $\Phi_\kappa(s_m,z_j)-\Phi_\kappa(z_m,z_j)=O\left(\frac{\kappa a}{d}+\frac{a}{d^2}\right)$ and the estimate of ${Q}_j$ given in \eqref{estofQm}.
\end{itemize}
\end{indeed}
In operator form we can write \eqref{dbarqcimsurfacefrm1}, using \eqref{def-of_Sik}, as 
\begin{eqnarray*}
(\mathcal{S}^{i_\kappa}_{D_m})^{*}\left(\sigma_{m}\hspace{-.07cm}-\hspace{-.07cm}\bar{\sigma}_{m}\right) (s_m):=\int_{\partial D_m}\hspace{-.2cm}\Phi_0(s_m,s)\left(\sigma_{m}\hspace{-.07cm}-\hspace{-.07cm}\bar{\sigma}_{m}\right) (s)ds=
O(\kappa a)+O\left((M\hspace{-.07cm}-\hspace{-.07cm}1)\left(\frac{\kappa a^2}{d}\hspace{-.07cm}+\hspace{-.07cm}\frac{a^2}{d^2}\right)\right),~ s_m\in \partial D_m.\nonumber
\end{eqnarray*} 
 Here, $(\mathcal{S}^{i_\kappa}_{D_m})^{*}:H^{-1}(\partial D_m)\rightarrow L^2(\partial D_m)$ is the adjoint of $\mathcal{S}^{i_\kappa}_{D_m}:L^{2}(\partial D_m)\rightarrow H^{1}(\partial D_m)$.  We know that, 
\begin{eqnarray*}
\left\|(\mathcal{S}^{i_\kappa}_{D_m})^{*}\right\|_{\mathcal{L}\left(H^{-1}(\partial D_m),L^2(\partial D_m)\right)}=
 \left\|\mathcal{S}^{i_\kappa}_{D_m}\right\|_{\mathcal{L}\left(L^2(\partial D_m),H^1(\partial D_m)\right)}
\end{eqnarray*}
and 
\begin{eqnarray*}
\left\|{((\mathcal{S}^{i_\kappa}_{D_m})^{*})}^{-1}\right\|_{\mathcal{L}\left(L^2(\partial D_m),H^{-1}(\partial D_m)\right)}=
 \left\|(\mathcal{S}^{i_\kappa}_{D_m})^{-1}\right\|_{\mathcal{L}\left(H^{1}(\partial D_m),L^2(\partial D_m)\right)},\nonumber
\end{eqnarray*}
 then from \eqref{optnormsik} of Lemma \ref{lemmanrm1singulayer31}, we obtain 
$\left\|{((\mathcal{S}^{i_\kappa}_{D_m})^{*})}^{-1}\right\|_{\mathcal{L}\left(L^2(\partial D_m),H^{-1}(\partial D_m)\right)}=O(a^{-1})$. 
Hence, we get the required results in the following manner. First, we have
\begin{eqnarray*}
 \|\sigma_{m}-\bar{\sigma}_{m}\|_{H^{-1}(\partial D_m)}&\leq&\left\|{((\mathcal{S}^{i_\kappa}_{D_m})^{*})}^{-1}\right\|_{\mathcal{L}\left(L^2(\partial D_m),H^{-1}(\partial D_m)\right)}
                                                    \left\| O(\kappa a)+O\left((M-1)\left(\frac{\kappa a^2}{d}+\frac{a^2}{d^2}\right)\right)\right\|_{L^2(\partial D_m)}\nonumber\\
                                               &=&O\left(\kappa a+(M-1)\left(\frac{\kappa a^2}{d}+\frac{a^2}{d^2}\right)\right).\nonumber
\end{eqnarray*}
and second, we have
\begin{eqnarray*}
 |Q_m-\bar{Q}_m|&=&\left|\int_{\partial D_m}\left(\sigma_{m}-\bar{\sigma}_{m}\right) (s)ds\right|\nonumber\\
                &\leq& \|\sigma_{m}-\bar{\sigma}_{m}\|_{H^{-1}(\partial D_m)} \|1\|_{H^1(\partial D_m)}\nonumber\\
                &=& O\left(\kappa a^2+(M-1)\left(\frac{\kappa a^3}{d}+\frac{a^3}{d^2}\right)\right).
\end{eqnarray*}
\end{proofn}
\begin{lemma}\label{lemmadifssbQQbCCb1}
For every $m,\,1\leq m\leq M$, the capacitance $\bar{C}_m$  and the charge $\bar{Q}_m$ are of the form;
\begin{eqnarray}\label{asymptotCap}
\bar{C}_m\,=\,\frac{\bar{C}_{B_m}}{\max\limits_{1\leq m \leq M} diam(B_m)}a
 & \mbox{ and }&
\bar{Q}_m\,=\,\frac{\bar{Q}_{B_m}}{\max\limits_{1\leq m \leq M} diam(B_m)}a,
 \end{eqnarray}
where $\bar{C}_{B_m}$ and $\bar{Q}_{B_m}$ are the capacitance and the charge of $B_m$ respectively.
\end{lemma}
\begin{proofn}{\it{of Lemma \ref{lemmadifssbQQbCCb1}}.} 
Take $0<\epsilon\leq 1$, $z\in \mathbb{R}^3$ and 
write, $D_\epsilon:=\epsilon B+z\subset \mathbb{R}^3$. 
For $\psi_\epsilon\in L^2(\partial D_\epsilon)$ and $\psi\in L^2(\partial B)$, the operators $\mathcal{S}^{i_{\kappa}}_{D_\epsilon}:L^{2}(\partial D_\epsilon)\rightarrow H^{1}(\partial D_\epsilon)$ 
and $\mathcal{S}^{i_{\kappa}}_B:L^{2}(\partial B)\rightarrow H^{1}(\partial B)$ are defined in \eqref{def-of_Sik}.
These two operators define the corresponding potentials $\bar{U}_\epsilon$, $\bar{U}_B$ on the surfaces $\partial D_\epsilon$ and $\partial B$ with respect to the 
surface charge distributions $\psi_\epsilon$ and $\psi$ respectively. Let these potentials be equal to some constant $D$. 
 Let also the total charge of these conductors $D_\epsilon$, $B$ be $\bar{Q}_\epsilon$ and $\bar{Q}_B$, and the capacitances be $\bar{C}_\epsilon$ and $\bar{C}_B$ respectively.
Then we can write
\begin{eqnarray*}
 \bar{U}_\epsilon:=\mathcal{S}^{i_{\kappa}}_{D_\epsilon} \psi_\epsilon(x)=D,
&\,&\bar{U}_B:=\mathcal{S}^{i_{\kappa}}_B \psi(\xi)=D,
\hspace{.10cm} \forall x\in\partial D_\epsilon,\forall \xi\in\partial B.
\end{eqnarray*}
We have by definitions,
$
\bar{Q}_\epsilon=\int_{\partial D_\epsilon} \psi_\epsilon(y) dy,\,\bar{Q}_B=\int_{\partial B} \psi(\eta) d\eta,\text{ and }
\bar{C}_\epsilon=\frac{\bar{Q}_\epsilon}{\bar{U}_\epsilon},\, \bar{C}_B=\frac{\bar{Q}_B}{\bar{U}_B}.
$
\newline
Observe that,
\begin{align*}
D=\mathcal{S}^{i_{\kappa}}_{D_\epsilon} \psi_\epsilon(x)=&\int_{\partial D_\epsilon} \frac{1}{4\pi|x-y|}\psi_\epsilon(y) dy
       & D=\mathcal{S}^{i_{\kappa}}_B \psi(\xi)=&\int_{\partial B} \frac{1}{4\pi|\xi-\eta|}\psi(\eta) d\eta  
\nonumber\\
                                        =&\int_{\partial B} \frac{1}{4\pi\epsilon|\xi-\eta|}\psi_\epsilon(\epsilon\eta+z)\epsilon^2 d\eta
       &                                =&\int_{\partial D_\epsilon} \frac{1}{\frac{4\pi}{\epsilon}|x-y|}\psi\left(\frac{y-z}{\epsilon}\right)\epsilon^{-2} dy
\nonumber\\
                                        =&\,\epsilon\int_{\partial B} \frac{1}{4\pi|\xi-\eta|}\hat{\psi}_\epsilon(\eta) d\eta
       &                                =&\,\epsilon^{-1}\int_{D_\epsilon} \frac{1}{4\pi|x-y|}\check{\psi}(y) d\eta
\nonumber\\
                                        =&\,\epsilon\,\mathcal{S}^{i_{\kappa}}_B \hat{\psi}_\epsilon(\xi),\hspace{.25cm} \left[\hat{\psi}_\epsilon(\eta):=\psi_\epsilon(\epsilon\eta+z)\right],
      &                                 =&\,\epsilon^{-1}\mathcal{S}^{i_{\kappa}}_{D_\epsilon} \check{\psi}(x),\hspace{.25cm} \left[\check{\psi}(y):=\psi\left(\frac{y-z}{\epsilon}\right)\right].
\end{align*}
Hence, 
$
 \psi_\epsilon=\frac{1}{\epsilon}\check{\psi}
\mbox{ and }  
\psi=\epsilon\hat{\psi}_\epsilon.
$
Now we have,
\begin{eqnarray*}
 \bar{Q}_\epsilon&=&\int_{\partial D_\epsilon} \psi_\epsilon(y) dy
                 =\int_{\partial D_\epsilon}\frac{1}{\epsilon}\check{\psi}(y) dy,\nonumber\\
                 &=&\int_{\partial B}\frac{1}{\epsilon}\check{\psi}(\epsilon\eta+z) \epsilon^2 d\eta
                 =\epsilon\int_{\partial B}\check{\psi}(\epsilon\eta+z)d\eta,\nonumber\\
                 &=&\epsilon\int_{\partial B}\hat{\check{\psi}}(\eta) d\eta
                 =\epsilon\int_{\partial B}\psi(\eta) d\eta,\nonumber\\
                 &=&\epsilon\,\bar{Q}_B\nonumber\\
\end{eqnarray*}
which gives us,
\begin{eqnarray*}       
 \bar{C}_\epsilon&=&\frac{\bar{Q}_\epsilon}{\bar{U}_\epsilon}\,=\,\frac{\epsilon\bar{Q}_B}{D}
                 \,=\,\epsilon\frac{\bar{Q}_B}{D}\,=\,\epsilon\frac{\bar{Q}_B}{\bar{U}_B}
                 \,=\,\epsilon\bar{C}_B.        
\end{eqnarray*}
As we have $D_m=\epsilon B_m+z_m$ and $a=\max\limits_{1\leq m \leq M} diam D_m=\epsilon\max\limits_{1\leq m \leq M} diam(B_m)$, we obtain 
\begin{eqnarray*}
\bar{Q}_m\,=\,\epsilon\bar{Q}_{B_m}\,=\,\frac{\bar{Q}_{B_m}}{\max\limits_{1\leq m \leq M} diam(B_m)}a & \mbox{ and }&
\bar{C}_m\,=\,\epsilon\bar{C}_{B_m}\,=\,\frac{\bar{C}_{B_m}}{\max\limits_{1\leq m \leq M} diam(B_m)}a.
\end{eqnarray*}
\end{proofn}
\begin{proposition} \label{fracqfracc-ac}For $m=1,2,\dots,M$, the total charge $\bar{Q}_m$ on each surface $\partial D_m$ of the small scatterer $D_m$ can be calculated from the algebraic system 
  \begin{eqnarray}\label{fracqfrac}
 \frac{\bar{Q}_m}{\bar{C}_m} &=&-U^{i}(z_m)-\sum_{\substack{j=1 \\ j\neq m}}^{M}\bar{C}_j \Phi_\kappa(z_m,z_j)\frac{\bar{Q}_j}{\bar{C}_j}+ 
O\left((M-1)\frac{\kappa a^2}{d}+(M-1)^2\left(\frac{\kappa a^3}{d^2}+\frac{a^3}{d^3}\right)\right).
  \end{eqnarray}
\end{proposition}
\begin{proofn}{\it{of Proposition  \ref{fracqfracc-ac}}.}
We can rewrite \eqref{def-acoust-singl-revise-defntn-barum} as
\begin{eqnarray*}
 \frac{\bar{Q}_m}{\bar{C}_m} &=&-U^{i}(z_m)-\sum_{\substack{j=1 \\ j\neq m}}^{M} \Phi_\kappa(z_m,z_j)Q_j\nonumber\\
                             &=&-U^{i}(z_m)-\sum_{\substack{j=1 \\ j\neq m}}^{M} \Phi_\kappa(z_m,z_j)\bar{Q}_j-\sum_{\substack{j=1 \\ j\neq m}}^{M} \Phi_\kappa(z_m,z_j)(Q_j-\bar{Q}_j)\nonumber\\
                             &=&-U^{i}(z_m)-\sum_{\substack{j=1 \\ j\neq m}}^{M} \Phi_\kappa(z_m,z_j)\bar{Q}_j+O\left((M-1)\frac{\kappa a^2}{d}+(M-1)^2\left(\frac{\kappa a^3}{d^2}+\frac{a^3}{d^3}\right)\right),\nonumber\\
\end{eqnarray*}

where we used \eqref{difQmd} and the fact $\Phi_\kappa(z_m,z_j)=O\left(\frac{ 1}{d}\right)$.
\end{proofn}
\subsection{Using DLP representation}\label{DLPR}
\subsubsection{The representation via DLP}\label{DLPR-1}
 \ ~ \ 
\begin{proposition} \label{existence-of-sigmasdbl}
For $m=1,2,\dots,M$, there exists $\sigma_m\in H^{r}(\partial D_m),\,r\in[0,\,1]$ such that the solution of the problem (\ref{acimpoenetrable}-\ref{radiationc}) is of the form
 \begin{equation}\label{qcimprequiredfrm1dbl}
  U^{t}(x)=U^{i}(x)+\sum_{m=1}^{M}\int_{\partial D_m}\frac{\partial\Phi_\kappa(x,s)}{\partial \nu_m(s)}\sigma_{m} (s)ds,~x\in\mathbb{R}^{3}\backslash\left(\mathop{\cup}_{m=1}^M \bar{D}_m\right), 
\end{equation}
with $\nu_m$ denoting the outward unit normal vector of $\partial D_m$.
\end{proposition}
\begin{proofn}{\it{of Proposition  \ref{existence-of-sigmasdbl}}.}
  We look for the solution of the problem (\ref{acimpoenetrable}-\ref{radiationc}) of the form \eqref{qcimprequiredfrm1dbl}, 
  then from the Dirichlet boundary condition \eqref{acgoverningsupport}, we obtain
\begin{equation}\label{qcimprequiredfrm1bddbl}
\frac{\sigma_{j} (s_j)}{2}+\int_{\partial D_j}\frac{\partial\Phi_\kappa(s_j,s)}{\partial \nu_j(s)}\sigma_{j} (s)ds+\sum_{\substack{m=1\\m\neq j}}^{M}\int_{\partial D_m}\frac{\partial\Phi_\kappa(s_j,s)}{\partial \nu_m(s)}\sigma_{m} (s)ds=-U^i{(s_j)},\,\forall s_j\in \partial D_j,\, j=1,\dots,M.
\end{equation}
One can write it in compact form as 
$(\frac{1}{2}\textbf{I}+DL+DK)\sigma=-U^{In}$ with $\textbf{I}:=(\textbf{I}_{mj})_{m,j=1}^{M}$, $DL:=(DL_{mj})_{m,j=1}^{M}$ and $DK:=(DK_{mj})_{m,j=1}^{M}$, where
\begin{eqnarray}\label{definition-DL_DK}
\hspace{-.3cm}\textbf{I}_{mj}=\left\{\begin{array}{ccc}
            I,\, \text{Identity operator} & m=j\\
            0,\, \text{zero~operator}~~~~~~ & else
           \end{array}\right.,\,
&
DL_{mj}=\left\{\begin{array}{ccc}
            \mathcal{D}_{mj} & m=j\\
            0 & else
           \end{array}\right.,\,
&
DK_{mj}=\left\{\begin{array}{ccc}
            \mathcal{D}_{mj} & m\neq j\\
            0 & else
           \end{array}\right., 
\end{eqnarray}
%
$U^{In}=U^{In}(s_1,\dots,s_M):=\left(U^i(s_1),\dots,U^i(s_M)\right)^T$ and $\sigma=\sigma(s_1,\dots,s_M):=\left(\sigma_1(s_1),\dots,\sigma_M(s_M)\right)^T$. 
Here, for the indices  $m$ and $j$ fixed, $\mathcal{D}_{mj}$ is the integral operator acting as
\begin{eqnarray}\label{defofDmjed}
 \mathcal{D}_{mj}(\sigma_j)(t):=\int_{\partial D_j}\frac{\partial\Phi_\kappa(t,s)}{\partial \nu_j(s)}\sigma_j(s)ds,\quad t\in\partial D_m. 
\end{eqnarray}
The operator $\frac{1}{2}I+\mathcal{D}_{mm}:H^{r}(\partial D_m)\rightarrow H^{r}(\partial D_m)$ 
is Fredholm with zero index and for $m\neq j$, $\mathcal{D}_{mj}:H^{r}(\partial D_j)\rightarrow H^{r}(\partial D_m)$ 
is compact for $0\leq r\leq 1$, when $\partial D_m$ has a Lipschitz regularity, see \cite{MD:InteqnsaOpethe1997}.\footnote{In \cite{MD:InteqnsaOpethe1997}, this property is proved for the case $\kappa=0$. 
As in Section \ref{SLPR-1}, concerning the single layer operator, by a perturbation argument, we have the same results for every $\kappa$ in $[0,\,\kappa_{\max}]$, assuming that 
$a<\frac{1}{\kappa_{\max}}\sqrt[3]{\frac{4\pi}{3}}{\rm j}_{1/2,1}$.} Remark here that, for the scattering by single obstacle $DK$ is zero  operator.
\par
So, $(\frac{1}{2}\textbf{I}+DL+DK):\prod\limits_{m=1}^{M}H^{r}(\partial D_m)\rightarrow \prod\limits_{m=1}^{M}H^{r}(\partial D_m)$ is Fredholm with zero index.  
We induce the product of spaces by the maximum of the norms of the space.
To show that $(\frac{1}{2}\textbf{I}+DL+DK)$ is invertible it is enough to show that it is injective. i.e. $(\frac{1}{2}\textbf{I}+DL+DK)\sigma=0$ implies $\sigma=0$.
\newline
Write, 
$$\tilde{U}(x)=\sum_{m=1}^{M}\int_{\partial D_m}\frac{\partial\Phi_\kappa(x,s)}{\partial \nu_m(s)}\sigma_m(s)ds, \mbox{ in } \mathbb{R}^3\backslash \left(\mathop{\cup}_{m=1}^{M}\bar{D}_m\right)$$
and
$$\tilde{\tilde{U}}(x)=\sum_{m=1}^{M}\int_{\partial D_m}\frac{\partial\Phi_\kappa(x,s)}{\partial \nu_m(s)}\sigma_m(s)ds, \mbox{ in } \mathop{\cup}_{m=1}^{M}D_m.$$

Then $\tilde{U}$ satisfies $\Delta\tilde{U}+\kappa^2\tilde{U}=0$ for $x\in\mathbb{R}^{3}\backslash\left(\mathop{\cup}\limits_{m=1}^{M}\bar{D}_m\right)$,  
with S.R.C and $\tilde{U}(x)=0$ on $\mathop{\cup}\limits_{m=1}^{M}\partial D_m$.  
Similarly, $\tilde{\tilde{U}}$ satisfies $\Delta\tilde{\tilde{U}}+\kappa^2\tilde{\tilde{U}}=0$ for $x\in\mathop{\cup}\limits_{m=1}^{M}D_m$ with  $\tilde{\tilde{U}}(x)=0$ 
on $\mathop{\cup}\limits_{m=1}^{M}\partial D_m$.
By the jump relations, we have
\begin{eqnarray}\label{nbc1dbl}
\tilde{U}(s)=0&\Longrightarrow&(\mathbf{K}\sigma_m)(s)+\frac{\sigma_m(s)}{2}+\sum\limits^{M}_{\substack{j=1\\j\neq m}}\mathcal{D}_{mj}(\sigma_j)(s)=0
\end{eqnarray} 
and 
\begin{eqnarray}\label{nbc2dbl}
\tilde{\tilde{U}}(s)=0&\Longrightarrow&(\mathbf{K}\sigma_m)(s)-\frac{\sigma_m(s)}{2}+\sum\limits^{M}_{\substack{j=1\\j\neq m}}\mathcal{D}_{mj}(\sigma_j)(s)=0
\end{eqnarray}
 for $s\in \partial D_m$ and for $m=1,\dots,M$. Here we recall that
\[ (\mathbf{K}\sigma_m)(s)=\int_{\partial D_m}\frac{\partial\Phi_\kappa(s,t)}{\partial \nu_m(t)} \sigma_m(t) dt,\mbox{ for } m=1,\dots,M.\]
Difference between \eqref{nbc1dbl} and \eqref{nbc2dbl} provide us with $\sigma_m=0$ for all $m$. \\
We conclude then that $\frac{1}{2}\textbf{I}+DL+DK=:\frac{1}{2}\textbf{I}+\mathcal{D}:\prod\limits_{m=1}^{M}H^{r}(\partial D_m)\rightarrow \prod\limits_{m=1}^{M}H^{r}(\partial D_m)$ is invertible.
\end{proofn}
\subsubsection{An appropriate estimate of the densities $\sigma_m,\,m=1,\dots,M$}\label{DLPR-2}  
From the above theorem, we have the following representation of $\sigma$:
\begin{eqnarray}\label{invDLplusDK}
 \sigma&=&(\frac{1}{2}\textbf{I}+DL+DK)^{-1}U^{In} \nonumber\\
       &=&(\frac{1}{2}\textbf{I}+DL)^{-1}(\textbf{I}+(\frac{1}{2}\textbf{I}+DL)^{-1}DK)^{-1}U^{In} \nonumber\\
       &=&(\frac{1}{2}\textbf{I}+DL)^{-1}\sum_{l=0}^{\infty}\left(-(\frac{1}{2}\textbf{I}+DL)^{-1}DK\right)^{l}U^{In}, \hspace{.5cm} \mbox{if } \left\|(\frac{1}{2}\textbf{I}+DL)^{-1}DK\right\|<1.
\end{eqnarray}
\noindent
The operator $\frac{1}{2}\textbf{I}+DL$ is invertible since it is Fredholm of index zero and injective. 
This implies that
\begin{eqnarray}\label{nrminvDLplusDK}
 \left\|\sigma\right\| 
                                  &\leq&\frac{\left\|(\frac{1}{2}\textbf{I}+DL)^{-1}\right\|}{1-\left\|(\frac{1}{2}\textbf{I}+DL)^{-1}\right\|\left\|DK\right\|}\left\|U^{In}\right\|.
\end{eqnarray}
We use the following notations
\begin{eqnarray}
\left\|DK\right\|&:= &\left\|DK\right\|_{\mathcal{L}\left(\prod\limits_{m=1}^{M}L^{2}(\partial D_m),\prod\limits_{m=1}^{M}L^{2}(\partial D_m)\right)}\nonumber\\
    &\equiv&
 \max\limits_{1\leq m \leq M}\sum_{j=1}^{M}\left\|DK_{mj}\right\|_{\mathcal{L}\left(L^{2}(\partial D_j),L^{2}(\partial D_m)\right)}\nonumber\\
    &=&\max\limits_{1\leq m \leq M}\sum_{\substack{j=1\\j\neq\,m}}^{M}\left\|\mathcal{D}_{mj}\right\|_{\mathcal{L}\left(L^{2}(\partial D_j),L^{2}(\partial D_m)\right)},\label{DKnrm}\\
\left\|(\frac{1}{2}\textbf{I}+DL)^{-1}\right\|
&:=& \left\|(\frac{1}{2}\textbf{I}+DL)^{-1}\right\|_{\mathcal{L}\left(\prod\limits_{m=1}^{M}L^{2}(\partial D_m),\prod\limits_{m=1}^{M}L^{2}(\partial D_m)\right)}\nonumber\\
     &\equiv&\max\limits_{1\leq m \leq M}\sum_{j=1}^{M}\left\|{\left((\frac{1}{2}\textbf{I}+DL)^{-1}\right)_{mj}}\right\|_{\mathcal{L}\left(L^{2}(\partial D_m),L^{2}(\partial D_j)\right)}\nonumber\\
     &=&
\max\limits_{1\leq m \leq M}\left\|(\frac{1}{2}I+\mathcal{D}_{mm})^{-1}\right\|_{\mathcal{L}\left(L^{2}(\partial D_m),L^{2}(\partial D_m)\right)},\label{invDLnrm}\\
\left\|\sigma\right\|&:=& \left\|\sigma\right\|_{\prod\limits_{m=1}^{M}L^{2}(\partial D_m)}
    \,\quad\equiv\,\max\limits_{1\leq m \leq M}\left\|\sigma_{m}\right\|_{L^{2}(\partial D_m)}\label{sigmaU^{In}nrmdbl1}\\
\,\mbox{and}\qquad
\left\|U^{In}\right\|&:= &\left\|U^{In}\right\|_{\prod\limits_{m=1}^{M}L^{2}(\partial D_m)}
    \equiv\,\max\limits_{1\leq m \leq M}\left\|U^{i}\right\|_{L^{2}(\partial D_m)}. \label{sigmaU^{In}nrmdbl}
\end{eqnarray}
In the following proposition, we provide conditions under which $\left\|L^{-1}\right\|\left\|K\right\|<1$ and then estimate $\left\|\sigma\right\|$ via \eqref{nrminvDLplusDK}.

\begin{proposition}\label{normofsigmastmtdbl}
There exists a constant $\grave{c}$ depending only on the 
Lipschitz character of $B_m,m=1,\dots,M$, $d_{\max}$
and $\kappa_{\max}$ such that 
if $\sqrt{M-1}\epsilon<\grave{c}d$, 
 then the $L^2$-norm of densities $\sigma_m$ appearing in the solution \eqref{qcimprequiredfrm1dbl} 
of the problem (\ref{acimpoenetrable}-\ref{radiationc}) are bounded by a uniform constant times $\epsilon$.
 \end{proposition}
 \bigskip
\noindent
\ ~ \ \\
Here as well, we divide the proof of Proposition \ref{normofsigmastmtdbl} into two steps. In the first step, we assume we have a single obstacle and then in the second step we deal with the multiple obstacle case.
\paragraph{The case of a single obstacle}\label{singleobscasedbl}
Let us consider a single obstacle $D_\epsilon:=\epsilon B+z$. 
Then define the operator 
$\mathcal{D}_{ D_\epsilon}:L^2(\partial D_{\epsilon})\rightarrow L^2(\partial D_{\epsilon})$ by
\begin{eqnarray}\label{defofDpartialDe}
\left(\mathcal{D}_{ D_\epsilon}\psi\right) (s):=\int_{\partial D_\epsilon}\frac{\partial\Phi_\kappa(s,t)}{\partial\nu(t)}\psi(t)dt.
\end{eqnarray}
Following the arguments in the proof of Proposition  \ref{existence-of-sigmasdbl}, the integral operator 
$\frac{1}{2}I+\mathcal{D}_{ D_\epsilon}:L^2(\partial D_{\epsilon})\rightarrow L^2(\partial D_{\epsilon})$ is invertible. 
If we consider the problem (\ref{acimpoenetrable}-\ref{radiationc}) 
in $\mathbb{R}^{3}\backslash \bar{D}_\epsilon$, we obtain
  $$\sigma=(\frac{1}{2}I+\mathcal{D}_{ D_\epsilon})^{-1}U^{i},\,\text{where}\,DL=\mathcal{D}_{ D_\epsilon}$$ 
and then 
\begin{equation}\label{estsigm1dbl}
 \|\sigma \|_{L^2(\partial D_{\epsilon})}\leq  \|(\frac{1}{2}I+\mathcal{D}_{ D_\epsilon})^{-1} \|_{\mathcal{L}\left(L^2(\partial D_{\epsilon}),L^2(\partial D_{\epsilon})\right)} \|U^i \|_{L^2(\partial D_{\epsilon})}.
\end{equation}
\begin{lemma}\label{rep1dbllayer}
 Let $\phi,\psi \in L^{2}(\partial D_\epsilon)$. Then, 
\begin{equation}\label{rep1dbllayer0}
 \mathcal{D}_{ D_\epsilon}\psi= (\mathcal{D}^\epsilon_B \hat{\psi})^\vee,
\end{equation}
\begin{equation}\label{rep1dbllayer1}
 \left(\frac{1}{2}I+\mathcal{D}_{ D_\epsilon}\right)\psi=  \left(\left(\frac{1}{2}I+\mathcal{D}^\epsilon_B \right)\hat{\psi}\right)^\vee,
\end{equation}
 \begin{equation}\label{rep1dbllayer2}
 {\left(\frac{1}{2}I+\mathcal{D}_{ D_\epsilon}\right)}^{-1}\phi= \left({\left(\frac{1}{2}I+\mathcal{D}^\epsilon_B\right)}^{-1} \hat{\phi}\right)^\vee
\end{equation}
\begin{equation}\label{nrm1dbllayer2}
 \left\|{\left(\frac{1}{2}I+\mathcal{D}_{ D_\epsilon}\right)}^{-1}\right\|_{\mathcal{L}\left(L^2(\partial D_\epsilon), L^2(\partial D_\epsilon) \right)}=
\left\|{\left(\frac{1}{2}I+\mathcal{D}^\epsilon_B\right)}^{-1}\right\|_{\mathcal{L}\left(L^2(\partial B), L^2(\partial B) \right)}
\end{equation}
and
\begin{equation}\label{nrm1dbllayer2-1}
 \left\|{\left(\frac{1}{2}I+\mathcal{D}_{ D_\epsilon}\right)}^{-1}\right\|_{\mathcal{L}\left(H^1(\partial D_\epsilon), H^1(\partial D_\epsilon) \right)}\leq 
\epsilon^{-1}\left\|{\left(\frac{1}{2}I+\mathcal{D}^\epsilon_B\right)}^{-1}\right\|_{\mathcal{L}\left(H^1(\partial B), H^1(\partial B) \right)}
\end{equation}
with $\mathcal{D}^\epsilon_B \hat{\psi}(\xi):=\int_{\partial B}\frac{\partial\Phi^\epsilon(\xi,\eta)}{\partial \nu(\eta)}\hat{\psi}(\eta) d\eta$ and
$\Phi^{\epsilon}(\xi,\eta):=\frac{e^{i\kappa\epsilon|\xi-\eta|}}{4\pi|\xi-\eta|}$.
\end{lemma}
\begin{proofn}{\it{of Lemma \ref{rep1dbllayer}}.}
\begin{itemize}
\item We have,
 \begin{eqnarray*}
 \mathcal{D}_{ D_\epsilon}\psi(s)&=&\int_{\partial D_\epsilon}\frac{\partial \Phi_\kappa(s,t)}{\partial \nu(t)}\psi(t) dt\\
          &=&\int_{\partial D_\epsilon}\frac{e^{i\kappa|s-t|}}{4\pi|s-t|}\left[\frac{1}{s-t}-ik\right]\frac{s-t}{|s-t|}\cdot \nu(t) \psi(t) dt\\
          &=&\int_{\partial B}\frac{e^{i\kappa\epsilon|\xi-\eta|}}{4\pi\epsilon|\xi-\eta|}\left[\frac{1}{\epsilon(\xi-\eta)}-ik\right]\frac{\xi-\eta}{|\xi-\eta|}\cdot \nu(\eta)\psi(\epsilon\eta+z)\epsilon^{2}d\eta\\
          &=&\int_{\partial B}\frac{e^{i\kappa\epsilon|\xi-\eta|}}{4\pi|\xi-\eta|}\left[\frac{1}{\xi-\eta}-ik\epsilon\right]\frac{\xi-\eta}{|\xi-\eta|}\cdot \nu(\eta)\psi(\epsilon\eta+z)d\eta\\          
          &=&\int_{\partial B}\frac{\partial \Phi^{\epsilon}(\xi,\eta)}{\partial \nu(\eta)}\psi(\epsilon\eta+z) d\eta\\
          &=&\mathcal{D}^\epsilon_B \hat{\psi}(\xi).
 \end{eqnarray*}
The above gives us \eqref{rep1dbllayer0}. From \eqref{rep1dbllayer0}, we can obtain \eqref{rep1dbllayer1} in the following way
\begin{eqnarray*}
 \left(\frac{1}{2}I+\mathcal{D}_{ D_\epsilon}\right)\psi(s)&=&\frac{1}{2}\psi(s)+\mathcal{D}_{ D_\epsilon}\psi(s)
          \,=\,\frac{1}{2}\hat{\psi}(\xi)+\mathcal{D}^\epsilon_B \hat{\psi}(\xi)
          \,=\,\left(\frac{1}{2}I+\mathcal{D}^\epsilon_B\right)\hat{\psi}(\xi).
 \end{eqnarray*}
\item The following equalities
\begin{eqnarray*}
\left(\frac{1}{2}I+\mathcal{D}_{ D_\epsilon}\right)\left({\left(\frac{1}{2}I+\mathcal{D}^\epsilon_B\right)}^{-1} \hat{\phi}\right)^\vee 
                   \,\substack{=\\\eqref{rep1dbllayer1}}\, \left(\left(\frac{1}{2}I+\mathcal{D}^\epsilon_B\right){\left(\frac{1}{2}I+\mathcal{D}^\epsilon_B\right)}^{-1} \hat{\phi}\right)^\vee 
                   \,= \,\hat{\phi}^{\vee}
                  \,=\, \phi,
\end{eqnarray*}
provide us \eqref{rep1dbllayer2}.
\item We have from the estimate,
\begin{eqnarray*}
 \left\|\left(\frac{1}{2}I+\mathcal{D}_{ D_\epsilon}\right)^{-1}\right\|_{\mathcal{L}\left(L^2(\partial D_\epsilon), L^2(\partial D_\epsilon) \right)}&:=&
\substack{Sup\\ \phi(\neq0)\in L^{2}(\partial D_\epsilon)} \frac{\left\|\left(\frac{1}{2}I+\mathcal{D}_{ D_\epsilon}\right)^{-1}\phi\right\|_{L^2(\partial D_\epsilon)}}{ \|\phi \|_{L^2(\partial D_\epsilon)}}\\
&\substack{=\\\eqref{habib2*},\eqref{habib1*}}\,&\substack{Sup\\ \phi(\neq0)\in L^{2}(\partial D_\epsilon)} 
\frac{\epsilon~\left\|\left(\left(\frac{1}{2}I+\mathcal{D}_{ D_\epsilon}\right)^{-1}\phi\right)^{\wedge}\right\|_{L^2(\partial B)}}{\epsilon~ \|\hat{\phi} \|_{L^2(\partial B)}}\\
&\substack{=\\\eqref{rep1dbllayer2}}&\substack{Sup\\ \hat{\phi}(\neq0)\in L^{2}(\partial D_\epsilon)} 
\frac{\left\|\left(\left({\left(\frac{1}{2}I+\mathcal{D}^\epsilon_B\right)}^{-1} \hat{\phi}\right)^\vee\right)^{\wedge}\right\|_{L^2(\partial B)}}{ \|\hat{\phi} \|_{L^2(\partial B)}}\\
&=&\substack{Sup\\ \hat{\phi}(\neq0)\in L^{2}(\partial D_\epsilon)} 
\frac{\left\| {\left(\frac{1}{2}I+\mathcal{D}^\epsilon_B\right)}^{-1} \hat{\phi}\right\|_{L^2(\partial B)}}{ \|\hat{\phi} \|_{L^2(\partial B)}} \\
&=&\left\|{\left(\frac{1}{2}I+\mathcal{D}^\epsilon_B\right)}^{-1}\right\|_{\mathcal{L}\left(L^2(\partial B), L^2(\partial B) \right)}. 
\end{eqnarray*}
It provides us \eqref{nrm1dbllayer2}. By proceeding in the similar manner we can obtain \eqref{nrm1dbllayer2-1} as mentioned below,
\begin{eqnarray*}
 \left\|\left(\frac{1}{2}I+\mathcal{D}_{ D_\epsilon}\right)^{-1}\right\|_{\mathcal{L}\left(H^1(\partial D_\epsilon), H^1(\partial D_\epsilon) \right)}&:=&
\substack{Sup\\ \phi(\neq0)\in H^{1}(\partial D_\epsilon)} \frac{\left\|\left(\frac{1}{2}I+\mathcal{D}_{ D_\epsilon}\right)^{-1}\phi\right\|_{H^1(\partial D_\epsilon)}}{ \|\phi \|_{H^1(\partial D_\epsilon)}}\\
&\substack{\leq\\\eqref{habib2*},\eqref{habib1*}}\,&\substack{Sup\\ \phi(\neq0)\in H^{1}(\partial D_\epsilon)} 
\frac{\left\|\left(\left(\frac{1}{2}I+\mathcal{D}_{ D_\epsilon}\right)^{-1}\phi\right)^{\wedge}\right\|_{H^1(\partial B)}}{\epsilon~ \|\hat{\phi} \|_{H^1(\partial B)}}\\
&\substack{=\\\eqref{rep1dbllayer2}}&\epsilon^{-1}\substack{Sup\\ \hat{\phi}(\neq0)\in H^{1}(\partial D_\epsilon)} 
\frac{\left\|\left(\left({\left(\frac{1}{2}I+\mathcal{D}^\epsilon_B\right)}^{-1} \hat{\phi}\right)^\vee\right)^{\wedge}\right\|_{H^1(\partial B)}}{ \|\hat{\phi} \|_{H^1(\partial B)}}\\
&=&\epsilon^{-1}\substack{Sup\\ \hat{\phi}(\neq0)\in H^{1}(\partial D_\epsilon)} 
\frac{\left\| {\left(\frac{1}{2}I+\mathcal{D}^\epsilon_B\right)}^{-1} \hat{\phi}\right\|_{H^1(\partial B)}}{ \|\hat{\phi} \|_{H^1(\partial B)}} \\
&=&\epsilon^{-1}\left\|{\left(\frac{1}{2}I+\mathcal{D}^\epsilon_B\right)}^{-1}\right\|_{\mathcal{L}\left(H^1(\partial B), H^1(\partial B) \right)}. 
\end{eqnarray*}
\end{itemize} 
\end{proofn}
\noindent
\begin{lemma}\label{lemmanrm1dbllayer31}
  The operator norm of the inverse of $\frac{1}{2}I+\mathcal{D}_{ D_\epsilon}:L^{2}(\partial D_\epsilon)\rightarrow L^{2}(\partial D_\epsilon)$ defined by
$\mathcal{D}_{ D_\epsilon}\psi(s):=\int_{\partial D_\epsilon}\frac{\partial\Phi_\kappa(s,t)}{\partial\nu(t)}\psi(t)dt$ in \eqref{defofDpartialDe}, is bounded by a constant, 
i.e.\begin{eqnarray}\label{nrm1dbllayer31}
 \left\|\left(\frac{1}{2}I+\mathcal{D}_{ D_\epsilon}\right)^{-1}\right\|_{\mathcal{L}\left(L^2(\partial D_\epsilon), L^2(\partial D_\epsilon) \right)}
&\leq&\grave{C}_6,
\end{eqnarray}
with 
$\grave{C}_6:=\frac{4\pi\left\|\left(\frac{1}{2}I+{\mathcal{D}_{B}^{i_{\kappa}}}\right)^{-1}\right\|_{\mathcal{L}\left(L^2(\partial B), L^2(\partial B) \right)}}
{4\pi-\kappa^2\epsilon^2|\partial B|\left\|\left(\frac{1}{2}I+{\mathcal{D}_{B}^{i_{\kappa}}}\right)^{-1}\right\|_{\mathcal{L}\left(L^2(\partial B), L^2(\partial B) \right)}}$. 
Here, $\mathcal{D}_{B}^{i_{\kappa}}:L^2(\partial B)\rightarrow L^2(\partial B)$ is the double layer potential with the wave number zero.
\end{lemma}
\noindent

Here we should mention that if 
$\epsilon^2\leq\frac{\pi}{\kappa^2|\partial B|\left\|\left(\frac{1}{2}I+{\mathcal{D}_{B}^{i_{\kappa}}}\right)^{-1}\right\|_{\mathcal{L}\left(L^2(\partial B), L^2(\partial B) \right)}}$,
then $\grave{C}_6$ is bounded by 
\\
$\frac{4}{3}\left\|\left(\frac{1}{2}I+{\mathcal{D}_{B}^{i_{\kappa}}}\right)^{-1}\right\|_{\mathcal{L}\left(L^2(\partial B), L^2(\partial B) \right)}$, 
which is a universal constant depending only on $\partial B$ through its Lipschitz character, see Remark \ref{Lipschitz-character}.
\begin{proofn}{\it{of Lemma \ref{lemmanrm1dbllayer31}}.} To estimate the operator norm of $\left(\frac{1}{2}I+\mathcal{D}_{ D_\epsilon}\right)^{-1}$ we decompose 
$\mathcal{D}_{ D_\epsilon}=:\mathcal{D}_{ D_\epsilon}^\kappa=\mathcal{D}_{ D_\epsilon}^{i_{\kappa}}+\mathcal{D}_{ D_\epsilon}^{d_{\kappa}}$ into 
two parts $\mathcal{D}_{ D_\epsilon}^{i_{\kappa}}$ ( independent of $\kappa$ ) and  $\mathcal{D}_{ D_\epsilon}^{d_{\kappa}}$ ( dependent of $\kappa$ ) given by
\begin{eqnarray}
\mathcal{D}_{ D_\epsilon}^{i_{\kappa}} \psi(s):=\int_{\partial D_\epsilon} \left(\frac{\partial}{\partial \nu(t)}\Phi_0(s,t)\right)\psi(t) dt,\label{def-of_Dik}\\
\mathcal{D}_{ D_\epsilon}^{d_{\kappa}} \psi(s):=\int_{\partial D_\epsilon} \left(\frac{\partial}{\partial \nu(t)}[\Phi_\kappa(s,t)-\Phi_0(s,t)]\right)\psi(t) dt.\label{def-of_Ddk}
\end{eqnarray}
 With this definition, $\frac{1}{2}I+\mathcal{D}_{ D_\epsilon}^{i_{\kappa}}:L^2(\partial D_\epsilon)\rightarrow L^2(\partial D_\epsilon)$ is invertible, see \cite{MD:InteqnsaOpethe1997}. Hence, 
$\frac{1}{2}I+\mathcal{D}_{ D_\epsilon}=\left(\frac{1}{2}I+\mathcal{D}_{ D_\epsilon}^{i_{\kappa}}\right)\left(I+\left(\frac{1}{2}I+\mathcal{D}_{ D_\epsilon}^{i_{\kappa}}\right)^{-1}\mathcal{D}_{ D_\epsilon}^{d_{\kappa}}\right)$ and so
\begin{eqnarray}
 \left\|\left(\frac{1}{2}I+\mathcal{D}_{ D_\epsilon}\right)^{-1}\right\|_{\mathcal{L}\left(L^2(\partial D_\epsilon), L^2(\partial D_\epsilon) \right)}
&\leq& \left\|\left(I+\left(\frac{1}{2}I+\mathcal{D}_{ D_\epsilon}^{i_{\kappa}}\right)^{-1}\mathcal{D}_{ D_\epsilon}^{d_{\kappa}}\right)^{-1}\right\|_{\mathcal{L}\left(L^2(\partial D_\epsilon), L^2(\partial D_\epsilon) \right)}\nonumber \\
 &&                   \times       \left\|\left(\frac{1}{2}I+\mathcal{D}_{ D_\epsilon}^{i_{\kappa}}\right)^{-1}\right\|_{\mathcal{L}\left(L^2(\partial D_\epsilon), L^2(\partial D_\epsilon) \right)}.\label{nrm1dbllayer3}
\end{eqnarray}
So, to estimate the operator norm of $\left(\frac{1}{2}I+\mathcal{D}_{ D_\epsilon}\right)^{-1}$ one needs to estimate the 
operator norm of $\left(I+\left(\frac{1}{2}I+\mathcal{D}_{ D_\epsilon}^{i_{\kappa}}\right)^{-1}\mathcal{D}_{ D_\epsilon}^{d_{\kappa}}\right)^{-1}$, 
in particular one needs to have the knowledge about the operator norms of
 $\left(\frac{1}{2}I+\mathcal{D}_{ D_\epsilon}^{i_{\kappa}}\right)^{-1}$ and $\mathcal{D}_{ D_\epsilon}^{d_{\kappa}}$ to apply the Neumann series. 
 For that purpose,  we can estimate the operator norm of  $\left(\frac{1}{2}I+\mathcal{D}_{ D_\epsilon}^{i_{\kappa}}\right)^{-1}$ from \eqref{nrm1dbllayer2} by
\begin{eqnarray}\label{optnormdik}
  \left\|\left(\frac{1}{2}I+\mathcal{D}_{ D_\epsilon}^{i_{\kappa}}\right)^{-1}\right\|_{\mathcal{L}\left(L^2(\partial D_\epsilon), L^2(\partial D_\epsilon) \right)}
   &=& \left\|\left(\frac{1}{2}I+\mathcal{D}_{B}^{i_{\kappa}}\right)^{-1}\right\|_{\mathcal{L}\left(L^2(\partial B), L^2(\partial B) \right)}.
 \end{eqnarray}
Here $\mathcal{D}^{i_{\kappa}}_B \hat{\psi}(\xi):=\int_{\partial B} \frac{1}{4\pi}\frac{\xi-\eta}{|\xi-\eta|^3}\cdot \nu(\eta)\,\hat{\psi}(\eta) d\eta$. 
From the definition of the operator $\mathcal{D}_{ D_\epsilon}^{d_{k}}$ in \eqref{def-of_Ddk}, by performing the similar calculations made in $\eqref{TanS^{d_{k}}}$, we deduce that
\begin{eqnarray}\label{D^{d_{k}}}
 \mathcal{D}_{ D_\epsilon}^{d_{\kappa}} \psi(s)
                                              &=&\int_{\partial D_\epsilon} \left(\frac{\partial}{\partial \nu(t)}\frac{e^{i\kappa|s-t|}-1}{4\pi|s-t|}\right)\psi(t) dt\nonumber\\
&=&\int_{\partial D_\epsilon} 
     \left(\frac{e^{i\kappa|s-t|}}{4\pi|s-t|}\left[\frac{1}{|s-t|}-i\kappa\right]\frac{s-t}{|s-t|}-\frac{1}{4\pi}\frac{s-t}{|s-t|^3}\right)\cdot \nu(t) \psi(t) dt\nonumber \\
&=&\frac{\kappa^2\epsilon^{2}}{4\pi}\int_{\partial B} 
     \left[\sum_{l=1}^{\infty}\left(\frac{1}{l!}-\frac{1}{l+1!}\right)(i\kappa\epsilon)^{l-1}|\xi-\eta|^{l-1}\right]
\frac{\xi-\eta}{|\xi-\eta|}\cdot \nu(\eta) \hat{\psi}(\eta) d\eta.
\end{eqnarray}
Now, by performing the similar calculations made in \eqref{modTanS^{d_{k}}}, \eqref{D^{d_{k}}} will direct us to calculate the below
 \begin{eqnarray}\label{modD^{d_{k}}}
\left|\mathcal{D}_{ D_\epsilon}^{d_{\kappa}} \psi(s)\right|&\leq&\frac{\kappa^2\epsilon^{2}}{4\pi}
 \left\|\sum_{l=1}^{\infty}\left(\frac{1}{l!}-\frac{1}{l+1!}\right)(i\kappa\epsilon)^{l-1}|\xi-\cdot|^{l-1}\right\|_{L^2(\partial B)} \|\hat{\psi}\|_{L^2(\partial B)} \nonumber\\
&\leq&\frac{\kappa^2\epsilon^{2}}{4\pi}|\partial B|^{\frac{1}{2}}
\|\hat{\psi}\|_{L^2(\partial B)},\,\text{for}\,\left(\epsilon\leq\frac{1}{\kappa_{\max}diam(B_m)}\right)\equiv (\kappa_{\max}~a\leq1).
 \end{eqnarray}
We set $\grave{C}_1:=\frac{|\partial B|^{\frac{1}{2}}}{4\pi}$.
From \eqref{modD^{d_{k}}}, by doing the similar calculations used to dervie $\eqref{normS^{d_{k}}}$, we obtain
\begin{eqnarray}\label{normD^{d_{k}}}
 \left\|\mathcal{D}_{ D_\epsilon}^{d_{\kappa}} \psi\right\|_{L^2(\partial D_\epsilon)}&\leq&
\grave{C}_1 \kappa^2\epsilon^{3}|\partial B|^\frac{1}{2}
  \|\hat{\psi}\|_{L^2(\partial B)}.
\end{eqnarray}
We estimate norm of the operator $\mathcal{D}_{ D_\epsilon}^{d_{\kappa}}$ as
\begin{eqnarray}\label{fnormD^{d_{k}}L2}
 \left\|\mathcal{D}_{ D_\epsilon}^{d_{\kappa}}\right\|_{\mathcal{L}\left(L^2(\partial D_\epsilon), L^2(\partial D_\epsilon)\right)}&=&
\substack{Sup\\ \psi(\neq0)\in L^{2}(\partial D_\epsilon)} \frac{ \|\mathcal{D}_{ D_\epsilon}^{d_{\kappa}}\psi \|_{L^2(\partial D_\epsilon)}}{ \|\psi \|_{L^2(\partial D_\epsilon)}}\nonumber\\
&\leq&
 \substack{Sup\\ \hat{\psi}(\neq0)\in L^{2}(\partial B)} \frac{\grave{C}_1\kappa^2\epsilon^{3}|\partial B|^\frac{1}{2}\|\hat{\psi}\|_{L^2(\partial B)}}
{\epsilon~\|\hat{\psi}\|_{L^2(\partial B)}}~~(\mbox{ By }\eqref{normD^{d_{k}}},\,\eqref{habib2*})\nonumber\\
&=&\grave{C}_1\kappa^2\epsilon^{2}|\partial B|^\frac{1}{2}.
\end{eqnarray}
Hence, we get
 \begin{eqnarray}\label{norms0dk}
  \left\|\left(\frac{1}{2}I+{\mathcal{D}_{ D_\epsilon}^{i_{\kappa}}}\right)^{-1}\mathcal{D}_{ D_\epsilon}^{d_{\kappa}}\right\|_{\mathcal{L}\left(L^2(\partial D_\epsilon), L^2(\partial D_\epsilon) \right)}
&\substack{\leq\\ 
\eqref{optnormdik},\eqref{fnormD^{d_{k}}L2}}&\left\|\left(\frac{1}{2}I+{\mathcal{D}_{ B}^{i_{\kappa}}}\right)^{-1}\right\|_{\mathcal{L}\left(L^2(\partial B), L^2(\partial B) \right)}\grave{C}_1\kappa^2\epsilon^{2}|\partial B|^\frac{1}{2}\nonumber\\
&\leq&\grave{C}_2\kappa^2\epsilon^2 ,
 \end{eqnarray}
where $\grave{C}_2:=\grave{C}_1|\partial B|^\frac{1}{2}\left\|\left(\frac{1}{2}I+{\mathcal{D}_{B}^{i_{\kappa}}}\right)^{-1}\right\|_{\mathcal{L}\left(L^2(\partial B), L^2(\partial B) \right)}
=\frac{|\partial B|}{4\pi}\left\|\left(\frac{1}{2}I+{\mathcal{D}_{B}^{i_{\kappa}}}\right)^{-1}\right\|_{\mathcal{L}\left(L^2(\partial B), L^2(\partial B) \right)}.$
Assuming $\epsilon$ to satisfy the condition $\epsilon<\frac{1}{\sqrt{\grave{C}_2}\kappa_{\max}}$, then 
$\left\|\left(\frac{1}{2}I+{\mathcal{D}_{ D_\epsilon}^{i_{\kappa}}}\right)^{-1}\mathcal{D}_{ D_\epsilon}^{d_{\kappa}}\right\|_{\mathcal{L}\left(L^2(\partial D_\epsilon), L^2(\partial D_\epsilon) \right)}<1$
and hence by using the Neumann series we obtain the following
\begin{eqnarray*}
 \left\|\left(I+\left(\frac{1}{2}I+{\mathcal{D}_{ D_\epsilon}^{i_{\kappa}}}\right)^{-1}\mathcal{D}_{ D_\epsilon}^{d_{\kappa}}\right)^{-1}\right\|_{\mathcal{L}\left(L^2(\partial D_\epsilon), L^2(\partial D_\epsilon) \right)}
 &\leq&
 \frac{1}{1-\left\|\left(\frac{1}{2}I+{\mathcal{D}_{ D_\epsilon}^{i_{\kappa}}}\right)^{-1}\mathcal{D}_{ D_\epsilon}^{d_{\kappa}}\right\|_{\mathcal{L}\left(L^2(\partial D_\epsilon), L^2(\partial D_\epsilon) \right)}}\\
 &\substack{\leq\\ \eqref{norms0dk}}&\grave{C}_3:=\frac{1}{1-\grave{C}_2 \kappa^2\epsilon^2}.\\
\end{eqnarray*}
 By substituting the above and \eqref{optnormdik} in \eqref{nrm1dbllayer3}, we obtain the required result \eqref{nrm1dbllayer31}.
\end{proofn}

\paragraph{The multiple obstacle case}\label{gbmulobscasedbl}
\ ~ \ \\
\begin{proposition} \label{propsmjsmmestdbl}
 For $m,j=1,2,\dots,M$, the operator $\mathcal{D}_{mj}:L^2(\partial D_j)\rightarrow L^{2}(\partial D_m)$ defined 
in Proposition  \ref{existence-of-sigmasdbl}, see \eqref{defofDmjed}, 
satisfies the following estimates,
\begin{itemize}
\item
For $j=m$,
\begin{eqnarray}\label{estinvsmmdbl}
 \left\|\left(\frac{1}{2}I+\mathcal{D}_{mm}\right)^{-1}\right\|_{\mathcal{L}\left(L^2(\partial D_m), L^2(\partial D_m) \right)}&\leq&\grave{C}_{6m}, \label{invfnormD_{ii}1}
\end{eqnarray}
where $\grave{C}_{6m}:=\frac{4\pi\left\|\left(\frac{1}{2}I+\mathcal{D}^{i_{\kappa}}_{B_m}\right)^{-1}\right\|_{\mathcal{L}\left(L^2(\partial B_m), L^2(\partial B_m) \right)}}
{4\pi-\kappa^2\epsilon^2|\partial B_m|\left\|\left(\frac{1}{2}I+\mathcal{D}^{i_{\kappa}}_{B_m}\right)^{-1}\right\|_{\mathcal{L}\left(L^2(\partial B_m), L^2(\partial B_m) \right)}}$.
\item For $j\neq m$,
\begin{eqnarray}\label{estinvsmjdbl}
\left\|\mathcal{D}_{mj}\right\|_{\mathcal{L}\left(L^2(\partial D_j), L^2(\partial D_m)\right)}
                          &\leq&\frac{1}{4\pi}\left(\frac{\kappa}{d}+\frac{1}{d^2}\right)\left|\partial \c{B} \right|\epsilon^{2},\label{fnormD_{ij}21}
\end{eqnarray}
where $\left|\partial \c{B} \right|:=\max\limits_m \partial B_m$.

\end{itemize}
\end{proposition}
\begin{proofn}{\it{of Proposition \ref{propsmjsmmestdbl}}.}
The estimate \eqref{estinvsmmdbl} is nothing but \eqref{nrm1dbllayer31} of Lemma \ref{lemmanrm1dbllayer31}, replacing $B$ by $B_m$, $z$ by $z_m$ and $D_\epsilon$ by $D_m$ respectively. It remains to prove the estimate \eqref{estinvsmjdbl}.
We have
\begin{eqnarray}\label{fnormD_{ij}}
 \left\|\mathcal{D}_{mj}\right\|_{\mathcal{L}\left(L^2(\partial D_j), L^2(\partial D_m)\right)}&=&
\substack{Sup\\ \psi(\neq0)\in L^{2}(\partial D_j)} \frac{ \|\mathcal{D}_{mj}\psi \|_{L^2(\partial D_m)}}{ \|\psi \|_{L^2(\partial D_j)}}.
\end{eqnarray}
\noindent
Let $\psi\in L^2(\partial D_j)$ then for $s\in \partial D_m$,  we have
\begin{eqnarray}\label{modD_{ij}}
\left|\mathcal{D}_{mj}\psi(s)\right|&=&\left|\int_{\partial D_j}\frac{\partial\Phi_\kappa(s,t)}{\partial \nu_j(t)}\psi(t)dt\right|\nonumber\\
                                                               &=&\left|\int_{\partial D_j}\nabla_t\Phi_\kappa(s,t)\cdot\nu_j(t)\psi(t) dt\right|\nonumber\\
                                                               &\leq&\int_{\partial D_j}\left|\nabla_t\Phi_\kappa(s,t)\right|~\left|\psi(t)\right| dt\nonumber\\    
                                                               &\leq&\frac{1}{4\pi}\left(\frac{\kappa}{d_{mj}}+\frac{1}{d_{mj}^2}\right)\epsilon
                                                                             \left|\partial B_j\right|^\frac{1}{2} \|\psi \|_{L^2(\partial D_j)}.                  
\end{eqnarray}
\noindent
Here we used the similar calculations made in \eqref{modTS_{ij}}. From \eqref{modD_{ij}}, by proceeding further in the way of \eqref{fL2normTS_{ij}}, we get
\begin{eqnarray}\label{fL2normD_{ij}}
\left\|\mathcal{D}_{mj}\psi\right\|_{L^2(\partial D_m)}
&\leq&\frac{1}{4\pi}\left(\frac{\kappa}{d_{mj}}+\frac{1}{d_{mj}^2}\right)\epsilon^{2}
                                                                                             \left|\partial B_j\right|^\frac{1}{2}\left|\partial B_m\right|^\frac{1}{2} \|\psi \|_{L^2(\partial D_j)}.
\end{eqnarray}
\noindent
Substitution of \eqref{fL2normD_{ij}} in \eqref{fnormD_{ij}} gives us
\begin{eqnarray*}\label{fnormD_{ij}1}
 \left\|\mathcal{D}_{mj}\right\|_{\mathcal{L}\left(L^2(\partial D_j), L^2(\partial D_m)\right)}
                          &\leq&\frac{1}{4\pi}\left(\frac{\kappa}{d_{mj}}+\frac{1}{d_{mj}^2}\right)\epsilon^2
                                                                   \left|\partial B_j\right|^\frac{1}{2}\left|\partial B_m\right|^\frac{1}{2}
                          \,\leq\,\frac{1}{4\pi}\left(\frac{\kappa}{d}+\frac{1}{d^2}\right)\left|\partial \c{B}\right|\epsilon^{2}.
\end{eqnarray*}
\end{proofn}

\begin{proofe}
\textbf{\textit{End of the proof of Proposition  \ref{normofsigmastmtdbl}}.}
By substituting \eqref{invfnormD_{ii}1} in \eqref{invDLnrm} and \eqref{fnormD_{ij}21} in \eqref{DKnrm}, we obtain
\begin{eqnarray}
\left\|DK\right\|
    &\equiv&
\max\limits_{1\leq m \leq M}\sum_{\substack{j=1\\j\neq\,m}}^{M}\left\|\mathcal{D}_{mj}\right\|_{\mathcal{L}\left(L^{2}(\partial D_j),L^{2}(\partial D_m)\right)}\nonumber\\
    &\leq&\frac{M-1}{4\pi}\left(\frac{\kappa}{d}+\frac{1}{d^2}\right)\left|\partial \c{B}\right|\epsilon^{2}\label{DKnrm1}
\end{eqnarray}
and
\begin{eqnarray}
\left\|\left(\frac{1}{2}\textbf{I}+DL\right)^{-1}\right\|
     &\equiv&\max\limits_{1\leq m \leq M}\left\|\left(\frac{1}{2}I+\mathcal{D}_{mm}\right)^{-1}\right\|_{\mathcal{L}\left(L^{2}(\partial D_m),L^{2}(\partial D_m)\right)}\nonumber\\ 
     &\equiv&\max\limits_{1\leq m \leq M}\grave{C}_{6m}.\label{invDLnrm1}
\end{eqnarray}
Hence, \eqref{DKnrm1} and \eqref{invDLnrm1} jointly provide
\begin{eqnarray}
\left\|\left(\frac{1}{2}\textbf{I}+DL\right)^{-1}\right\|\left\|DK\right\|
    &\leq&\underbrace{\frac{M-1}{4\pi}\left(\max\limits_{1\leq m \leq M}\grave{C}_{6m}\right)\left|\partial \c{B}\right|\left(\frac{\kappa}{d}+\frac{1}{d^2}\right)\epsilon^2}_{=:\grave{C}_s} ,\label{invDLDKnrm1}
\end{eqnarray}
By imposing the condition $\left\|\left(\frac{1}{2}\textbf{I}+DL\right)^{-1}\right\|\left\|DK\right\|<1$, we get the following from \eqref{nrminvDLplusDK} and (\ref{sigmaU^{In}nrmdbl1}-\ref{sigmaU^{In}nrmdbl});
\begin{eqnarray}\label{nrminvDLplusDK2}
 \left\|\sigma_m\right\|_{L^{2}(\partial D_m)}\leq\left\|\sigma\right\| 
                                  &\leq&\frac{\left\|\left(\frac{1}{2}\textbf{I}+DL\right)^{-1}\right\|}{1-\left\|\left(\frac{1}{2}\textbf{I}+DL\right)^{-1}\right\|\left\|DK\right\|}\left\|U^{In}\right\|\nonumber\\
                   &\leq&\grave{C}_p\left\|\left(\frac{1}{2}\textbf{I}+DL\right)^{-1}\right\| \max\limits_{1\leq m \leq M}\left\|U^{i}\right\|_{L^{2}(\partial D_m)}\hspace{.25cm} \left( \grave{C}_p\geq\frac{1}{1-\grave{C}_s}\right)\nonumber\\
                   &\substack{\leq\\ \eqref{invDLnrm1} }&\grave{\mathrm{C}}  \max\limits_{1\leq m \leq M}\left\|U^{i}\right\|_{L^{2}(\partial D_m)}\hspace{.25cm} \left(\grave{\mathrm{C}}:=\grave{C}_p\max\limits_{1\leq m \leq M}\grave{C}_{6m}\right),
\end{eqnarray}
for all $m\in\{1,2,\dots,M\}$.
Since, $U^{i}$ denotes the plane incident wave given by $U^{i}(x,\theta)=e^{i\kappa{x}\cdot\theta}$, we have
\begin{eqnarray}\label{L2normU^i}
\left\|U^{i}\right\|_{L^{2}(\partial D_m)}&=&\epsilon~\left|\partial B_m\right|^{\frac{1}{2}}, \forall m=1,2,\dots,M.
\end{eqnarray}
Now by substituting \eqref{L2normU^i} in \eqref{nrminvDLplusDK2}, for each $m=1,\dots,M$, we obtain
\begin{eqnarray}\label{nrmsigmafdbl}
  \left\|\sigma_m\right\|_{L^{2}(\partial D_m)}&\,\leq\,&  \grave{\mathcal{C}}(\kappa)\epsilon,
\end{eqnarray}
where $\hspace{.25cm}\grave{\mathcal{C}}(\kappa):=\grave{\mathrm{C}} \left|\partial \c{B}\right|^{\frac{1}{2}}$. 
\par
The condition $\left\|\left(\frac{1}{2}I+DL\right)^{-1}\right\|\left\|DK\right\|<1$ is
satisfied if 
\begin{eqnarray}\label{Accond-invLK-singl-smalldbl}
\frac{M-1}{4\pi}\left|\partial \c{B}\right|\left(\frac{\kappa}{d}+\frac{1}{d^2}\right)\left(\max\limits_{1\leq m \leq M}\grave{C}_{6m}\right)\epsilon^2<1.
\end{eqnarray}
Since $\kappa\,d<\bar{c}$ for $\bar{c}:=\kappa_{\max}d_{\max}$, then \eqref{Accond-invLK-singl-smalldbl} 
reads as $(M-1)\epsilon^2<cd^2$, where we set $c:=\left[\frac{(\bar{c}+1)}{4\pi}|\partial\c{B}|\max\limits_{1\leq m \leq M}\grave{C}_{6m}\right]^{-1}$ and $\grave{c}:=\sqrt{c}$ 
will serve our purpose in Proposition  \ref{normofsigmastmtdbl} and hence in Theorem \ref{Maintheorem-ac-small-sing} for the case $\sqrt{M-1}\epsilon<\grave{c}d$.
\par Remark here that if we do not have the condition \eqref{def-dmax}, then \eqref{Accond-invLK-singl-smalldbl} holds
if $(M-1)\epsilon<cd$ with $c:=\frac{1}{2}\left[\frac{|\partial\c{B}|}{4\pi}\max\limits_{1\leq m \leq M}\grave{C}_{6m}\right]^{-1}$
(for $M$ large enough). Then $c_2:=c\,\max\limits_{1\leq m \leq M } diam (B_m)$ 
will serve our purpose in Remark \eqref{remarkof-Maintheorem-ac-small-sing}.
\end{proofe}

For $m=1,2,\dots,M$, let $U^{\sigma_m}$ be the solution of the problem
\begin{equation}\label{acimpenetrableUsigma}
\begin{cases}
(\Delta + \kappa^{2})U^{\sigma_m}=0& \mbox{ in }D_m,\\
U^{\sigma_m}=\sigma_m& \mbox{ on } \partial D_m.
\end{cases}
\end{equation}
The function $\sigma_m$ is in $H^{1}(\partial D_m)$, see Proposition \ref{existence-of-sigmasdbl}. Hence $U^{\sigma_m}\in\,H^{\frac{3}{2}}(D_m)$ and then $\left.\frac{\partial\,U^{\sigma_{m}}}{\partial \nu_m}\right|_{\partial D_m}\in\,L^2(\partial D_m)$.
From Proposition  \ref{existence-of-sigmasdbl}, the solution of the problem (\ref{acimpoenetrable}-\ref{radiationc}) has the form
 \begin{equation}\label{qcimprequiredfrm2dbl}
  U^{t}(x)=U^{i}(x)+\sum_{m=1}^{M}\int_{\partial D_m}\frac{\partial\Phi_\kappa(x,s)}{\partial \nu_m(s)}\sigma_{m} (s)ds,~x\in\mathbb{R}^{3}\backslash\left(\mathop{\cup}_{m=1}^M \bar{D}_m\right).
\end{equation}
It can be written in terms of single layer potanetial using Gauss theorem as
 \begin{equation}\label{qcimprequiredfrm3dbl}
  U^{t}(x)=U^{i}(x)+\sum_{m=1}^{M}\int_{\partial D_m}\Phi_\kappa(x,s)\frac{\partial\,U^{\sigma_{m}} (s)}{\partial \nu_m(s)}ds,~x\in\mathbb{R}^{3}\backslash\left(\mathop{\cup}_{m=1}^M \bar{D}_m\right).
\end{equation}
\begin{description}
\item Indeed,
\begin{eqnarray}\label{revised-greensform-DtN-acoustic-small}
 \hspace{-2cm}\int_{\partial D_m} \hspace{-.2cm}\frac{\partial\Phi_\kappa(x,s)}{\partial \nu_m(s)}\sigma_{m} (s)ds=\int_{\partial D_m}\hspace{-.3cm}\Phi_\kappa(x,s)\frac{\partial\,U^{\sigma_{m}} (s)}{\partial \nu_m(s)}ds\hspace{-.03cm}+\hspace{-.03cm}
                       \int_{D_m}\hspace{-.1cm}[U^{\sigma_{m}} (y)\Delta\,\Phi_\kappa(x,y)-\Phi_\kappa(x,y)\Delta\,U^{\sigma_{m}} (y)]dy.
\end{eqnarray}
\end{description}
\bigskip

\begin{lemma}\label{thmestmdhoudhonu}
For $m=1,2,\dots,M$, $U^{\sigma_m}$, the solutions of the problem \eqref{acimpenetrableUsigma} satisfy the estimate
\begin{equation}\label{estmdhoudhonu}
 \left\|\frac{\partial\,U^{\sigma_{m}} (s)}{\partial \nu_m(s)}\right\|_{H^{-1}(\partial D_m)}\le C_7,
\end{equation}
for some constant $C_7$ depending on the  Lipschitz character of $B_m$ but it is independent of $\epsilon$.
\end{lemma}
\begin{proofn}{\it{of Lemma \ref{thmestmdhoudhonu}}.}
 For $m=1,2,\dots,M$, write $$\mathcal{U}^{m}(x):=U^{\sigma_{m}}(\epsilon\,x+z_m),\forall\,x\in\,B_m.$$
 Then we obtain
\begin{equation}\label{mathcalUmeqn}
\begin{cases}
 (\Delta+\epsilon^2\kappa^2)\mathcal{U}^{m}(x)&=\epsilon^2(\Delta+\kappa^2)U^{\sigma_{m}}(\epsilon\,x+z_m)=0,\text{ for } x\in\,B_m,\\
\hspace{1cm}\mathcal{U}^{m}(\xi)&=U^{\sigma_{m}}(\epsilon\,\xi+z_m)=\sigma(\epsilon\,\xi+z_m),\text{ for } \xi\in\partial\,B_m,
\end{cases}
\end{equation}
and also $\frac{\partial\mathcal{U}^{m}(\xi)}{\partial\nu_m(\xi)}=\nabla\mathcal{U}^{m}(\xi)\cdot\nu_m(\xi)=\epsilon\nabla\,U^{\sigma_{m}}(\epsilon\,\xi+z_m)\cdot\nu_m(\epsilon\,\xi+z_m)=\epsilon\frac{\partial\mathcal{U}^{\sigma_m}(\epsilon\,\xi+z_m)}{\partial\nu_m}$.
Hence,
\begin{eqnarray*}
 \left\|\frac{\partial\mathcal{U}^{m}}{\partial\nu_m}\right\|^2_{L^2(\partial B_m)}&=&\int_{\partial B_m}\left|\frac{\partial\mathcal{U}^{m}(\eta)}{\partial\nu_m(\eta)}\right|^2d\eta\nonumber\\
&=&\int_{\partial D_m}\epsilon^2\left|\frac{\partial\,U^{\sigma_{m}}(s)}{\partial\nu_m(s)}\right|^2\epsilon^{-2}ds,~[s:=\epsilon\eta+z_m]\nonumber\\
&=& \left\|\frac{\partial\,U^{\sigma_{m}}}{\partial\nu_m}\right\|^2_{L^2(\partial D_m)},
\end{eqnarray*}
 which gives us
\begin{eqnarray}\label{D2Nbasic}
 \frac{\left\|\frac{\partial\,U^{\sigma_{m}}}{\partial\nu_m}\right\|_{L^{2}(\partial D_m)}}{ \|U^{\sigma_{m}} \|_{H^1(\partial D_m)}}
&\substack{\le \\ \eqref{habib1*}}&\frac{ \|\frac{\partial\,\mathcal{U}^{{m}}}{\partial\nu_m} \|_{L^2(\partial B_m)}}{\epsilon \|\mathcal{U}^{m} \|_{H^1(\partial B_m)}}.
\end{eqnarray}
For every function $\zeta_m\in\,H^1(\partial D_m)$, the corresponding $U^{\zeta_{m}}$ exists on $D_m$ as mentioned in \eqref{acimpenetrableUsigma} and then the corresponding functions on $B_m$ and the inequality \eqref{D2Nbasic} 
will be satisfied by all these functions. Let $\varLambda_{D_m}:H^1(\partial D_m)\rightarrow\,L^2(\partial D_m)$ and $\varLambda_{B_m}:H^1(\partial\, B_m)\rightarrow\,L^2(\partial\, B_m)$ be the Dirichlet to Neumann maps. 
Then we get the following estimate from \eqref{D2Nbasic}.
$$\|\varLambda_{D_m}\|_{\mathcal{L}\left(H^1( \partial D_m), L^{2}( \partial D_m) \right)}\leq\frac{1}{\epsilon}\|\varLambda_{B_m}\|_{\mathcal{L}\left(H^1(\partial\, B_m), L^{2}(\partial\, B_m) \right)}.$$
This implies that,
\begin{eqnarray}
 \frac{\left\|\frac{\partial\,U^{\sigma_{m}}}{\partial\nu_m}\right\|_{H^{-1}(\partial D_m)}}{ \|U^{\sigma_{m}} \|_{L^2(\partial D_m)}}&\le&\|\varLambda_{D_m}^{*}\|_{\mathcal{L}\left(L^2( \partial D_m), H^{-1}( \partial D_m) \right)}\nonumber\\
&=&\|\varLambda_{D_m}\|_{\mathcal{L}\left(H^1( \partial D_m), L^{2}( \partial D_m) \right)}\nonumber\\
&\leq&\frac{1}{\epsilon}\|\varLambda_{B_m}\|_{\mathcal{L}\left(H^1(\partial\, B_m), L^{2}(\partial\, B_m) \right)}.
\end{eqnarray}
Now, by \eqref{nrmsigmafdbl} and \eqref{acimpenetrableUsigma}, we obtain
\begin{eqnarray}
 \left\|\frac{\partial\,U^{\sigma_{m}}}{\partial\nu_m}\right\|_{H^{-1}(\partial D_m)}
&\leq&\grave{\mathcal{C}}(\kappa)\|\varLambda_{B_m}\|_{\mathcal{L}\left(H^1(\partial\, B_m), L^{2}(\partial\, B_m) \right)}.
\end{eqnarray}
Hence the result is true with $C_7 :=\grave{\mathcal{C}}(\kappa)\|\varLambda_{B_m}\|_{\mathcal{L}\left(H^1(\partial\, B_m), L^{2}(\partial\, B_m) \right)}$ as  $\|\varLambda_{B_m}\|_{\mathcal{L}\left(H^1(\partial\, B_m), L^{2}(\partial\, B_m) \right)}$ is bounded by a 
constant depending only on $B_m$ through its size and Lipschitz character of $B_m$, see Remark \ref{Lipschitz-character}. 
\end{proofn}
\subsubsection{Further estimates on the total charge $\int_{ \partial D_m} \frac{\partial\,U^{\sigma_{m}} (s)}{\partial \nu_m(s)} ds,\,m=1,\dots\,M$}\label{DLPR-3}
\ ~ \
\begin{definition}
\label{Qmdefdbl}
Similarly to Definition \ref{Qmdef}, we call $\sigma_m\in L^2(\partial D_m)$ satisfying \eqref{qcimprequiredfrm1dbl}, the solution of the problem (\ref{acimpoenetrable}-\ref{radiationc}), as surface charge distributions. Using these surface charge distributions 
we define the total charge on each surface $\partial D_m$ denoted by $Q_m$ as
\begin{eqnarray}\label{defofQmdbl}
Q_m:=\int_{ \partial D_m} \frac{\partial\,U^{\sigma_{m}} (s)}{\partial \nu_m(s)} ds.
\end{eqnarray}
\end{definition}
\begin{lemma}\label{Qmestbigodbl}
For $m=1,2,\dots,M$, the absolute value of the total charge $Q_m$ on each surface $\partial D_m$ is bounded by $\epsilon$, i.e. 
\begin{equation}\label{estofQmdbl}
 |Q_m|\leq\grave{\tilde{c}}\epsilon,
\end{equation}
where $\grave{\tilde{c}}:=|\partial \c{B}|\grave{\mathrm{C}}\|\varLambda_{B_m}\|_{\mathcal{L}\left(H^1(\partial\, B_m), L^{2}(\partial\, B_m) \right)}$ 
with $\partial\c{B}$ and $\grave{\mathrm{C}}$ are defined in 
\eqref{estinvsmjdbl} and \eqref{nrminvDLplusDK2} respectively.
\end{lemma}
\begin{proofn}{\it{of Lemma \ref{Qmestbigodbl}}.}
From Proposition  \ref{normofsigmastmtdbl}, we have the estimate for the surface charge distributions $\sigma_m$ as
$
  \left\|\sigma_m\right\|_{L^{2}(\partial D_m)}\,\leq\,  \grave{\mathcal{C}}(\kappa)\epsilon,
$
 with $\grave{\mathcal{C}}(\kappa):=\grave{\mathrm{C}} \left|\partial \c{B}\right|^{\frac{1}{2}}$, which results the Lemma \ref{thmestmdhoudhonu} as $\left\|\frac{\partial\,U^{\sigma_{m}}}{\partial\nu_m}\right\|_{H^{-1}(\partial D_m)}
\leq\grave{\mathcal{C}}(\kappa)\|\varLambda_{B_m}\|_{\mathcal{L}\left(H^1(\partial\, B_m), L^{2}(\partial\, B_m) \right)}$.  Hence
\begin{eqnarray*}\label{Q2dbl}
 |Q_m|&=&\left|\int_{ \partial D_m} \frac{\partial\,U^{\sigma_{m}} (s)}{\partial \nu_m(s)} ds\right|\nonumber\\
    &\leq& \| 1  \|_{H^{1}(\partial D_m)}\left\| \frac{\partial\,U^{\sigma_{m}} (s)}{\partial \nu_m(s)} \right\|_{H^{-1}(\partial D_m)}\nonumber\\
    &\leq& \| 1  \|_{L^2(\partial D_m)}\grave{\mathcal{C}}(\kappa)\,
\|\varLambda_{B_m}\|_{\mathcal{L}\left(H^1(\partial\, B_m), L^{2}(\partial\, B_m) \right)}\nonumber\\
    &\leq&\epsilon\, |\partial \c{B}|\,\grave{\mathrm{C}} 
\,\|\varLambda_{B_m}\|_{\mathcal{L}\left(H^1(\partial\, B_m), L^{2}(\partial\, B_m) \right)}.
\end{eqnarray*}
\end{proofn}
\begin{proposition}\label{farfldthmdbl} The far-field pattern $U^\infty$ corresponding to the scattered solution of the problem (\ref{acimpoenetrable}-\ref{radiationc}) 
has the following asymptotic expansion
\begin{equation}\label{x oustdie1 D_m-dbl}
U^\infty(\hat{x})=\sum_{m=1}^{M}e^{-i\kappa\hat{x}\cdot z_{m}}[Q_m+O(\kappa\,a^2)],
\end{equation}
with $Q_m$ given by (\ref{defofQmdbl}), if $\kappa\,a<1$ where $O(\kappa\,a^2)\,\leq\,\grave{C}\kappa\,a^2$ and 
$\grave{C}:=\frac{|\partial \c{B}|\grave{\mathrm{C}}\|\varLambda_{B_m}\|_{\mathcal{L}\left(H^1(\partial B_m), L^{2}(\partial B_m) \right)} }{\max\limits_{1\leq m \leq M} diam(B_m)}$.
\end{proposition}
\begin{proofn}{\it{of Proposition \ref{farfldthmdbl}}.}
From \eqref{qcimprequiredfrm3dbl}, we know that 
\begin{eqnarray*}
 U^{s}(x)&=&\sum_{m=1}^{M}\int_{\partial D_m}\Phi_\kappa(x,s)\frac{\partial\,U^{\sigma_{m}} (s)}{\partial \nu_m(s)}ds,\,\text{ for }x\in\mathbb{R}^{3}\backslash\left(\mathop{\cup}\limits_{m=1}^M \bar{D}_m\right). \nonumber
\end{eqnarray*}
Hence
\begin{eqnarray}\label{xfarawayimpntdbl}
U^{\infty}(\hat{x})&=&\sum_{m=1}^{M}\int_{\partial D_m}e^{-i\kappa\hat{x}\cdot s}\frac{\partial\,U^{\sigma_{m}} (s)}{\partial \nu_m(s)}ds\nonumber\\
 &=&\sum_{m=1}^{M}\left(e^{-i\kappa\hat{x}\cdot\,z_{m}}Q_m+\int_{\partial D_m}[e^{-i\kappa\hat{x}\cdot\,s}-e^{-i\kappa\hat{x}\cdot\,z_{m}}]\frac{\partial\,U^{\sigma_{m}} (s)}{\partial \nu_m(s)}ds\right).
\end{eqnarray}
As in Lemma \ref{Qmestbigodbl},  we have from Lemma \ref{thmestmdhoudhonu};
\begin{eqnarray}\label{estimationofintsigmadbl}
\int_{\partial D_m}\left|\frac{\partial\,U^{\sigma_{m}} (s)}{\partial \nu_m(s)}\right| ds
&\leq&\grave{C}a,\,\forall\,m=1,2,\dots,M
\end{eqnarray}
with $\grave{C}:=\frac{\grave{\tilde{c}}}{\max\limits_{1\leq m \leq M} diam(B_m)}=\frac{|\partial\c{B}|\grave{\mathrm{C}}\|\varLambda_{B_m}\|_{\mathcal{L}\left(H^1(\partial\, B_m), L^{2}(\partial\, B_m) \right)} }{\max\limits_{1\leq m \leq M} diam(B_m)}$.
It gives us the following estimate in the similar lines of \eqref{acestimateforexponentsdif-pre};
\begin{eqnarray}\label{acestimateforexponentsdif-predbl}
\left|\int_{\partial D_m}[e^{-i\kappa\hat{x}\cdot\,s}-e^{-i\kappa\hat{x}\cdot\,z_{m}}]\frac{\partial\,U^{\sigma_{m}} (s)}{\partial \nu_m(s)}ds\right|
&\leq&\frac{1}{2}\grave{C}\kappa\,a^2\frac{1}{1-\frac{1}{2}\kappa\,a},\,\text{if}\,a<\frac{2}{\kappa_{\max}}\left(\leq\frac{2}{\kappa}\right)
\end{eqnarray}
which means
\begin{eqnarray}\label{acestimateforexponentsdifdbl}
\int_{\partial D_m}[e^{-i\kappa\hat{x}\cdot\,s}-e^{-i\kappa\hat{x}\cdot\,z_{m}}]\frac{\partial\,U^{\sigma_{m}} (s)}{\partial \nu_m(s)}ds&\leq&\grave{C}\kappa\,a^2,\,\text{for}\,a\leq\frac{1}{\kappa_{\max}}.
\end{eqnarray}
Now substitution of \eqref{acestimateforexponentsdifdbl} in \eqref{xfarawayimpntdbl} gives the required result \eqref{x oustdie1 D_m-dbl}.
\end{proofn}

\noindent
Let us derive a formula for $Q_m$. For $s_m\in \partial D_m$, using the Dirichlet boundary condition \eqref{acgoverningsupport}, we have
\begin{eqnarray}\label{Q_mintdbl}
 0&=&U^{t}(s_m)=U^{i}(s_m)+\sum_{j=1}^{M}\int_{\partial D_j}\Phi_\kappa(s_m,s)\frac{\partial\,U^{\sigma_{j}} (s)}{\partial \nu_j(s)}ds\nonumber \\
&=&U^{i}(s_m)+\sum_{\substack{j=1 \\ j\neq m}}^{M}\left(\Phi_\kappa(s_m,z_j)Q_j+\int_{\partial D_j}\hspace{-.1cm}[\Phi_\kappa(s_m,s)-\Phi_\kappa(s_m,z_j)]\frac{\partial\,U^{\sigma_{j}} (s)}{\partial \nu_j(s)}ds\right)+\int_{\partial D_m}\hspace{-.2cm}\Phi_\kappa(s_m,s)\frac{\partial\,U^{\sigma_{m}} (s)}{\partial \nu_m(s)}ds.\nonumber\\
\end{eqnarray}
Now let us estimate $\int_{\partial D_j}[\Phi_\kappa(s_m,s)-\Phi_\kappa(s_m,z_j)]\frac{\partial\,U^{\sigma_{j}} (s)}{\partial \nu_j(s)} (s)ds$ for $j\neq\,m$. 

For $m,j=1,\dots,M$, and $j\neq\,m$, by making use of \eqref{taylorphifar}, \eqref{estimatofRxs} and \eqref{estimationofintsigmadbl} we obtain the 
following in a similar way to \eqref{estphismzjdif};
\begin{eqnarray}\label{estphismzjdifdbl}
 \left|\int_{\partial D_j}[\Phi_\kappa(s_m,s)-\Phi_\kappa(s_m,z_j)]\frac{\partial\,U^{\sigma_{j}} (s)}{\partial \nu_j(s)}ds\right|
                                                               &<\,&\grave{C}\frac{a}{d}\left(\kappa+\frac{1}{d}\right)a.
\end{eqnarray}
\noindent
Then \eqref{Q_mintdbl} can be written as
\begin{equation}\label{qcimsurfacefrmdbl}
\begin{split}
 0=U^{i}(s_m)+\sum_{\substack{j=1 \\ j\neq m}}^{M}\Phi_\kappa(s_m,z_j)Q_j&+O\left((M-1)\left(\frac{\kappa a^2}{d}+\frac{a^2}{d^2}\right)\right) \\
 &+\int_{\partial D_m}\Phi_0(s_m,s)\left[1+(e^{i\kappa|s_m-s|}-1)\right]
                                                                        \frac{\partial\,U^{\sigma_{m}} (s)}{\partial \nu_m(s)}ds.
\end{split}
\end{equation}
By using the Taylor series expansions of the exponential term $e^{i\kappa|s_m-s|}$, the above can also be written as,
\begin{eqnarray}\label{qcimsurfacefrm1dbl}
\hspace{-.7cm}\int_{\partial D_m}\hspace{-.2cm}\Phi_0(s_m,s)\frac{\partial\,U^{\sigma_{m}} (s)}{\partial \nu_m(s)}ds+O(\kappa a)&=&-U^{i}(s_m)-\sum_{\substack{j=1 \\ j\neq m}}^{M}\Phi_\kappa(s_m,z_j)Q_j+O\left((M-1)\left(\frac{\kappa a^2}{d}+\frac{a^2}{d^2}\right)\right).
\end{eqnarray}
\begin{indeed}
for $m=1,\dots,M$, we have the below estimate in the similar lines of \eqref{phinotexp-1est};
\begin{eqnarray}\label{phinotexp-1estdbl}
\left|\int_{\partial D_m}\Phi_0(s_m,s)\left(e^{i\kappa|s_m-s|}-1\right)\frac{\partial\,U^{\sigma_{m}} (s)}{\partial \nu_m(s)}ds\right|
                             &\leq&\grave{C} \kappa\,a, \text{ for }a\leq\frac{1}{\kappa_{\max}}.
\end{eqnarray}
\end{indeed}
Define $U_m:=\int_{\partial D_m}\Phi_0(s_m,s)\frac{\partial\,U^{\sigma_{m}} (s)}{\partial \nu_m(s)}ds,\,s_m\in\,\partial D_m$. Then \eqref{qcimsurfacefrm1dbl} can be written as
\begin{eqnarray}\label{qcimsurfacefrm2dbl}
U_m&=&-U^{i}(s_m)-\sum_{\substack{j=1 \\ j\neq m}}^{M}\Phi_\kappa(s_m,z_j)Q_j+O(\kappa a)+O\left((M-1)\left(\frac{\kappa a^2}{d}+\frac{a^2}{d^2}\right)\right).
\end{eqnarray}
We set
\begin{eqnarray}\label{def-acoust-dbl-revise-defntn-barum}
\bar{U}_m:=-U^{i}(z_m)-\sum\limits_{\substack{j=1 \\ j\neq m}}^{M}\Phi_\kappa(z_m,z_j)Q_j, \,\mbox{for}\,m=1,\dots,M.
\end{eqnarray}
For $m=1,\dots,M$, let $\bar{\sigma}_m\in L^2(\partial D_m)$ be the solutions of following the integral equation;
\begin{eqnarray}\label{nbc1dbl-revise}
\frac{\sigma_m(s)}{2}+\int_{\partial D_m}\frac{\partial\Phi_0(x,s)}{\partial \nu_m(s)}\sigma_{m} (s)ds=\bar{U}_m \,\mbox{on}\, \partial D_m.
\end{eqnarray} 
Remark here that the left hand side of \eqref{nbc1dbl-revise} is the trace, on $\partial{D_m}$, of the double layer potential $\int_{\partial D_m}\frac{\partial\Phi_0(x,s)}{\partial \nu_m(s)}\sigma_{m} (s)ds$, 
$x\in\mathbb{R}^3\backslash\bar{D}_m$. Dealing in the similar way as we derived \eqref{revised-greensform-DtN-acoustic-small}, we obtain
\begin{eqnarray}\label{nbc1dbl-1-revise}
 \int_{\partial D_m}\frac{\partial\Phi_0(x,s)}{\partial \nu_m(s)}\sigma_{m} (s)ds=\int_{\partial D_m}\Phi_0(x,s)\frac{\partial\,U^{\bar{\sigma}_{m}} (s)}{\partial \nu_m(s)}ds,
\end{eqnarray}
with $U^{\bar{\sigma}_{m}}$ are the solutions of \eqref{acimpenetrableUsigma} replacing the wave number $\kappa$ by zero. As single layer potential is continuous up to the boundary, 
combining \eqref{nbc1dbl-revise} and \eqref{nbc1dbl-1-revise}, we deduce that
and 
the constant potentials $\bar{U}_m,\,m=1,\dots,M$ satisfy,
\begin{equation}\label{barqcimsurfacefrm1dbl}
\int_{\partial D_m}\Phi_0(s_m,s)\frac{\partial\,U^{\bar{\sigma}_{m}} (s)}{\partial \nu_m(s)} ds=\bar{U}_m,\,s_m\in\,\partial D_m.
\end{equation}
Then, the total charge on the surface $\partial D_m$ is given as
$$\bar{Q}_m:=\int_{\partial D_m}\frac{\partial\,U^{\bar{\sigma}_{m}} (s)}{\partial \nu_m(s)}ds.$$
We set the electrical capacitance $\bar{C}_m$ for $1\leq\,m\leq\,M$ as
$$\bar{C}_m=\frac{\bar{Q}_m}{\bar{U}_m}.$$

\par Following in the similar lines of the proofs of Lemma \ref{lemmadifssbQQbCCb}, Lemma \ref{lemmadifssbQQbCCb1} and Proposition \ref{fracqfracc-ac} 
concerning the single layer potentials, we can prove the following results in the case of double layer potentials.
 
\begin{itemize}
\item
The consecutive pairwise difference between $\frac{\partial\,U^{\sigma_{m}}}{\partial \nu_m}$, $\frac{\partial\,U^{\bar{\sigma}_{m}}}{\partial \nu_m}$, $Q_m$, $\bar{Q}_m$ have the following behaviour.
 \begin{eqnarray}
  \left\|\frac{\partial\,U^{\sigma_{m}} }{\partial \nu_m}-\frac{\partial\,U^{\bar{\sigma}_{m}} }{\partial \nu_m}\right\|_{H^{-1}(\partial D_m)}&=&O\left(\kappa a+(M-1)\left(\frac{\kappa a^2}{d}+\frac{a^2}{d^2}\right)\right),\label{difssbdbl}\\
  Q_m-\bar{Q}_m&=&O\left(\kappa a^2+(M-1)\left(\frac{\kappa a^3}{d}+\frac{a^3}{d^2}\right)\right).\label{difQmddbl}
 \end{eqnarray}
\item
For every $1\leq m\leq M$, we have
\begin{eqnarray}\label{asymptotCapdbl}
 \bar{C}_m\,=\,\frac{\bar{C}_{B_m}}{\max\limits_{1\leq m \leq M} diam(B_m)}a
 & \mbox{ and }&
\bar{Q}_m\,=\,\frac{\bar{Q}_{B_m}}{\max\limits_{1\leq m \leq M} diam(B_m)}a.
 \end{eqnarray}
\item For $m=1,2,\dots,M$, the total charge $\bar{Q}_m$ on each surface $\partial D_m$ of the small scatterer $D_m$ can be calculated from the algebraic system 
  \begin{eqnarray}\label{fracqfracdbl}
 \frac{\bar{Q}_m}{\bar{C}_m} &=&-U^{i}(z_m)-\sum_{\substack{j=1 \\ j\neq m}}^{M}\bar{C}_j \Phi_\kappa(z_m,z_j)\frac{\bar{Q}_j}{\bar{C}_j},
  \end{eqnarray}
  which is valid with an error of order $O\left((M-1)\frac{\kappa a^2}{d}+(M-1)^2\left(\frac{\kappa a^3}{d^2}+\frac{a^3}{d^3}\right)\right)$.
\end{itemize}

\subsection{The algebraic system}\label{sec-algebraicsys-small}
Define the  algebraic system,
  \begin{eqnarray}\label{fracqcfrac}
 \frac{\tilde{Q}_m}{\bar{C}_m} &:=&-U^{i}(z_m)-\sum_{\substack{j=1 \\ j\neq m}}^{M}\bar{C}_j \Phi_\kappa(z_m,z_j)\frac{\tilde{Q}_j}{\bar{C}_j},
  \end{eqnarray}
 for all $m=1,2,\dots,M$. 
It can be written in a compact form as
\begin{equation}\label{compacfrm1}
 \mathbf{B}\tilde{Q}=\mathrm{U}^I,
\end{equation}
\noindent
where $\tilde{Q},\mathrm{U}^I \in \mathbb{C}^{M\times 1}\mbox{ and } \mathbf{B}\in\mathbb{C}^{M\times M}$ are defined as
\begin{eqnarray}
\mathbf{B}:=\left(\begin{array}{ccccc}
   -\frac{1}{\bar{C}_1} &-\Phi_\kappa(z_1,z_2)&-\Phi_\kappa(z_1,z_3)&\cdots&-\Phi_\kappa(z_1,z_M)\\
-\Phi_\kappa(z_2,z_1)&-\frac{1}{\bar{C}_2}&-\Phi_\kappa(z_2,z_3)&\cdots&-\Phi_\kappa(z_2,z_M)\\
 \cdots&\cdots&\cdots&\cdots&\cdots\\
-\Phi_\kappa(z_M,z_1)&-\Phi_\kappa(z_M,z_2)&\cdots&-\Phi_\kappa(z_M,z_{M-1}) &-\frac{1}{\bar{C}_M}
   \end{array}\right),\label{mainmatrix-acoustic-small}\\
\nonumber\\
 \tilde{Q}:=\left(\begin{array}{cccc}
    \tilde{Q}_1 & \tilde{Q}_2 & \ldots  & \tilde{Q}_M
   \end{array}\right)^\top \text{ and } 
\mathrm{U}^I:=\left(\begin{array}{cccc}
     U^i(z_1) & U^i(z_2)& \ldots &  U^i(z_M)
   \end{array}\right)^\top.\label{coefficient-and-incidentvectors-acoustic-small}
\end{eqnarray}
The above linear algebraic system is solvable for $\tilde{Q}_j,~1\leq j\leq M$, when the matrix $\mathbf{B}$ is invertible. Next we discuss the possibilities of this invertibility,

$\bullet$ If
\begin{equation}\label{Diagdmntcndt}
\max\limits_{1\leq m\leq M}\sum_{\substack{j=1\\ j\neq m}}^{M}\frac{\bar{C}_m}{|z_m-z_j|}< 4\pi, 
\end{equation}
  then $\mathbf{B}$ satisfies the diagonally dominant condition and hence the system is solvable. Since, $\bar{C}_m\,:=\,\epsilon\bar{C}_{B_m}\,=\,\frac{\bar{C}_{B_m}}{\max\limits_{1\leq m \leq M} diam(B_m)}a$ 
then \eqref{Diagdmntcndt} is valid if $(M-1)\frac{a}{d}<4\pi\left(\max\limits_{1\leq m \leq M}\bar{C}_{B_m}\right)^{-1}
\max\limits_{1\leq m \leq M} diam(B_m)$, where $\bar{C}_{B_m}$ depends only the
 Lipschitz character of $B_m$.

$\bullet$ If $a<cd$ for some constant $c$, i.e. precisely if $\max\limits_{1\leq\,m\leq\,M}\bar{C}_m<\frac{5\pi}{3}d$ and 
$t:=\min\limits_{j\neq\,m,1\leq\,j,m\leq\,M}\cos(\kappa|z_m-z_j|) \geq 0$, 
then the matrix $\mathbf{B}$ is invertible. Remark that $c\leq\frac{5\pi}{3}\left(\max\limits_{1\leq m \leq M}\bar{C}_{B_m}\right)^{-1}\max\limits_{1\leq m \leq M} diam(B_m)$.
We state this invertibility property in the following lemma.
\begin{lemma}\label{Mazyawrkthm}
If $\max_{1\leq\,m\leq\,M}\bar{C}_m<\frac{5\pi}{3}d$ and $t:=\min\limits_{j\neq\,m,1\leq\,j,m\leq\,M}\cos(\kappa|z_m-z_j|) \geq 0$, then the matrix $\mathbf{B}$ 
is invertible and the solution vector $\tilde{Q}$ of \eqref{compacfrm1} satisfies the estimate
\begin{equation}\label{mazya-fnlinvert-small-ac-2}
\begin{split}
 \sum_{m=1}^{M}|\tilde{Q}_m|^{2}\bar{C}_m^{-1}
\leq4\left(1-\frac{3t}{5\pi\,d}\right.&\left.\max\limits_{1\leq m \leq M}\bar{C}_m\right)^{-2}\sum_{m=1}^{M}\left|U^i(z_m)\right|^2\bar{C}_m.
\end{split}
\end{equation}

\end{lemma}
\begin{proofn}{\it{of Lemma \ref{Mazyawrkthm}}.} The idea of the proof of this lemma is given by Maz'ya and Movchan in \cite{M-M:MathNach2010}
for the case where $\Phi$ is the Green's function of the Dirichlet-Laplacian in a bounded domain. We adapt their
argument to the case of Helmholtz, $\kappa \neq 0$, on the whole space. For the reader's convenience, we give the detailed proof in the appendix.
\end{proofn}

\subsection{End of the proof of Theorem \ref{Maintheorem-ac-small-sing}}\label{sec-mainthrmproof-small}
\begin{proofe}
We can rewrite the equation \eqref{mazya-fnlinvert-small-ac-2} using norm inequalities as
\begin{equation}\label{mazya-fnlinvert-small-ac-3}
\begin{split}
 \sum_{m=1}^{M}|\tilde{Q}_m|
\leq2\left(1-\frac{3t}{5\pi\,d}\right.&\left.\max\limits_{1\leq m \leq M}\bar{C}_m\right)^{-1}M\max\limits_{1\leq m \leq M}|\bar{C}_m|\max\limits_{1\leq m \leq M}\left|U^i(z_m)\right|
\end{split}
\end{equation}
with $t:=\min\limits_{j\neq\,m,1\leq\,j,m\leq\,M}\cos(\kappa|z_m-z_j|)$.
 The difference between \eqref{fracqfrac}/\eqref{fracqfracdbl} and \eqref{fracqcfrac} produce the following
  \begin{eqnarray}\label{qcdiftilde}
   \frac{\bar{Q}_m-\tilde{Q}_m}{\bar{C}_m} &=&-\sum_{\substack{j=1 \\ j\neq m}}^{M} \Phi_\kappa(z_m,z_j)\left(\bar{Q}_j-\tilde{Q}_j\right)+O\left((M-1)\frac{\kappa a^2}{d}+(M-1)^2\left(\frac{\kappa a^3}{d^2}+\frac{a^3}{d^3}\right)\right),
  \end{eqnarray}
for $m=1,\dots,M$. Comparing the above system of equations \eqref{qcdiftilde} with \eqref{fracqcfrac} and by making use of the estimate \eqref{mazya-fnlinvert-small-ac-3}, we obtain
 \begin{eqnarray}\label{unncmaybe}
\sum_{m=1}^{M}({\bar{Q}_m-\tilde{Q}_m})&=&O\left(M(M-1)\frac{\kappa a^3}{d}+M(M-1)^2\left(\frac{\kappa a^4}{d^2}+\frac{a^4}{d^3}\right)\right).
\end{eqnarray}
We can evaluate the $\tilde{Q}_m$'s from the algebraic system \eqref{fracqcfrac}. Hence, by using \eqref{difQmd}/\eqref{difQmddbl} 
and \eqref{unncmaybe} in 
\eqref{x oustdie1 D_m}/\eqref{x oustdie1 D_m-dbl} we can represent the far-field pattern in terms of 
$\tilde{Q}_m$ as below:
\begin{eqnarray}\label{x oustdie1 D_m1}
\hspace{-.7cm}U^{\infty}(\hat{x})&=&\sum_{m=1}^{M}[e^{-i\kappa\hat{x}\cdot z_{m}}Q_m+O(\kappa\,a^2)]\nonumber\\
                   &=&\sum_{m=1}^{M}\left[e^{-i\kappa\hat{x}\cdot z_{m}}[\tilde{Q}_m+(Q_m-\bar{Q}_m)+(\bar{Q}_m-\tilde{Q}_m)]+O(\kappa\,a^2)\right]\nonumber\\
 &=&\sum_{m=1}^{M}e^{-i\kappa\hat{x}\cdot z_{m}}\tilde{Q}_m+O\left(M\kappa a^2+M(M-1)\left(\frac{\kappa a^3}{d}+\frac{a^3}{d^2}\right)+M(M-1)^2\frac{a}{d}\left(\frac{\kappa a^3}{d}+\frac{a^3}{d^2}\right)\right).
\end{eqnarray}
Hence Theorem \ref{Maintheorem-ac-small-sing} is proved by setting $\bar{\sigma}_m:=\frac{\bar{\sigma}_m}{\bar{U}_m}$ 
as the surface density which defines $\tilde{Q}_m$. Finally, let us stress that
\begin{enumerate}
  \item The constant 
$\grave{c}=\left[\frac{(\bar{c}+1)}{4\pi}|\partial\c{B}|\max\limits_{1\leq m \leq M}\grave{C}_{6m}\right]^{-\frac{1}{2}}$
 appearing in Proposition \ref{normofsigmastmtdbl} 
will serve our purpose in Theorem \ref{Maintheorem-ac-small-sing} by defining $c_0:=\grave{c}\, \max\limits_{1\leq m \leq M } diam (B_m)$ when we use double 
layer potentials. Also observe that
the constant $c$  appearing in Proposition \ref{normofsigmastmt} 
can be used to prove Theorem \ref{Maintheorem-ac-small-sing} if we replace  the second condition in 
\eqref{conditions} by the stronger one $(M-1)\frac{a}{d^2}\leq c_0$ with $c_0:=c\,\max\limits_{ m } diam (B_m)$.

\item The coefficients $\bar{\sigma}_m {\bar{U}_m}^{-1}, \, \tilde{Q}_m,\, \bar{C}_m$ play the roles of $\sigma_m,\, Q_m,\, C_m$ respectively in Theorem \ref{Maintheorem-ac-small-sing}.
\item The constant appearing in $O\left(M\kappa a^2+M(M-1)\left(\frac{\kappa a^3}{d}+\frac{a^3}{d^2}\right)+M(M-1)^2\frac{a}{d}\left(\frac{\kappa a^3}{d}+\frac{a^3}{d^2}\right)\right)$ is 
 $\frac{\max\limits_{1\leq m \leq M}\left|\left|{\mathcal{S}^{i_{\kappa}}_{B_m}}^{-1}\right|\right|_{\mathcal{L}\left(H^1(\partial B_m), L^2(\partial B_m) \right)}|\partial\c{B}|}{\max\limits_{1\leq m \leq M} diam(B_m)}(\grave{C}+1)
 \max\left\{
 1+\frac{\grave{C}}{(\grave{C}+1)}\frac{\max\limits_{1\leq m \leq M} diam(B_m)}{\max\limits_{1\leq m \leq M}\left|\left|{\mathcal{S}^{i_{\kappa}}_{B_m}}^{-1}\right|\right|_{\mathcal{L}\left(H^1(\partial B_m), L^2(\partial B_m) \right)}|\partial\c{B}|},\right.$
 $\left.1+\frac{\max\limits_{1\leq m \leq M} \bar{C}_{B_m}}{4\pi\max\limits_{1\leq m \leq M} diam(B_m)}\right\}$. 
The constants $|\partial\c{B}|$ ,
and $\grave{C}$ are defined in
 Lemma \ref{propsmjsmmest} and Proposition \ref{farfldthmdbl}.
 \item The constant $a_0$ appearing in \eqref{conditions} of Theorem \ref{Maintheorem-ac-small-sing} is the minimum among
$$\frac{1}{\kappa_{\max}}\min \left\{1, \sqrt[3]{\frac{4\pi}{3}}{\rm j}_{1/2,1}\right\}, ~~ 
 \frac{2\pi\max\limits_{1\leq m \leq M} diam(B_m)}
{\kappa_{\max}(1+\kappa_{\max})|\partial\c{B}|\max\limits_{1\leq m \leq M}\left\|{\mathcal{S}^{i_{\kappa}}_{B_m}}^{-1}\right\|_{\mathcal{L}\left(H^1(\partial B_m), L^2(\partial B_m) \right)}}
$$
and  
 $\frac{2\sqrt{\pi}\max\limits_{1\leq m \leq M} diam(B_m)}{\kappa_{\max}\left(|\partial\c{B}|\max\limits_{1\leq m \leq M}\left\|\left(\frac{1}{2}I+{\mathcal{D}_{B_m}^{i_{\kappa}}}\right)^{-1}\right\|_{\mathcal{L}\left(L^2(\partial B_m), L^2(\partial B_m) \right)}\right)^\frac{1}{2}}$.

 \item The constant $c_1$ appearing in \eqref{invertibilityconditionsmainthm} of Theorem \ref{Maintheorem-ac-small-sing} is 
 $\frac{5\pi}{3}\frac{\max\limits_{1\leq m \leq M} diam(B_m)}{\max\limits_{1\leq m \leq M}\bar{C}_{B_m}}$ and it follows from Lemma \ref{lemmadifssbQQbCCb1} and Lemma \ref{Mazyawrkthm}.
\end{enumerate}
\end{proofe}
From the last points, we see that the constants appearing in Theorem \ref{Maintheorem-ac-small-sing} depend only on  
 $d_{\max}$, $\kappa_{\max}$ and $B_m$'s through their diameters, capacitances and the norms of the boundary operators ${\mathcal{S}^{i_{\kappa}}_{B_m}}^{-1}:
H^1(\partial B_m) \rightarrow L^2(\partial B_m)$, $\left(\frac{1}{2}I+
{\mathcal{D}_{B_m}^{i_{\kappa}}}\right)^{-1}:L^2(\partial B_m)\rightarrow L^2(\partial B_m)$ and
$\varLambda_{B_m}: H^1(\partial\, B_m) \rightarrow L^{2}(\partial\, B_m)$. In the following remark, 
we show how the dependency on $B_m$'s are actually only through their Lipschitz character.

\begin{remark}\label{Lipschitz-character}
\begin{enumerate}
\item Let $B$ be a bounded, simply connected and Lipschitz domain in $\mathbb{R}^3$. The quantities
$\left\|\left(\frac{1}{2}I+
{\mathcal{D}_{B}^{i_{\kappa}}}\right)^{-1}\right\|_{\mathcal{L}\left(L^2(\partial B), L^2(\partial B) \right)}$
and $\|\varLambda_{B}\|_{\mathcal{L}\left(H^1(\partial\, B), L^{2}(\partial\, B) \right)} $
depend only on the Lipschitz character of $B$. Indeed, we first remark that $\varLambda_{B}=\left(\frac{1}{2}I+
{\mathcal{D}_{B}^{i_{\kappa}}}\right) (S^{i_\kappa})^{-1}$ hence 
$\|\varLambda_{B}\|_{\mathcal{L}\left(H^1(\partial\, B), L^{2}(\partial\, B) \right)}  
\leq \left\|\left(\frac{1}{2}I+
{\mathcal{D}_{B}^{i_{\kappa}}}\right)\right\|_{\mathcal{L}\left(L^2(\partial B), L^2(\partial B) \right)} 
\Vert (S^{i_\kappa})^{-1}\Vert_{\mathcal{L}\left(H^1(\partial B), L^2(\partial B) \right)}$. The operators $\frac{1}{2}I+
{\mathcal{D}_{B}^{i_{\kappa}}}$ and $S^{i_\kappa}$ are isomorphism in the mentioned spaces and their norms 
depend only on the Lipschitz character of $B$, see Section 1 in \cite{A-K:2007} for instance. 
\item The capacitance $C_B$ of a bounded, connected and Lipschitz domain $B$ of $\mathbb{R}^3$ depends only on the Lipschitz character of $B$.
Indeed, since $C_B:=\int_{\partial B}\sigma(s)ds$ where $\int_{\partial B}\frac{\sigma (t)}{4\pi\vert s-t\vert}dt=1, \;~ s \in \partial B$, 
then from the invertibility of the single layer potential $S^{i_\kappa}: L^2(\partial B) \rightarrow H^1(\partial B)$, we deduce that
\begin{equation}\label{r-estimate}
C_B \leq \vert \partial B \vert^{\frac{1}{2}} \Vert \sigma \Vert_{L^2(\partial B)} \leq \vert \partial B \vert^{\frac{1}{2}} 
\Vert (S^{i_\kappa})^{-1}\Vert_{\mathcal{L}\left(H^1(\partial B), L^2(\partial B) \right)} \vert \partial B \vert^{\frac{1}{2}}=
\Vert (S^{i_\kappa})^{-1}\Vert_{\mathcal{L}\left(H^1(\partial B), L^2(\partial B) \right)} \vert \partial B \vert
\end{equation}
On the other hand, we recall the following lower estimate, see Theorem 3.1 in \cite{MRAMAG-book1} for instance,
\begin{equation}\label{C-J}
C_B \geq \frac{4\pi \vert \partial B \vert^2}{J} 
\end{equation}
where $J:=\int_{\partial B}\int_{\partial B}\frac{1}{\vert s-t\vert}ds dt$. Remark that $J=4\pi\int_{\partial B} S^{i_\kappa}(1)(s)ds$. Hence
$$ J\leq  4 \pi \vert \partial B\vert^{\frac{1}{2}} \Vert S^{i_\kappa}\Vert_{\mathcal{L}\left(L^2(\partial B), H^1(\partial B)\right)}
\Vert 1\Vert_{H^1(\partial B)}\leq 4 \pi \Vert S^{i_\kappa}\Vert_{\mathcal{L}\left(L^2(\partial B), H^1(\partial B)\right)} \vert
\partial B\vert
$$
and using (\ref{C-J}) we obtain the lower bound
\begin{equation}\label{l-estimate}
C_{B}\geq \Vert S^{i_\kappa}\Vert^{-1}_{\mathcal{L}\left(L^2(\partial B), H^1(\partial B) \right)} \vert \partial B\vert.
\end{equation}
Finally combining (\ref{r-estimate}) and (\ref{l-estimate}), we derive the estimated
\begin{equation}\label{Capacitance-size-estimate}
 \Vert S^{i_\kappa}\Vert^{-1}_{\mathcal{L}\left(L^2(\partial B), H^1(\partial B) \right)} \vert \partial B\vert \leq C_B \leq
\Vert (S^{i_\kappa})^{-1}\Vert_{\mathcal{L}\left(H^1(\partial B), L^2(\partial B) \right)} \vert \partial B \vert
\end{equation}
which shows, in particular, that the capacitance of $B$ depends only on the Lipschitz character $B$, knowing that also $\vert \partial B\vert$
can be estimated only with the Lipschitz character of $B$, see for instance the observations in Remark 2.5 in \cite{A-M-R:2002}.
\end{enumerate}
Related estimates in terms of the Lipschitz character for a different problem can be found in \cite{AH-SJK:AAM2003}.
\end{remark}
\subsection{Proof of remark \ref{corMaintheorem-ac-small-sing}} 
\begin{proofe}
 For $m=1,\dots,M$ fixed, we distinguish between the obstacles $D_j$, $j\neq\,m$ which are near to $D_m$ from the ones which are far from $D_m$ as follows.  Let $\Omega_m$, $1\leq\,m\leq\,M$ be the balls
of center $z_m$ and of radius $(\frac{a}{2}+d^\alpha)$ with $0<\alpha\leq1$. The bodies lying in $\Omega_m$ will fall into the category, $N_m$, of near by obstacles 
 and the others into the category, $F_m$, of far obstacles to $D_m$. Since the obstacles $D_m$ are balls with same diameter, the number of obstacles
 near by $D_m$ will not exceed $\left(\frac{a+2d^\alpha}{a+d}\right)^3$   $\left[=\frac{\frac{4}{3}\pi\left((a+2d^\alpha)/2\right)^3}{\frac{4}{3}\pi\left((a+d)/2\right)^3}\right]$.\\
 \par With this observation, instead of \eqref{x oustdie1 D_m farmain}\slash\eqref{x oustdie1 D_m1}, the far field will have the asymptotic expansion \eqref{x oustdie1 D_m farmain-recent-near}.
 \begin{indeed}\\
 $\bullet$ For the bodies $D_j\in\,N_m,j\neq{m}$ we have the estimate \eqref{estphismzjdif}
  but for the bodies $D_j\in F_m$, we obtain the following estimate 
\begin{eqnarray}\label{estphismzjdifnear}
 \left|\int_{\partial D_j}[\Phi_\kappa(s_m,s)-\Phi_\kappa(s_m,z_j)]\sigma_{j} (s)ds\right|
                                                                &<\,&C\frac{a}{d^\alpha}\left(\kappa+\frac{1}{d^\alpha}\right)a.
\end{eqnarray}
$\bullet$ Due to the estimates \eqref{estphismzjdif} and \eqref{estphismzjdifnear}, while representing the scattered field in terms of single layer potential, 
corresponding changes will take place in \eqref{qcimsurfacefrm}, \eqref{qcimsurfacefrm1}, 
 \eqref{qcimsurfacefrm2}, \eqref{difssb}, \eqref{difQmd} and in \eqref{fracqfrac} which inturn modify \eqref{qcdiftilde}, \eqref{unncmaybe} and hence the asymptotic expansion \eqref{x oustdie1 D_m farmain}\slash\eqref{x oustdie1 D_m1} modify as follows
 \begin{eqnarray}\label{x oustdie1 D_m farmain-recent}
U^\infty(\hat{x},\theta)&=&\sum_{m=1}^{M}e^{-i\kappa\hat{x}\cdot z_m}Q_m+
O\left(M\kappa a^2+M(M-1)\left(\frac{\kappa a^3}{d^\alpha}+\frac{a^3}{d^{2\alpha}}\right)+M\left(\frac{a+2d^\alpha}{a+d}\right)^3\left(\frac{\kappa a^3}{d}+\frac{a^3}{d^2}\right)\right.\nonumber\\
&&\left.+M(M-1)^2\frac{a}{d^\alpha}\left(\frac{\kappa a^3}{d^\alpha}+\frac{a^3}{d^{2\alpha}}\right)+M(M-1)\frac{a}{d^\alpha}\left(\frac{a+2d^\alpha}{a+d}\right)^3\left(\frac{\kappa a^3}{d}+\frac{a^3}{d^{2}}\right)\right.\nonumber\\
&&\left.+M(M-1)\frac{a}{d}\left(\frac{a+2d^\alpha}{a+d}\right)^3\left(\frac{\kappa a^3}{d^\alpha}+\frac{a^3}{d^{2\alpha}}\right)+M\left(\frac{a+2d^\alpha}{a+d}\right)^6\frac{a}{d}\left(\frac{\kappa a^3}{d}+\frac{a^3}{d^2}\right)\right).
 \end{eqnarray}
  The same thing can also be done while representing the scattered field in terms of double layer potential.\\
$\bullet$ Since $\kappa\leq\kappa_{\max}$, $d\leq d^\alpha,\, 0<\alpha\leq1$ and $\frac{a}{d}<\infty$, we have
$$\left(\frac{a+2d^\alpha}{a+d}\right) =d^{\alpha-1}\frac{ad^{-\alpha}+2}{ad^{-1}+1}=O(d^{\alpha-1}),$$  
which can be used to derive \eqref{x oustdie1 D_m farmain-recent-near} from \eqref{x oustdie1 D_m farmain-recent}.\\
 \end{indeed}
\end{proofe}
\section{The inverse problem}\label{The inverse problem}
From \eqref{x oustdie1 D_m1}, we can write the far-field pattern as
\begin{equation}\label{x oustdie1 D_m far}
U^\infty(\hat{x})=\sum_{m=1}^{M}e^{-i\kappa\hat{x}\cdot z_m}\tilde{Q}_m
\end{equation}
with the error of order $O\left(M\kappa a^2+M(M-1)\left(\frac{\kappa a^3}{d}+\frac{a^3}{d^2}\right)+M(M-1)^2\frac{a}{d}\left(\frac{\kappa a^3}{d}+\frac{a^3}{d^2}\right)\right)$ 
and $\tilde{Q}_m$ can be obtained from the Foldy type 
system \eqref{compacfrm1}. 
Let us denote the inverse of $\mathbf{B}$ defined in \eqref{mainmatrix-acoustic-small} by $\mathcal{B}$ and recall $\mathrm{U}^I$ defined in \eqref{coefficient-and-incidentvectors-acoustic-small}, 
then from \eqref{compacfrm1} and by keeping in mind that $U^i$ are plane incident waves given by $U^i(x):=U^i(x,\theta)=e^{i\kappa{x}\cdot\theta}$, we can state the following expression of the far field pattern;
\begin{equation}\label{x oustdie1 D_m far3}
 U^\infty(\hat{x},\theta)=\sum_{m=1}^{M}\sum_{j=1}^{M}\mathcal{B}_{mj}e^{-i\kappa\hat{x}\cdot z_m}e^{i\kappa\theta\cdot z_j}+O\left(M\kappa a^2+M(M-1)\left(\frac{\kappa a^3}{d}+\frac{a^3}{d^2}\right)+M(M-1)^2\frac{a}{d}\left(\frac{\kappa a^3}{d}+\frac{a^3}{d^2}\right)\right)
 \end{equation}
for a given incident direction $\theta$ and observation direction $\hat{x}$.
\par Our main focus is the following inverse problem.\\
\textbf{Inverse Problem :} Given the far-field pattern $U^{\infty}(\hat{x},\theta)$ for several incident and observation directions $\theta\mbox{ and } \hat{x}$,  
find the locations $z_1,z_2,\dots,z_M$ and reconstruct the  capacitances $\bar{C}_{1},\bar{C}_{2},\dots,\bar{C}_{M}$ of the small scatterers $D_1,D_2,\dots,D_M$ respectively.\\
Let us consider that all the small bodies are included in a large bounded domain $\Omega$. We distinguish the following two cases in terms of the minimum distance $d:=d(\epsilon)$:
\begin{itemize}
 \item Case 1: We have $d \ll 1$ and $(M-1)\frac{a}{d} \ll 1$. In this case, using the formula \eqref{x oustdie1 D_m far3} 
we deduce that
\begin{equation}\label{Multi-scat}
 U^\infty(\hat{x},\theta)=\sum_{m=1}^{M}\sum_{j=1}^{M}\mathcal{B}_{m, j}e^{-i\kappa\hat{x}\cdot z_m}e^{i\kappa\theta\cdot z_j}+o\left((M-1)\frac{a^2}{d}\right)
\end{equation}
as $\epsilon \rightarrow 0$. Indeed, since $\tilde{Q}_m=\sum_{j=1}^{M}\mathcal{B}_{m, j}e^{i\kappa\theta\cdot z_j}$, then
\begin{eqnarray*}
U^{\infty}(\hat{x},\theta)
 &=&\sum_{m=1}^{M}e^{-i\kappa\hat{x}\cdot z_{m}}\tilde{Q}_m+O\left(M\kappa a^2+M(M-1)\left(\frac{\kappa a^3}{d}+\frac{a^3}{d^2}\right)+M(M-1)^2\frac{a}{d}\left(\frac{\kappa a^3}{d}+\frac{a^3}{d^2}\right)\right)\\
&=&\sum_{m=1}^{M}e^{-i\kappa\hat{x}\cdot z_{m}}\sum_{j=1}^{M}\mathcal{B}_{m, j}e^{i\kappa\theta\cdot z_j}+O\left( Ma^2+M(M-1)\frac{a^3}{d^2}\right)\\
&=&\sum_{m=1}^{M}e^{-i\kappa\hat{x}\cdot z_{m}}\sum_{j=1}^{M}\mathcal{B}_{m, j}e^{i\kappa\theta\cdot z_j}+o\left((M-1)\frac{a^2}{d}\right),
\end{eqnarray*}
where the second and the third equalities holds due to the conditions $\kappa\leq\kappa_{\max}$, $d\ll1$, $(M-1)\frac{a}{d}\ll1$ and by writing $\frac{\kappa a^3}{d}=\kappa\frac{ a^3}{d^2}d$ and $Ma^2=(M-1)\frac{a^2}{d}\frac{M}{M-1}d$.

Remark that $\tilde{Q}_m:=-\bar{C}_m\, U^{i}(z_m)+\bar{C}_m\sum\limits_{\substack{j=1 \\ j\neq m}}^{M}\Phi_\kappa(z_m,z_j)\bar{C}_j\,U^{i}(z_j)+O\left((M-1)^2\frac{a^3}{d^2}\right)$, 
where $\bar{C}_m\,U^{i}(z_m)$ behaves as $a$ and $\bar{C}_m\sum_{\substack{j=1 \\ j\neq m}}^{M}\Phi_\kappa(z_m,z_j)\bar{C}_j\,U^{i}(z_j)$ behaves as $(M-1)\frac{a^2}{d}$. 
This means that \eqref{Multi-scat} reduces to 
\begin{equation*}
  U^\infty(\hat{x},\theta)=-\sum_{m=1}^{M}\bar{C}_{m}e^{-i\kappa\hat{x}\cdot z_m}e^{i\kappa\theta\cdot z_m}+\sum_{m=1}^{M}\sum_{\substack{j=1\\j\neq\,m}}^{M}\bar{C}_{m}\bar{C}_{j}\Phi_\kappa(z_m,z_{j})e^{-i\kappa\hat{x}\cdot z_m}e^{i\kappa\theta\cdot z_j}
 +o\left(M(M-1)\frac{a^2}{d}\right),
\end{equation*}
where the first term models the Born approximation and second term models the first order interaction between the scatterers. As a conclusion, when we use \eqref{Multi-scat} we compute the 
field generated by the first interaction between the collection of the scatterers $z_m,\,m=1,\dots,M$.

 \item Case 2: We assume that there exists a positive constant $d_0$ such that $d_0 \leq d$. In this case, using the formula \eqref{x oustdie1 D_m far3} 
we have
\begin{equation}\label{Born 1}
 U^\infty(\hat{x},\theta)=-\sum_{m=1}^{M}\bar{C}_{m}e^{-i\kappa\hat{x}\cdot z_m}e^{i\kappa\theta\cdot z_m}+O\left(M\kappa\,a^2\right)
\end{equation}
as $\epsilon \rightarrow 0$. In this case, we have only the Born approximation. 

\end{itemize}
Based on the formulas (\ref{Multi-scat}) and (\ref{Born 1}), we setup a MUSIC type algorithm to locate the points 
$z_j, \; j=1, ..., M$ and then estimate the sizes of the scatterers $D_j$. 
An example of distribution of the scatterers in case 1 is $M=\epsilon^{-\frac{1}{3}}$ and $d(\epsilon)=\epsilon^{\frac{1}{3}}$.
Note that in case 2 we have a lower bound on the distances between the scatterers. 
This explains why we are in the Born regime. We remark also that, in this case, $M$ is uniformly bounded since the obstacles are included in the 
bounded domain $\Omega$.
\subsection{Localisation of $D_{m}$'s via the MUSIC algorithm}\label{LocalisationviaMUSIC-smallac-sdlp}

The MUSIC algorithm is a method to determine the locations $z_{m},m=1,2,\dots,M$, of the scatterers $D_m,m=1,2,\dots,M$ from the measured far-field pattern $U^{\infty}(\hat{x},\theta)$ for 
a finite set of incidence and observation directions, i.e. $\hat{x},\theta \in \{ \theta_{j},j=1,\dots,N\}\subset\mathbb{S}^{2}$.  
We refer the reader to the monographs \cite{Ammari:2008,AH-IE-LD:MMS2005,AH-KH-KE-Lk-VMS:NM2008} and \cite{K-G:2008} for more information about this algorithm.  We follow the way presented in \cite{K-G:2008}.  
We assume that the number of scatterers is not larger than the number of incident and observation directions, i.e. $N\geq M$.  
We define the response matrix $F\in\mathbb{C}^{N\times N}$ by
\begin{equation}\label{respomatdef}
 F_{jl} := U^{\infty}(\theta_{j},\theta_{l}).
\end{equation}
From \eqref{x oustdie1 D_m far3} and \eqref{respomatdef}, we can write
\begin{eqnarray}\label{respomatdef1}
 F_{jl}&=&\sum_{m=1}^{M}\sum_{j=1}^{M}\mathcal{B}_{mj}e^{-i\kappa\theta_j\cdot z_m}e^{i\kappa\theta_l\cdot z_j}\nonumber\\
       &=&\left[e^{-i\kappa\theta_j\cdot z_1},e^{-i\kappa\theta_j\cdot z_2},\cdots,e^{-i\kappa\theta_j\cdot z_M}\right]\cdot\mathcal{B}\cdot
                                  \left[e^{i\kappa\theta_l\cdot z_1},e^{i\kappa\theta_l\cdot z_2},\cdots,e^{i\kappa\theta_l\cdot z_M}\right]^{T}
\end{eqnarray}
for all $j,l=1,\dots,N$. We can factorize the response matrix $F$ as below,
\begin{eqnarray}\label{respomatfact1}
 F= H^{*}\mathcal{B}H
\end{eqnarray}
where $H$ is a complex matrix of order $M\times N$ given by
\begin{eqnarray*}
 {H}_{jl}:=e^{i\kappa\theta_l\cdot z_j}
\end{eqnarray*}
for all $j,l=1,\dots,N$.
In order to determine the locations $z_{m}$, we consider a 3D-grid of sampling points $z\in\mathbb{R}^{3}$ in a region containing the scatterers $D_{1},D_{2},\dots,D_{M}$.
For each point $z$, we define the vector $\phi_{z}\in\mathbb{C}^{N}$ by
\begin{equation}\label{actotalphiws}
 \phi_{z}:=(e^{-i\kappa\theta_{1}\cdotp z}, e^{-i\kappa\theta_{2}\cdotp z},\dots, e^{-i\kappa\theta_{N}\cdotp z})^{T}.
\end{equation}
\noindent
\subsubsection{MUSIC characterisation of the response matrix}\label{MUSICchar-smallac-sdlp} Recall that MUSIC is essentially based on characterizing the range of the response matrix $F$ (signal space), forming projections
onto its null (noise) spaces, and computing its singular value decomposition.  In other words,
the MUSIC algorithm is based on the property that $\phi_{z}$ is in the range $\mathcal{R}(F)$ of $F$ if and only if $z$ is at one of locations of the scatterers. 
\textit{Precisely, let $\mathcal{P}$ be the projection onto the null space $\mathcal{N}(F^{*})=\mathcal{R}(F)^{\bot}$ of the adjoint matrix $F^{*}$ of $F$, then 
$$
z\in\{z_{1},z_{2},\dots,z_{M}\} \Longleftrightarrow \mathcal{P}\phi_{z}=0.
$$} 
This fact can be proved based on the non-singularity of the scattering matrix $\mathcal{B}$ in the factorization \eqref{respomatfact1} of $F$. Due to this, the standard linear algebraic argument yields that,
if $N\geq M$ and if the matrix $H$ has maximal rank $M$, then the ranges $\mathcal{R}(H^*)$ and $\mathcal{R}(F)$ coincide. 
\par For a sufficiently large number $N$ of incident and the observational directions by following the same lines as in \cite{K-G:2008,A-C-I:SIAM2008,DPC-SM:InvPrbm2012,G-S-T:JCAM2012}, 
the maximal rank property of $H$ can be justified. In this case MUSIC algorithm is applicable for our response matrix $F$. 
\par From the above discussion,  MUSIC characterization of the locations of the small scatterers in acoustic exterior Dirichlet problem can be written as the following 
and is valid if the number of the observational and incidental directions is sufficiently large.
\begin{theorem}\label{MUSIC-small}
Suppose $N\geq M$; then
\begin{eqnarray}\label{acsmalltheoremstatement1}
  z \in \{z_{1},\dots,z_{M}\}\Longleftrightarrow~\phi^{j}_{z}\in \mathcal{R}({H}^{*})=\mathcal{R}(F)\Longleftrightarrow \mathcal{P}\phi_{z}=0,
\end{eqnarray}
 where $\mathcal{P} : \mathbb{C}^{N}\rightarrow \mathcal{R}(F)^{\bot}= \mathcal{N}({F}^{*})$ is the orthogonal projection onto the null space $\mathcal{N}({F}^{*})$ of ${F}^{*}$. 
\end{theorem}
\ ~ \ \par
Let us point out that the orthogonal projection $\mathcal{P}\phi_{z}$ in the MUSIC algorithm does not contain any information about the shape and the orientation of the small scatterers. 
Yet, if the locations $z_m$ of the scatterers $D_m$ are found, approximately, via the observation of the pseudo norms of $\mathcal{P}\phi_{z}$, using the factorization 
\eqref{respomatfact1} of the response matrix $F$, then one can retrieve the capacitances $\bar{C}_{m}$ for $m=1,\dots,M$.

\subsection{Recovering the capacitances and estimating the sizes of the scatterers}\label{Recovering the capacitances-smallac-sdlp}
 Once we locate the scatterers from the given far-field patterns using the MUSIC algorithm, we can recover the capacitances $\bar{C}_{m}$ of $D_m$ 
from the factorization \eqref{respomatfact1} of $F\in\mathbb{C}^{N\times N}$ as hinted earlier in Section \ref{MUSICchar-smallac-sdlp}. 

Indeed,  we know that the matrix $H$ has maximal rank, see Theorem 4.1 of \cite{K-G:2008} for instance.  So, the matrix $HH^{*}\in\mathbb{C}^{M\times M}$ is invertible.  
Let us denote its inverse by $I_{H}$. Once we locate the scatterers through finding the locations $z_{1},z_{2},\dots,z_{M}$ by using the MUSIC algorithm for the given far-field patterns, 
we can recover $I_{H}$ and hence the matrix $\mathcal{B}\in C^{M\times M}$ given by $\mathcal{B}=I_{H}HFH^{*}I_{H}$, where $I_{H}H$ is the pseudo inverse of $H^{*}$.
As we know the structure of $\mathbf{B}\in\mathbb{C}^{M\times M}$, the inverse of $\mathcal{B}\in\mathbb{C}^{M\times M}$, we can recover the capacitances $\bar{C}_{1},\bar{C}_{2},\dots,\bar{C}_{M}$ 
of the small scatterers $D_1,D_2,\dots,D_M$ from the diagonal entries of $\mathbf{B}$. 
From these capacitances, we can estimate the size of the obstacles. Indeed, assume that $D_j$'s are balls of radius $\rho_j$, 
and center $0$ for simplicity, then we know that $\int_{\{y:|y|=\rho_j\}}\frac{dS_y}{|x-y|}=4\pi\rho_j$, for $\vert x \vert = \rho_j$, as observed 
in \cite[formula (5.12)]{M-M:MathNach2010}. Hence $\sigma_j(s)=\rho^{-1}_j$ and then $C_{j}=\int_{\partial D_j}\rho^{-1}_jds=4\pi \rho_j$ from which we can estimate 
the radius $\rho_j$. Other geometries, as cylinders, for which one can estimate exactly the size are shown in chapter 4 of \cite{MRAMAG-book1}. For general 
geometries, we can use the estimate (\ref{Capacitance-size-estimate}) to provide lower and upper estimates of the size of scatterers. For this,
 one should estimate the norms of the operators appearing in (\ref{Capacitance-size-estimate}) in terms of (only) the size of $\partial D_j$ and the Lipschitz 
smoothness character $L_j$ defined in the beginning of the introduction. This issue will be considered in a future work.

\subsection{Numerical results and discussions}\label{Numerical results and discussions-smallac-sdlp}
 
In order to illustrate relevant features of the method reviewed in the previous section, several computations have been performed, and the typical results acquired and presented here in.
To generate the far-field data we numerically solve the exterior acoustic Dirichlet problem (\ref{acimpoenetrable}-\ref{radiationc}) via a Galerkin method, 
see \cite{G-G:JCP2004,G_H:JANIAM2008}.
For our calculations, we considered $50$ incident and the observational directions obtained from the Gauss-Legendre polynomial. Let $d_{GL}$ stands for the degree of Gauss-Legendre polynomial.
Then these $2d_{GL}^2(=50)$ incident and the observational directions are obtained from the Gauss-Legendre polynomial of degree $d_{GL}(=5)$,
i.e. if we denote the zeros of the Gauss-Legendre polynomial of degree by $GL_{k}$, for $k=1,\dots,d_{GL}$ then the azimuth and the zenith angles $\theta$ and $\phi$ are given by
\begin{eqnarray*}\label{directions}
 \phi=\cos^{-1}(GL_{k}),~k=1,\dots,d_{GL}\\
 \theta=j*\frac{\pi}{d_{GL}},~j=0,1,\dots,2d_{GL}-1.
\end{eqnarray*}

\begin{figure}[htp]
\centering
\includegraphics[width=4.5cm,height =4.5cm]{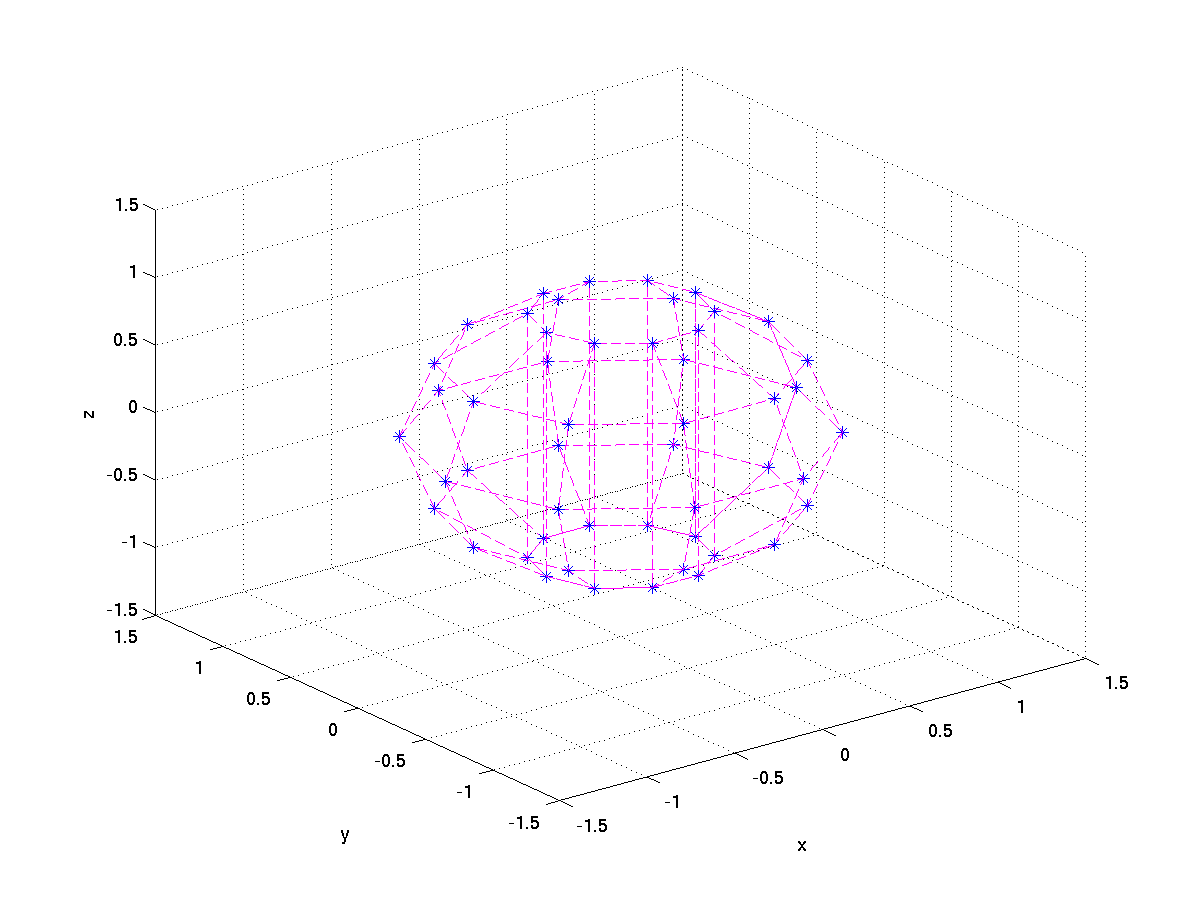}
\caption{ Incident and the observational directions.}\label{fig:1-acsmall}
\end{figure}

Combinations of these spherical coordinates will allow us to find the incident and the observational directions given by
$(\cos\theta\sin\phi, \sin\theta\sin\phi, \cos\phi)$. These directions are shown in Fig:\ref{fig:1-acsmall}.

\begin{figure}[]
\centering
\subfigure[]{
\includegraphics[width=5cm,height =4.5cm]{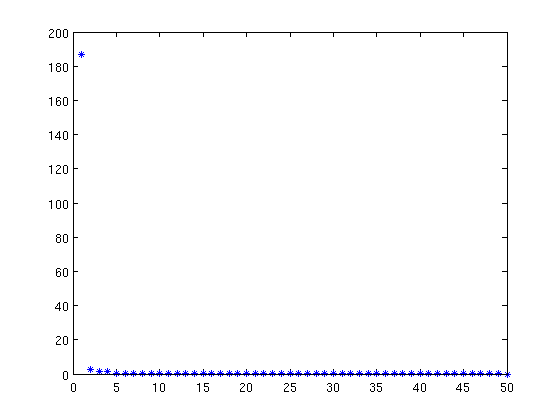}}\label{fig:2(a1)}
\subfigure[]{
\includegraphics[width=5cm,height =4.5cm]{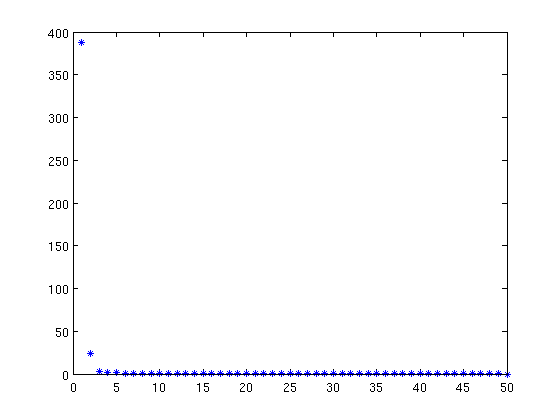}}\label{fig:2(c1)}
\subfigure[]{
\includegraphics[width=5cm,height =4.5cm]{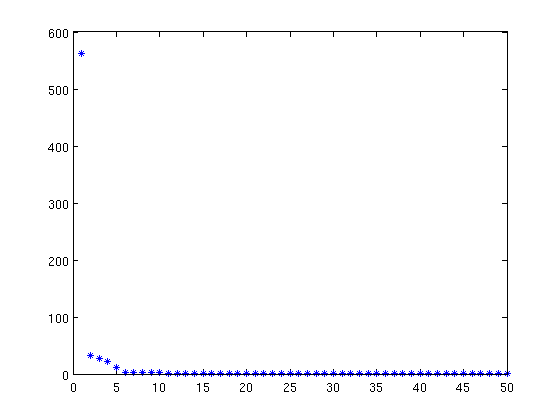}}\label{fig:2(b1)}\\

\subfigure[]{
\includegraphics[width=5cm,height=4.5cm]{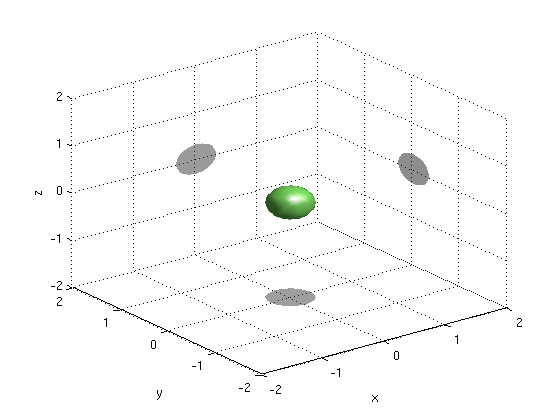}}\label{fig:2(g)}
\subfigure[]{
\includegraphics[width=5cm,height=4.5cm]{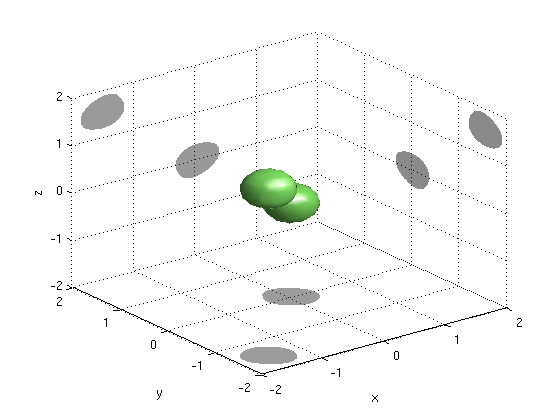}}\label{fig:2(i)}
\subfigure[]{
\includegraphics[width=5cm,height=4.5cm]{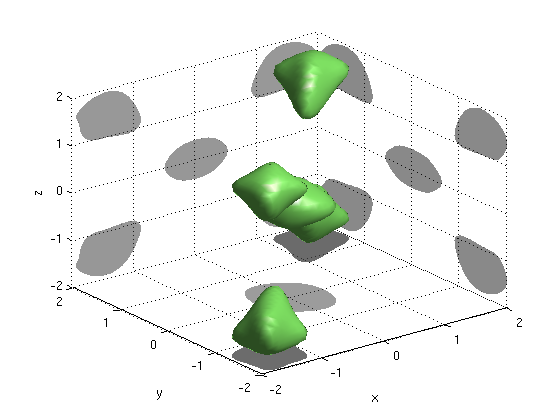}}\label{fig:2(h)}

\caption{ Reconstruction of the small scatterers with $30db$ SNR random noise}\label{fig:2}
\end{figure}

Due to the fact that MUSIC allows us to find the location of scatterers but not the size and shape, we presented the numerical results that are related to the balls of various size with various centers. 
Of course, we are able to locate the different type of scatterers as well. 
In Fig:\ref{fig:2} we have shown the results with $10\%$ random noise in the measured far-field 
performed on the balls with radius 0.5, and having the centers at A=(0,0,0), B=(1.5,1.5,1.5), C=(1.5,1.5,-1.5), D=(-1.5,-1.5,1.5) and at E=(-1.5,-1.5, -1.5).
Fig:\ref{fig:2}(a) and Fig:\ref{fig:2}(d) are results performed for the single small scatterer having center at origin A. Here the first one shows the distribution of the singular values and 
second one, iso-surface of the pseudonormal values,  shows the corresponding location of the scatterer. Fig:\ref{fig:2}(b) and Fig:\ref{fig:2}(c) are the distributions of the singular values 
while considering the  two and five small scatterers namely A,D and A,B,C,D,E respectively. Fig:\ref{fig:2}(e) and Fig:\ref{fig:2}(f) are the corresponding isosurface plots to locate the scatterers.

From these figures it can be seen that we have a good reconstruction of location of the scatterers. In general we can observe that, within our assumptions, MUSIC algorithm allows us to the locate 
the scatterers in a more finer way in the presence of less noise 
while as the noise increases due to the noise location of the scatterer will get disturbed.

To finish this section, let us mention that the reconstruction depends on the choice of the signal and noise subspaces of the multi scale response matrix. For small measurement noise 
(or) for higher SNR it is easy to choose these subspaces otherwise it become hard due to the smooth distribution of the singular values.
%
%
 
\section{Appendix: \it{Proof of Lemma \ref{Mazyawrkthm}}}
\begin{proofe}

We start by factorizing $\mathbf{B}$ as $\mathbf{B}=-(I+\mathbf{B}_{n}\mathbf{C})\mathbf{C}^{-1}$ where $\mathbf{C}:=Diag(\bar{C}_1,\bar{C}_2,\dots,\bar{C}_M)\in\mathbb{R}^{M\times\,M}$, $I$ is the identity matrix and 
$\mathbf{B}_{n}:=-\mathbf{C}^{-1}-\mathbf{B}$. Hence, the solvability of the system \eqref{compacfrm1}, depends on the existence of the inverse of $(I+\mathbf{B}_{n}\mathbf{C})$. 
We have $(I+\mathbf{B}_{n}\mathbf{C}):\mathbb{C}^{M}\rightarrow\mathbb{C}^{M}$, so it is enough to prove the injectivity in order to prove its invetibility. For this purpose,  let $X,Y$ are vectors in $\mathbb{C}^{M}$ and
 consider the system
\begin{eqnarray}\label{systemsolve1-small-ac}(I+\mathbf{B}_{n}\mathbf{C})X&=&Y.\end{eqnarray}
Let ${(\cdot)}^{real}$ and ${(\cdot)}^{img}$ denotes the real and the imaginary parts of the corresponding complex number/vector/matrix. Now, the following can be written from \eqref{systemsolve1-small-ac};
\begin{eqnarray}
 (I+\mathbf{B}^{real}_{n}\mathbf{C})X^{real}-\mathbf{B}^{img}_{n}\mathbf{C}X^{img}&=&Y^{real}\label{systemsolve1-small-sub1-ac},\\
 (I+\mathbf{B}^{real}_{n}\mathbf{C})X^{img}+\mathbf{B}^{img}_{n}\mathbf{C}X^{real}&=&Y^{img}\label{systemsolve1-small-sub2-ac},
\end{eqnarray}
which leads to
\begin{eqnarray}
 \langle\,(I+\mathbf{B}^{real}_{n}\mathbf{C})X^{real},\mathbf{C}X^{real}\rangle\,-\langle\,\mathbf{B}^{img}_{n}\mathbf{C}X^{img},\mathbf{C}X^{real}\rangle&=&\langle\,Y^{real},\mathbf{C}X^{real}\rangle\label{systemsolve1-small-sub1-ac1},\\
 \langle\,(I+\mathbf{B}^{real}_{n}\mathbf{C})X^{img},\mathbf{C}X^{img}\rangle\,+\langle\,\mathbf{B}^{img}_{n}\mathbf{C}X^{real},\mathbf{C}X^{img}\rangle&=&\langle\,Y^{img},\mathbf{C}X^{img}\rangle\label{systemsolve1-small-sub2-ac1}.
\end{eqnarray}
By summing up \eqref{systemsolve1-small-sub1-ac1} and \eqref{systemsolve1-small-sub2-ac1} will give
\begin{equation}
\begin{split}
\langle\,X^{real},\mathbf{C}X^{real}\rangle\,+\langle\,\mathbf{B}^{real}_{n}\mathbf{C}X^{real},\mathbf{C}X^{real}\rangle\,+\langle\,X^{img},\mathbf{C}X^{img}\rangle\,+\langle\,\mathbf{B}^{real}_{n}\mathbf{C}X^{img},\mathbf{C}X^{img}\rangle\,\\
=\langle\,Y^{real},\mathbf{C}X^{real}\rangle+\langle\,Y^{img},\mathbf{C}X^{img}\rangle.\label{systemsolve1-small-sub1-2-ac1}
\end{split}
\end{equation}
Indeed,
\[ \langle\,\mathbf{B}^{img}_{n}\mathbf{C}X^{img},\mathbf{C}X^{real}\rangle=\langle\,\mathbf{C}X^{img},\mathbf{B}^{img^{*}}_{n}\mathbf{C}X^{real}\rangle
=\langle\,\mathbf{C}X^{img},\mathbf{B}^{img}_{n}\mathbf{C}X^{real}\rangle=\langle\,\mathbf{B}^{img}_{n}\mathbf{C}X^{real},\mathbf{C}X^{img}\rangle.\]
We can observe that, the right-hand side in \eqref{systemsolve1-small-sub1-2-ac1} does not exceed
\begin{equation}\label{{systemsolve1-small-sub1-2-ac2}}
\begin{split}
\langle\,X^{real},\mathbf{C}X^{real}\rangle^{1\slash\,2}\langle\,Y^{real},\mathbf{C}Y^{real}\rangle^{1\slash\,2}+\langle\,X^{img},\mathbf{C}X^{img}\rangle^{1\slash\,2}\langle\,Y^{img},\mathbf{C}Y^{img}\rangle^{1\slash\,2}\\
\leq2\langle\,X^{|\cdot|},\mathbf{C}X^{|\cdot|}\rangle^{1\slash\,2}\langle\,Y^{|\cdot|},\mathbf{C}Y^{|\cdot|}\rangle^{1\slash\,2}.
\end{split}
\end{equation}
Here $W^{|\cdot|}_m:={[W^{{real}^2}_m+W^{{img}^2}_m]}^{1\slash2}= |W_m|$, 
for $W=X,Y$ and $m=1,\dots,M$.
Consider the second term in the left-hand side of \eqref{systemsolve1-small-sub1-2-ac1}. Using the mean value theorem for harmonic functions we deduce
\begin{eqnarray*}
\langle\,\mathbf{B}^{real}_{n}\mathbf{C}X^{real},\mathbf{C}X^{real}\rangle&=&\sum_{\substack{1\leq\,j,m\leq\,M\\j\neq\,m}}\Phi^{real}(z_m,z_j)\bar{C}_j\bar{C}_mX^{real}_{j}X^{real}_{m}\\
&\geq&t\sum_{\substack{1\leq\,j,m\leq\,M\\j\neq\,m}}\frac{\bar{C}_j\bar{C}_mX^{real}_{j}X^{real}_{m}}{|B^{(j)} \|B^{(m)}|}\int_{B^{(j)}}\int_{B^{(m)}}\Phi_0(x,y)~dx~dy,
\end{eqnarray*}
Similarly, if we consider the fourth term in the left-hand side of \eqref{systemsolve1-small-sub1-2-ac1}, we deduce
\begin{eqnarray*}
\langle\,\mathbf{B}^{real}_{n}\mathbf{C}X^{img},\mathbf{C}X^{img}\rangle&=&\sum_{\substack{1\leq\,j,m\leq\,M\\j\neq\,m}}\Phi^{real}(z_m,z_j)\bar{C}_j\bar{C}_mX^{img}_{j}X^{img}_{m}\\
&\geq&t\sum_{\substack{1\leq\,j,m\leq\,M\\j\neq\,m}}\frac{\bar{C}_j\bar{C}_mX^{img}_{j}X^{img}_{m}}{|B^{(j)} \|B^{(m)}|}\int_{B^{(j)}}\int_{B^{(m)}}\Phi_0(x,y)~dx~dy,
\end{eqnarray*}
where $t:=\min\limits_{j\neq\,m,1\leq\,j,m\leq\,M}\cos(\kappa|z_m-z_j|)$, assumed to be non negative, $\Phi_0(x,y)=1\slash(4\pi|x-y|)$ is defined in \eqref{definition-ac-small-fundamentalkappa}
and $B^{(m)}:=\{x:|x-z_m|<d\slash2\},m=1,\dots,M$, are non-overlapping balls of radius $d\slash2$ with centers at $z_m$, and $|B^{(m)}|=\pi\,d^3\slash6$ are the volumes of the balls. 
Also, we use the notation $B_d$ to denote the balls of radius $d\slash2$   with the center at the origin.

Let $\Omega$ be a large ball with radius $R$. Also let $\Omega_s\subset\Omega$ be a ball with fixed radius $r(\leq\,R)$, 
which consists of all our small obstacles $D_m$ and also the balls $B^{(m)}$, for $m=1,\dots,M$.

Let $\Upsilon^{real}(x)$ and $\Upsilon^{img}(x)$ be piecewise constant functions defined on $\mathbb{R}^3$ as
\begin{equation}\label{def-upsilon}
 \Upsilon^{real(img)}(x)=\begin{cases}
\begin{array}{ccc}
              X^{real(img)}_m\bar{C}_m &\mbox{in }  B^{(m)},&m=1,\dots,M,\\
              0 &\mbox{otherwise}.
\end{array}
             \end{cases}
\end{equation}
Then 
\begin{equation}\label{systemsolve1-small-sub1-2-ac1-basedonupsilon-real}
\begin{split}
\langle\,\mathbf{B}^{real}_{n}\mathbf{C}X^{real},\mathbf{C}X^{real}\rangle\geq\frac{36t}{\pi^2\,d^6}&\left(\int_\Omega\int_\Omega\Phi_0(x,y)\Upsilon^{real}(x)\Upsilon^{real}(y)~dx~dy\right.\\
&\left.-\sum_{m=1}^{M}\bar{C}_m^2X^{real^2}_{m}\int_{B^{(m)}}\int_{B^{(m)}}\Phi_0(x,y)~dx~dy\right),
\end{split}
\end{equation}
\begin{equation}\label{systemsolve1-small-sub1-2-ac1-basedonupsilon-img}
\begin{split}
\langle\,\mathbf{B}^{real}_{n}\mathbf{C}X^{img},\mathbf{C}X^{img}\rangle\geq\frac{36t}{\pi^2\,d^6}&\left(\int_\Omega\int_\Omega\Phi_0(x,y)\Upsilon^{img}(x)\Upsilon^{img}(y)~dx~dy\right.\\
&\left.-\sum_{m=1}^{M}\bar{C}_m^2X^{img^2}_{m}\int_{B^{(m)}}\int_{B^{(m)}}\Phi_0(x,y)~dx~dy\right).
\end{split}
\end{equation}
Applying the mean value theorem to the harmonic function $\frac{1}{4\pi\vert x-y \vert}$, as done in \cite[p:109-110]{M-M:MathNach2010}, we have the following estimate
\begin{equation}
 \begin{split}\label{systemsolve1-small-sub1-2-ac1-basedonupsilon-2ndpart}
  \int_{B^{(m)}}\int_{B^{(m)}}\Phi_0(x,y)~dx~dy\,=\,\frac{1}{4\pi}\int_{B_d}\int_{B_d}\frac{1}{|x-y|}~dx~dy \,\leq\, \frac{\pi\,d^5}{60}.
 \end{split}
\end{equation}
Consider the first term in the right-hand side of \eqref{systemsolve1-small-sub1-2-ac1-basedonupsilon-real}, denote it by $A_R^{real}$, then by Green's theorem
\begin{eqnarray}\label{systemsolve1-small-sub1-2-ac1-basedonupsilon-1stpart}
A_R^{real}&:=&\int_\Omega\int_\Omega\Phi_0(x,y)\Upsilon^{real}(x)\Upsilon^{real}(y)~dx~dy\\
&=&\underbrace{\int_\Omega\left|\nabla_x\int_\Omega\Phi_0(x,y)\Upsilon^{real}(y)~dy\right|^2~dx}_{=:B_R^{real}\geq0}-
\underbrace{\int_{\partial\Omega}\left(\frac{\partial}{\partial\nu_x}\int_\Omega\Phi_0(x,y)\Upsilon^{real}(y)~dy\right)~\left(\int_\Omega\Phi_0(x,y)\Upsilon^{real}(y)~dy\right)dS_x}_{=:C_R^{real}}.\nonumber
\end{eqnarray}
We have
\begin{equation}\label{systemsolve1-small-sub1-2-ac1-basedonupsilon-1stpart-1}
 \begin{split}
  C_R^{real}
&=\int_{\partial\Omega}\left(\int_\Omega\frac{\partial}{\partial\nu_x}\Phi_0(x,y)\Upsilon^{real}(y)~dy\right)~\left(\int_\Omega\Phi_0(x,y)\Upsilon^{real}(y)~dy\right)dS_x\\
&=\int_{\partial\Omega}\left(\int_{\Omega_s}\frac{\partial}{\partial\nu_x}\Phi_0(x,y)\Upsilon^{real}(y)~dy\right)~\left(\int_{\Omega_s}\Phi_0(x,y)\Upsilon^{real}(y)~dy\right)dS_x\\
&=\int_{\partial\Omega}\left(\int_{\Omega_s}\frac{-(x-y)}{4\pi|x-y|^3}\Upsilon^{real}(y)~dy\right)~\left(\int_{\Omega_s}\frac{1}{4\pi|x-y|}\Upsilon^{real}(y)~dy\right)dS_x,
 \end{split}
\end{equation}
which gives the following estimate;
\begin{equation}\label{systemsolve1-small-sub1-2-ac1-basedonupsilon-1stpart-2}
 \begin{split}
|C_R^{real}|&\leq\frac{1}{16\pi^2}\int_{\partial\Omega}\frac{1}{|R-r|^3}\left(\int_{\Omega_s}|\Upsilon^{real}(y)|~dy\right)^2dS_x\\
&\leq\frac{1}{16\pi^2}\frac{1}{(R-r)^3}\int_{\partial\Omega}|\Omega_s|~ \|\Upsilon^{real} \|^2_{L^2({\Omega_s})}dS_x\\
&=\frac{r^3d^3}{72(R-r)^3}\sum_{m=1}^{M}X^{real^2}_m\bar{C}^2_m~|\partial\Omega|\\
&=\frac{{\pi}R^2r^3d^3}{18(R-r)^3}\sum_{m=1}^{M}X^{real^2}_m\bar{C}^2_m.\\
\end{split}
\end{equation}
Substitution of \eqref{systemsolve1-small-sub1-2-ac1-basedonupsilon-1stpart-2} in \eqref{systemsolve1-small-sub1-2-ac1-basedonupsilon-1stpart} gives
\begin{eqnarray}\label{systemsolve1-small-sub1-2-ac1-basedonupsilon-1stpart-3}
\int_\Omega\int_\Omega\Phi_0(x,y)\Upsilon^{real}(x)\Upsilon^{real}(y)~dx~dy
&\geq&\int_\Omega\left|\nabla_x\int_\Omega\Phi_0(x,y)\Upsilon^{real}(y)~dy\right|^2~dx-\frac{{\pi}R^2r^3d^3}{18(R-r)^3}\sum_{m=1}^{M}X^{real^2}_m\bar{C}^2_m.
\nonumber\\
\end{eqnarray}
By considering the first term in the right-hand side of \eqref{systemsolve1-small-sub1-2-ac1-basedonupsilon-img}, and following the  same procedure as mentioned in 
\eqref{systemsolve1-small-sub1-2-ac1-basedonupsilon-1stpart}, \eqref{systemsolve1-small-sub1-2-ac1-basedonupsilon-1stpart-1} and \eqref{systemsolve1-small-sub1-2-ac1-basedonupsilon-1stpart-2}, we obtain
\begin{eqnarray}\label{systemsolve1-small-sub1-2-ac1-basedonupsilon-1stpart-4}
\int_\Omega\int_\Omega\Phi_0(x,y)\Upsilon^{img}(x)\Upsilon^{img}(y)~dx~dy
&\geq&\int_\Omega\left|\nabla_x\int_\Omega\Phi_0(x,y)\Upsilon^{img}(y)~dy\right|^2~dx-\frac{{\pi}R^2r^3d^3}{18(R-r)^3}\sum_{m=1}^{M}X^{img^2}_m\bar{C}^2_m.
\nonumber\\
\end{eqnarray}
Under our assumption $t\geq0$,  \eqref{systemsolve1-small-sub1-2-ac1-basedonupsilon-real}, 
\eqref{systemsolve1-small-sub1-2-ac1-basedonupsilon-img}, \eqref{systemsolve1-small-sub1-2-ac1-basedonupsilon-2ndpart}, \eqref{systemsolve1-small-sub1-2-ac1-basedonupsilon-1stpart-3} and \eqref{systemsolve1-small-sub1-2-ac1-basedonupsilon-1stpart-4} lead to
\begin{eqnarray}
\langle\,\mathbf{B}^{real}_{n}\mathbf{C}X^{real},\mathbf{C}X^{real}\rangle\geq\frac{36t}{\pi^2\,d^6}
&\left(\int_\Omega\left|\nabla_x\int_\Omega\Phi_0(x,y)\Upsilon^{real}(y)~dy\right|^2~dx-\left[\frac{{\pi}R^2r^3d^3}{18(R-r)^3}+\frac{\pi\,d^5}{60}\right]\sum_{m=1}^{M}X^{real^2}_m\bar{C}^2_m\right),
\nonumber\\
\label{systemsolve1-small-sub1-2-ac1-basedonupsilon-real-1}\\
\langle\,\mathbf{B}^{real}_{n}\mathbf{C}X^{img},\mathbf{C}X^{img}\rangle\geq\frac{36t}{\pi^2\,d^6}
&\left(\int_\Omega\left|\nabla_x\int_\Omega\Phi_0(x,y)\Upsilon^{img}(y)~dy\right|^2~dx-\left[\frac{{\pi}R^2r^3d^3}{18(R-r)^3}+\frac{\pi\,d^5}{60}\right]\sum_{m=1}^{M}X^{img^2}_m\bar{C}^2_m\right).
\nonumber
\\\label{systemsolve1-small-sub1-2-ac1-basedonupsilon-img-1}
\end{eqnarray}
Then \eqref{systemsolve1-small-sub1-2-ac1}, \eqref{systemsolve1-small-sub1-2-ac1-basedonupsilon-real-1} and \eqref{systemsolve1-small-sub1-2-ac1-basedonupsilon-img-1} imply
\begin{equation}\label{fnlinvert-small-ac-1}
\begin{split}
 \left(1-\frac{36t}{\pi^2\,d^6}\right.&\left.\left[\frac{{\pi}R^2r^3d^3}{18(R-r)^3}+\frac{\pi\,d^5}{60}\right]\max\limits_{1\leq m \leq M}\bar{C}_m\right)\sum_{m=1}^{M}|X_m|^{2}\bar{C}_m\\
&=\left(1-\frac{36t}{\pi^2\,d^6}\left[\frac{{\pi}R^2r^3d^3}{18(R-r)^3}+\frac{\pi\,d^5}{60}\right]\max\limits_{1\leq m \leq M}\bar{C}_m\right)\left(\sum_{m=1}^{M}X^{real^2}_m\bar{C}_m+\sum_{m=1}^{M}X^{img^2}_m\bar{C}_m\right)\\
&\leq\left(\sum_{m=1}^{M}X^{real^2}_m\bar{C}_m\right)^{1/2}\left(\sum_{m=1}^{M}Y^{real^2}_m\bar{C}_m\right)^{1/2}+\left(\sum_{m=1}^{M}X^{img^2}_m\bar{C}_m\right)^{1/2}\left(\sum_{m=1}^{M}Y^{img^2}_m\bar{C}_m\right)^{1/2}\\
&\leq2\left(\sum_{m=1}^{M}|X_m|^2\bar{C}_m\right)^{1/2}\left(\sum_{m=1}^{M}|Y_m|^2\bar{C}_m\right)^{1/2}.
\end{split}
\end{equation}
As we have $R$ arbitrary, by tending $R$ to $\infty$, we can write 
\eqref{fnlinvert-small-ac-1}
as 
\begin{equation}\label{fnlinvert-small-ac-2}
\begin{split}
 \left(1-\frac{3t}{5\pi\,d}\right.\left.\max\limits_{1\leq m \leq M}\bar{C}_m\right)\sum_{m=1}^{M}|X_m|^{2}\bar{C}_m
\leq2\left(\sum_{m=1}^{M}|X_m|^2\bar{C}_m\right)^{1/2}\left(\sum_{m=1}^{M}|Y_m|^2\bar{C}_m\right)^{1/2},
\end{split}
\end{equation}
which  yields
\begin{equation}\label{fnlinvert-small-ac-2f}
\begin{split}
\sum_{m=1}^{M}|X_m|^{2}\bar{C}_m
\leq4 \left(1-\frac{3t}{5\pi\,d}\right.&\left.\max\limits_{1\leq m \leq M}\bar{C}_m\right)^{-2}\sum_{m=1}^{M}|Y_m|^2\bar{C}_m.
\end{split}
\end{equation}
Thus, if $\max\limits_{1\leq\,m\leq\,M}\bar{C}_m<\frac{5\pi}{3}d$ and $t \geq 0$, then the matrix $\mathbf{B}$ 
in algebraic system \eqref{compacfrm1} is invertible and the estimate \eqref{fnlinvert-small-ac-2} and so \eqref{mazya-fnlinvert-small-ac-2} holds.
\end{proofe}
\bibliographystyle{abbrv}

\end{document}